\documentclass[12pt]{amsart}
\usepackage{ amsmath, amsthm, amsfonts, amssymb, color}
 \usepackage{mathrsfs}
\usepackage{amsfonts, amsmath}
 \usepackage{amsmath,amstext,amsthm,amssymb,amsxtra}
 \usepackage{txfonts} 
 \usepackage[colorlinks, citecolor=blue,pagebackref,hypertexnames=false]{hyperref}
 \allowdisplaybreaks
 \usepackage{pgf,tikz}

\begin{document}

\hfuzz=6pt

\widowpenalty=10000

\newtheorem{cl}{Claim}
\newtheorem{theorem}{Theorem}[section]
\newtheorem{proposition}[theorem]{Proposition}
\newtheorem{coro}[theorem]{Corollary}
\newtheorem{lemma}[theorem]{Lemma}
\newtheorem{definition}[theorem]{Definition}
\newtheorem{assum}{Assumption}[section]
\newtheorem{example}[theorem]{Example}
\newtheorem{remark}[theorem]{Remark}
\renewcommand{\theequation}
{\thesection.\arabic{equation}}

\def\SL{\sqrt L}

\newcommand{\mar}[1]{{\marginpar{\sffamily{\scriptsize
        #1}}}}
\newcommand{\li}[1]{{\mar{LY:#1}}}
\newcommand{\el}[1]{{\mar{EM:#1}}}
\newcommand{\as}[1]{{\mar{AS:#1}}}

\newcommand\RR{\mathbb{R}}
\newcommand\CC{\mathbb{C}}
\newcommand\NN{\mathbb{N}}
\newcommand\ZZ{\mathbb{Z}}

\renewcommand\Re{\operatorname{Re}}
\renewcommand\Im{\operatorname{Im}}

\newcommand{\mc}{\mathcal}
\newcommand\D{\mathcal{D}}

\newcommand{\la}{\lambda}
\def \l {\lambda}
\newcommand{\eps}{\varepsilon}
\newcommand{\pl}{\partial}
\newcommand{\supp}{{\rm supp}{\hspace{.05cm}}}
\newcommand{\x}{\times}

\newcommand\wrt{\,{\rm d}}

\title[Restriction estimates and sharp spectral multipliers  ]{%
Restriction estimates, sharp spectral multipliers and endpoint estimates for
 Bochner-Riesz means}
\author{Peng Chen}
\author{El Maati Ouhabaz}
\author{Adam Sikora}
\author{Lixin Yan}
\address{Peng Chen, Department of Mathematics, Sun Yat-sen (Zhongshan)
University, Guangzhou, 510275, P.R. China}
\email{achenpeng1981@163.com}
\address{El Maati Ouhabaz,  Institut de Math\'ematiques de Bordeaux,  Universit\'e Bordeaux 1, UMR 5251,
351, Cours de la Lib\'eration 33405 Talence, France}
\email{Elmaati.Ouhabaz@math.u-bordeaux1.fr}
\address {Adam Sikora, Department of Mathematics, Macquarie University, NSW 2109, Australia}
\email{sikora@maths.mq.edu.au}
\address{Lixin Yan, Department of Mathematics, Sun Yat-sen (Zhongshan) University,
Guangzhou, 510275, P.R. China}
\email{mcsylx@mail.sysu.edu.cn}
\date{\today}

\begin{abstract}
We consider abstract non-negative  self-adjoint operators on $L^2(X)$ which satisfy the finite speed propagation
property for the corresponding wave equation. For such operators we introduce a restriction
type condition which in the case of  the standard Laplace operator is equivalent to $(p,2)$
restriction estimate of Stein and Tomas. Next we show that in the considered  abstract setting
our restriction type condition implies  sharp spectral multipliers and
endpoint estimates for the Bochner-Riesz summability.  We also observe that
this restriction estimate holds for operators satisfying
dispersive or Strichartz estimates. We obtain new spectral multiplier results for
several second order differential operators and recover some known results.
Our examples include    Schr\"odinger operators with inverse square potentials on
$\RR^n$, the harmonic oscillator, elliptic operators on compact manifolds  and Schr\"odinger operators  on asymptotically conic manifolds.

\end{abstract}

\maketitle

 \tableofcontents

\section{Introduction}
\setcounter{equation}{0}

A celebrated   theorem of  H\"ormander \cite{Ho4} states  that for a given bounded function $F: [0, \infty) \to \CC$,
the operator $F(-\Delta)$, initially defined by Fourier analysis  on $L^2(\RR^n)$,  extends to a bounded operator on
$L^p(\RR^n)$ for all $p \in (1, \infty)$ provided the function satisfies
\begin{equation}\label{I0}
\sup_{t > 0} \| \eta(\cdot) F(t \cdot) \|_{W^{s, 2}} < \infty
\end{equation}
for some $s >  \frac{n}{2}$. Here $\eta \in C_c^\infty(0,\infty)$ is a non-trivial auxiliary function.
This result is a sharp version of the well known Mikhlin's Fourier multiplier theorem \cite{Mi}.
These results have led to a fruitful research activity on spectral multipliers and
new perspectives in harmonic analysis. The H\"ormander-Mikhlin theorem has been extended by several
authors to other operators  than the Laplacian and  settings  that go beyond the Euclidean case.  The bibliography is
so broad that it is impossible  to provide complete
list here. We refer  the reader  to \cite{A, CS,
C1, C2, C3,  ChS, CowS, DOS, DSY, F1, F2, F3, GHS, H,  Ho2, Ho3, MM, MS, Ou, Se,  SS, S1, Sog1, Sog4,
St1, St2, T2, Th1} and the references therein.

Suppose that $X$ is a measure space and that $L$ is a non-negative
self-adjoint  operator on~$L^2(X)$. Such an operator admits a spectral
resolution  $E_L(\lambda)$  and for  any  bounded  Borel function $F\colon [0, \infty)
\to {\mathbb C}$, one can define the operator $F(L)$
\begin{equation}\label{eq1.1}
F(L)=\int_0^{\infty}F(\lambda) \wrt E_L(\lambda).
\end{equation}
By the spectral theorem, $F(L)$ is well defined and bounded on~$L^2(X)$. Spectral multiplier
theorems give sufficient conditions on~$F$ under which the operator~$F(L)$ extends to a
bounded operator on~$L^p(X)$ for some range of~$p$.

Most of the references mentioned before  deal with the case of sub-Laplacians   on some Lie groups.
The papers \cite{DOS, DSY} deal with a rather general situation
 where $(X, d, \mu)$  is metric measure space of homogeneous type  (or even a domain of such space).
One of the results there  is a spectral multiplier theorem under the sole assumption that the heat kernel
of the operator has a Gaussian upper  bound. The condition there is however stronger than (\ref{I0}) in
 the sense that the norm in $W^{s, 2}$ is replaced by the norm of $W^{s, \infty}$ where $s$ is any constant
 larger than half of homogeneous dimension. Under Plancherel type estimates, one obtains sharp results
 with condition (\ref{I0}).  Condition (\ref{I0}) with  norm $W^{s,2}$ is  better than the
 corresponding one  with norm $W^{s, \infty}$.  This can be seen from Bochner-Riesz summability  which we discuss now.

The theory of spectral multipliers is  related to  and motivated by  the study of convergence
of  Bochner-Riesz   means  of self-adjoint  operators.  Given a non-negative self-adjoint  operator $L$ and  set
 \begin{equation}\label{eq1.2}
 S^{\delta}_R(\lambda)=
      \left\{
       \begin{array}{cl}
       \left(1-\frac{\lambda}{R^2}\right)^{\delta}  &\mbox{for}\;\; \lambda \le R^2 \\ [8pt]
       0  &\mbox{for}\;\; \lambda > R^2. \\
       \end{array}
      \right.
   \end{equation}
We then define the operator $S^{\delta}_R(L)$ using (\ref{eq1.1}).
We call $S^{\delta}_R(L)$ the Riesz or the Bochner-Riesz means of
order $\delta$. The basic question in the theory of Bochner-Riesz means
is  to establish  the critical exponent for the uniform continuity
with respect to $R$ and convergence of the Riesz means on $L^p$
spaces for various  $p$ with $1\le p \le \infty$.

For ${\delta}=0$, this is the spectral projector $E_{\SL}[0, R]$, while for ${\delta}>0$, $S^{\delta}_R(L)$
can be seen as a smoothed version of this spectral projector. Bochner-Riesz summability describes the range of ${\delta}$ for which the
above operators are bounded on $L^p$,  uniformly in $R$.
If one proves a  spectral multiplier result which states that $F(L)$ is bounded on $L^p(X)$
for all $p \in (1, \infty)$ whenever $F$ satisfies (\ref{I0}), then
the Bochner-Riesz  mean  $S^\delta_1(L)$ is  bounded on all  $L^p$ spaces provided $ \delta > \frac{n-1}{2}$.
In the case where $L$ is the Euclidean Laplacian  and $ \delta > \frac{n-1}{2}$,
 the kernel of $S^{\delta}_R(L)$ is $L^1$ and hence Bochner-Riesz means  are bounded on $L^p$  for $1 \le p \le \infty$.
In this setting more is known. Indeed for
  $ \delta > \max\{ n  | \frac{1}{2} - \frac{1}{p} | - \frac{1}{2}, 0 \} $,  it was known for a
  long time that as a consequence of restriction estimates for  the Fourier transform,
  Bochner-Riesz summability holds on $L^p$ for all $p \le \frac{2n + 2}{n+3}$,  and
  by duality for $p \ge \frac{2n +2}{n-1}$. See  Stein \cite{St2}, p. 420.    This was
  extended by Lee \cite{Lee} to the case $p < \frac{2n +4}{n+4}$ (or $p > \frac{2n+4}{n}$) and recent  improvements
  are proved  by  Bourgain and Guth \cite{BoGu}.  The
  question whether  Bochner-Riesz summability holds  on $L^p(\RR^n)$ for all $p$
  and all $\delta > \max\{ n | \frac{1}{p} - \frac{1}{2} | - \frac{1}{2}, 0 \}$  is a long-standing  open problem
   (except for  $n =  2$, see Carleson and  Sj\"olin  \cite{CS} and
   H\"ormander \cite{Ho22}).
For all this, see Stein \cite{St2}, p. 420  and the review paper of Tao \cite{T1}. The later contains many
other  information and relation of the Bochner-Riesz problem to other open problems in harmonic analysis.

If $L$ is a second order elliptic operator on a compact Riemannian manifold $M$ with dimension $n$,
then the Bochner-Riesz means $\sum_{\lambda_j \le R} (1- \lambda_j/R)^\delta <\cdot, e_j> e_j$ are uniformly bounded
 on $L^p(M)$ provided  $p \le \frac{2n + 2}{n+3}$ or  by duality  for $p \ge \frac{2n +2}{n-1}$ for
 $ \delta > \max\{ n  | \frac{1}{2} - \frac{1}{p} | - \frac{1}{2}, 0 \} $.
Here $\lambda_0 \le \lambda_j \le \lambda_{j+1}\le \ldots$ and $e_j$ are the corresponding eigenvalues and normalized
$L^2$ eigenvectors, respectively.  See Sogge \cite{Sog1}.

The theory of Fourier multipliers and Bochner-Riesz analysis in the setting of the standard Laplace operator
 on $\RR^n$ is   related to the so-called \emph{sphere restriction problem} for the Fourier transform:
find the pairs $(p, q)$  for which $R_\la\in \mc{L}(L^p(\RR^n),L^q({\mathbb S}^{n-1}))$ where $R_\la$ is defined by
\[R_\la f(\omega)=\hat{f}(\la\omega), \, \, \omega\in {\mathbb S}^{n-1}, \la >0.\]
See for example \cite{F1, F3, St1, St2,  T1}.
For $q=2$ the full description of possible range of $p$ is due  to Stein and Tomas.
The theorem of Tomas \cite{Tom}, extended  by Stein to the endpoint, states that $(p,2)$ restriction estimates hold
if and only if   $1 \le p \le {2(n+1)/(n+3)}$.   The case $q\neq 2$ is not relevant to our discussion so
 we refer the interested reader to  Tao \cite{T1} on the subject.

Note that on $\RR^n$, the Schwartz kernel of the spectral measure $dE_{\sqrt{-\Delta}}(\la)$ of
$\sqrt{-\Delta}$ is given by
\[dE_{\sqrt{-\Delta}}(\la;z,z')=\frac{\la^{n-1}}{(2\pi)^n}\int_{S^{n-1}}e^{i(z-z')\cdot \la\omega}d\omega,
 \quad z, z' \in \RR^n, \]
therefore $dE_{\sqrt{-\Delta}}(\la)=\frac{\la^{n-1}}{(2\pi)^n}R_\la^*R_\la$ and the restriction
 theorem for $q=2$ is equivalent to
\begin{equation}\label{I1}
\|dE_{\sqrt{-\Delta}}(\lambda) \|_{p\to p'} \le C \lambda^{n(1/p-1/p')-1}
\end{equation}
for all $p \in [1, \frac{2n +2}{n+3}]$. In the sequel, we  refer to  (\ref{I1}) as  $(p,2)$ restriction estimate of Stein-Tomas.

In this paper  we follow  the  line of  research described  above. We  deal  with  the problem
of sharp spectral multipliers and Bochner-Riesz summability  for other operators than the Euclidean Laplacian
and elliptic operators on compact manifolds. Our aim is to build a theory which applies in a rather general
setting of self-adjoint operators on spaces of homogeneous type (i.e., metric measure  spaces which satisfy
 the volume doubling property).
Our approach allows  us to prove
sharp multiplier results and Bochner-Riesz summability in new settings
and also unifies  several previously known results.
In order to do so we introduce a restriction type estimate which in the case of the Laplacian on $\RR^n$ turns to be  equivalent
to  the $(p,2)$ restriction estimate of Stein-Tomas.

Our setting will be the following.
We consider a non-negative self-adjoint operator  $L$ on $L^2(X)$ where $(X, d, \mu)$ is a metric measure space which
satisfies the volume doubling condition
$$ V(x, \lambda r) \le C \lambda^n V(x,r)\,\;  \forall x \in X,\, \lambda \ge 1, \, r > 0,$$
where $C$ and $n$ are positive constants and $V(x,r)$ denotes the volume of the open ball $B(x,r)$  of centre $x$
and radius $r$.    We assume that $L$ satisfies the finite speed propagation property for the corresponding wave equation.
We introduce the  condition that for any $R>0$ and all Borel functions $F$ supported   in $ [0,R],$
$$
\big\|F(\SL)P_{B(x, r)} \big\|_{p\to s} \leq CV(x,
r)^{{1\over s}-{1\over p}} \big( Rr \big)^{n({1\over p}-{1\over
s})}\big\| F (R \cdot) \big\|_{q}
\leqno{\rm (ST^q_{p, s})}
$$
for all $x\in X$ and all $r\geq 1/R$.

We will see that if the volume is polynomial, i.e. $V(x,r) \sim r^n$, then ${\rm (ST^2_{p, 2})}$ is equivalent to $(p,2)$
restriction estimate of Stein-Tomas. For this reason, we call
${\rm (ST^q_{p, s})}$ a {\it Stein-Tomas restriction type condition}. One of our main results  on sharp spectral multipliers
 can be stated as follows.

\medskip
\noindent{\bf Theorem A.} {\it
Assume that  $X$ satisfies the volume doubling condition. Suppose  that  $L$ is a non-negative self-adjoint operator
 which satisfies  the finite speed propagation
property
and condition  ${\rm (ST^{q}_{p, s})}$ for some $p,s,q$ such that
$1\leq p<s\leq \infty$ and $1\leq q\leq \infty$.
\begin{itemize}
\item[(i)] {\bf Compactly supported multipliers:}  Let   $F$ be  an even function such that $\supp F\subseteq [-1,1]$
 and $F \in W^{\beta,q}(\RR)$ for some $\beta>n(1/p-1/s)$. Then $F(\sqrt{L})$ is bounded on $L^p(X)$, and
$$ \sup_{t>0}\|F(t\sqrt{L})\|_{p\to p} \leq
C\|F\|_{W^{\beta,q}}.$$

\item[(ii)]{\bf  General  multipliers:}  Suppose $s= 2$ and  $F$ is an even  bounded  Borel
function such that $\sup_{t>0}\|\eta (\cdot) F(t \cdot)  \|_{W^{\beta, q}}<\infty $ for some
$\beta>\max\{n(1/p-1/2),1/q\}$
 and some non-trivial  function $\eta \in C_c^\infty(0,\infty)$. Then
$F(\sqrt{L})$ is bounded on $L^r(X)$ for all $p<r<p'$.
In addition,
\begin{eqnarray*}
   \|F(\sqrt{L})  \|_{r\to r}\leq    C_\beta\Big(\sup_{t>0}\|\eta (\cdot) F(t \cdot) \|_{W^{\beta, q}}
   + |F(0)|\Big).
\end{eqnarray*}
\end{itemize}
}

Assertion ${\rm (i)}$ of the theorem is inspired by  Guillarmou, Hassell and Sikora \cite{GHS} where a related  result is proved
under the assumption that  the  volume is  polynomial. Assertion ${\rm (ii)}$ is in the spirit of   H\"ormander's multiplier
theorem  for the Euclidean Laplacian.  Here, if $p > 1$, the order of differentiability required on $F$ is  smaller
since we do not search for
boundedness of $F(\sqrt{L})$ on $L^r$ for all $r \in (1, \infty)$.

The proof of assertion ${\rm (i)}$ makes  heavy use of the finite speed propagation property. This property together with the
classical dyadic decomposition of $F$ allow to reduce the problem of boundedness of $F(\sqrt{L})$ on $L^p$  to boundedness
of certain  compactly supported operators. The Stein-Tomas restriction type  condition  will be used to obtain an
 $L^p-L^s$ estimate of these operators from which we recover the  boundedness of $F(\sqrt{L})$ on $L^p$.

Assertion ${\rm (ii)}$ appeals as expected to singular integral theory.  We shall also make use of the estimate from assertion ${\rm (i)}$ since
$F(\sqrt{L})$ can be written as  the sum $\sum F_j(\sqrt{L})$ with compactly supported functions $F_j$.  However the operators
$F_j(\sqrt{L})$ do not act independently of each other and hence $L^p$ estimate for $F(\sqrt{L})$ does not hold in a trivial way from
the corresponding estimates for $F_j(\sqrt{L})$.  As explained by Littman, McCarthy and Rivi\`ere \cite{LMR},
we may have $F_j(\sqrt{-\Delta})$ to be  uniformly bounded  on $L^1(\RR^n)$ but
$F(\sqrt{-\Delta})$ fails  to be a multiplier of any $L^p$ other than $L^2$.  This  problem of recovering $L^p$ bounds for  $F(\sqrt{-\Delta})$ from
those for  $F_j(\sqrt{-\Delta})$  is discussed  by Carbery, Seeger and Sogge in  \cite{Carb1, S3, SS}.
We shall follow  closely Carbery  \cite{Carb1} and adapt some ideas  there to our abstract setting.

Our restriction type estimate does not hold  when the set of point spectrum is not empty. In particular, it does not
 hold for elliptic operators on compact manifolds or for the harmonic oscillator.  
 In order to treat these situations as well we modify the restriction estimate  as follows:
 for a fixed natural number $\kappa$ and for all $N\in \NN$ and  all even
  Borel functions $F$  such that\, $\supp F\subseteq [-N, N]$,
$$
\big\|F(\SL)P_{B(x, r)} \big\|_{p\to s} \leq
 CV(x,r)^{{1\over s}-{1\over p}}(Nr)^{n({1\over p}-{1\over s})}\| F (N \cdot) \|_{N^\kappa,\, q},
\leqno{\rm (SC^{q, \kappa}_{p, s})}
$$
for all $x\in X$ and all $r\geq 1/N$ where 
 $$
\|F\|_{N,q}=\left({1\over 2N}\sum_{\ell=1-N}^{N} \sup_{\lambda\in
 [{\ell-1\over N}, {\ell\over N})} |F(\lambda)|^q\right)^{1/q}
 $$
for $F$  supported in $[-1, 1]$.  For $q=\infty$, we put
 $\|F\|_{N, \infty}=\|F\|_{{\infty}}$. The  norm $\|F\|_{N,q}$ was used by Cowling and Sikora \cite{CowS} and Duong,
Ouhabaz and Sikora \cite{DOS} in the setting of spectral multipliers.

In some situations,   ${\rm (SC^{2,1}_{p, 2})}$ is equivalent to the following condition introduced by Sogge (see \cite{Sog1, Sog3, Sog4})
 $$
\big\|E_{\sqrt{L}}[\lambda, \lambda+1)\big\|_{p\to p'} \leq
C (1+\lambda)^{n({1\over p}-{1\over p'})-1}. \leqno{\rm (S_p)}
$$
We call ${\rm (SC^{q, \kappa}_{p, s})}$ {\it Sogge's spectral cluster condition}.  In this context we shall prove the following result
(see Theorems \ref{th3.2} and \ref{th31.2} for  precise statements).

\medskip
\noindent{\bf Theorem B.} {\it
Suppose that $X$ has finite measure and satisfies the volume doubling condition. Let $L$ be  a non-negative self-adjoint operator which
satisfies  the finite speed propagation
property  and Sogge's spectral cluster  condition ${\rm (SC^{q,1}_{p, s})}$   for some $p,s,q$ such that
$1\leq p<s\leq \infty$ and $1\leq q\leq \infty$. Then both  assertions of Theorem A hold  provided
$\beta> \max\{n(1/p-1/s), 1/q\}$.

The same conclusion holds in the case where  $\mu(X) = \infty$ provided  ${\rm (SC^{q,\kappa}_{p, s})}$ and an a priori estimate
for $\| F(\sqrt{L}) \|_{p\to p}$ are  satisfied. }

\medskip

As for Theorem A, an  appropriate decomposition of $F(\sqrt{L})$ as the sum of operators with compact supports
is the backbone of our arguments in proving  boundedness of $F(\sqrt{L})$ on $L^p$ for compactly supported $F$. Passing from compactly supported
multipliers to the general case will be done in the same way as for Theorem A. The proof of this part does not make explicit  use of ${\rm (ST^{q}_{p, 2})}$
or ${\rm (SC^{q,\kappa}_{p, 2})}$ but the  rather weaker condition
$$\big\|(I+t\SL)^{-N}P_{B(x, r)}\big\|_{p\to 2} \leq
CV(x, r)^{{1\over 2}-{1\over p}} \left({r\over  {t}}\right)^{n({1\over p}-{1\over 2})}, \; \; x \in X, \; r \geq  t>0.$$

Starting now from Theorem A or  Theorem B  with  $s = 2$ and choosing  the function $ F = S^{\delta}_R$ yields Bochner-Riesz summability  on $L^p(X)$
for   $\delta  > \delta_q(p) $  where
$$
\delta_q(p)=\max\Big\{0,n\Big|\frac{1}{p}-\frac{1}{2}\Big|-\frac{1}{q}\Big\}.
$$
Now we address the question of endpoint estimates, i.e., estimate for $S_R^{\delta_q(p)}(L)$.
It turns out that  our  Stein-Tomas restriction type   condition or  Sogge's cluster condition
imply that  Bochner-Riesz means are weak-type $(p,p)$ operators  for
 $\delta = \delta_q(p)$.  More precisely we obtain

\medskip
\noindent{\bf Theorem C.} {\it
Assume that $X$ satisfies the doubling condition and operator $L$ satisfies  the
 finite speed propagation property.
 \begin{itemize}
\item[(i)]
  If the restriction condition ${\rm (ST^{q}_{p, 2})}$ holds  for some $p, q$ satisfying
 $1\leq p <2$ and $1\leq q \leq \infty$ then $S_R^{\delta_q(p)} (L)$ is of weak-type $(p,p)$
uniformly in $R$.

\item[(ii)] The same conclusion as in ${\rm (i)}$ holds if  $\mu(X) < \infty$ and ${\rm (SC^{q,1}_{p, 2})}$  is satisfied  for some $p, q$ satisfying
 $1\leq p <2$ and $1\leq q \leq \infty$.
 \end{itemize}}

 \medskip
In the Euclidean case, it is known that $ S_R^{\delta_2(p)} (-\Delta)$
is not bounded on $L^p(\RR^n)$ for $p \not=2$. This was observed  by Christ and Sogge \cite{ChS}
who also proved weak-type $(1,1)$ for $S_R^{\delta_2(1)} (-\Delta)$.
Weak-type $(p,p)$ estimates of  $ S_R^{\delta_2(p)} (-\Delta)$
are proved by Christ \cite{C1, C2} when  $p < \frac{2n + 2}{n+3}$.
The corresponding  result  on compact
 manifolds is proved by  Seeger \cite{Se1}. The endpoint estimates
 for $p = \frac{2n+2}{n+3}$ are proved by Tao \cite{T2} both for $\RR^n$ and compact manifolds.

Our approach for  endpoint estimates is inspired by
Christ   and Tao  \cite{C1, C2, T2}. It   is based on $L^2$ Calder\'on-Zygmund
techniques (as used in Fefferman \cite{F1}), a spacial decomposition of the
Bochner-Riesz multiplier and the fact that if   $F$ has its inverse Fourier transform supported on a
set of width $R$, then by the finite speed propagation property the
operator $F(\SL)$ is supported in a $CR-$ neighbourhood   of the diagonal.
 It is worth to note that our proof of endpoint estimates
does not require any cancellation argument. This allows us  to consider applications to
operators with non-smooth kernels.

The  previous theorems are proved in Part 1 of this paper. In Part 2, we investigate the relation  of
${\rm (ST^{2}_{p, 2})}$ to  dispersive or Strichartz estimates for the corresponding
Schr\"odinger equation
\begin{equation}\label{I3}
\partial_t  u +  i L u = 0, \, \; u(0) = f \in L^2.
\end{equation}
In the setting of  Euclidean Laplacian, Strichartz's original proof for $L^p(\RR \times \RR^n)$ estimates of the solution
$u$ of (\ref{I3}) uses  restriction estimates  of the Fourier transform. In some sense we want to do
 the converse here, we want to take advantage of known dispersive or  Strichartz estimates for  (\ref{I3})
to prove a   Stein-Tomas  restriction type condition and then obtain  sharp spectral multipliers by Theorem A.
We are able do this either directly from dispersive estimate for $e^{itL}$ or from endpoint Strichartz estimate.
We prove the following result.

\medskip
\noindent{\bf Theorem D.} {\it
\begin{itemize}
\item[(i)]   Suppose that $L$ satisfies the Strichartz estimate
$$\int_{\RR}  \| e^{it L} f \|_{{\frac{2n}{n-2}}}^2 dt \le C \| f \|_{2}^2,  \, \, f \in L^2 $$
for some $n > 2$.  Assume also that the smoothing  property
$$\|\exp(-tL)\|_{p \to \frac{2n}{n+2} }\le Kt^{-\frac{n}{2} (\frac{1}{p} - \frac{n+2}{2n} )},$$
holds for all  $p \in [1,   \frac{2n}{n+2}]$.
Then for all  $\lambda \ge 0$
$$\|d E_{\sqrt{L}}(\lambda)\|_{p \to p' }\le C\lambda^{n(\frac{1}{p} -\frac{1}{p'}) -1}.$$

\item[(ii)] Fix $p \in [1,   \frac{2n}{n+2}]$.  Suppose that $X$  satisfies the doubling condition and  that
there exists   a  positive constants $C > 0$ such that
$V(x, r)\leq C r^n$ for every $x\in X $ and $r>0$.
Assume that  $L$ satisfies  the finite speed propagation property together with  Strichartz  and smoothing estimates as in ${\rm (i)}$.
Then for every even compactly supported   bounded function $F$ such that
 $ \|   F \|_{W^{\beta,2}} < \infty$
 for some $\beta > n (\frac{1}{p} - \frac{1}{2})$, the operator
 $F(\sqrt{L})$ is bounded on $L^p$ and
 $$\sup_{t > 0} \| F(t\sqrt{L}) \|_{p\to p} \le C \| F \|_{W^{\beta,2}}.$$

\item[(iii)] Suppose that the  conditions of {\rm (ii)} are satisfied. Then for any even bounded Borel function $F$ such that
$\sup_{t>0}\|\eta(.)  F(t \cdot) \|_{W^{\beta, 2}}<\infty $ for some  $\beta>\max\{n(1/p-1/2),1/2\}$ and some
non-trivial function $\eta \in C_c^\infty(0,\infty)$, the operator
$F(\SL)$ is bounded on $L^r(X)$ for all $r \in (p, p')$.
\end{itemize}}

\medskip

The main assertion here is ${\rm (i)}$. Indeed,  once ${\rm (i)}$  is proved we obtain a  Stein-Tomas restriction type
estimate  and then appeal to Theorem A to prove  assertions
${\rm (ii)}$ and ${\rm (iii)}$. We can also  replace the Strichartz estimate by a dispersive estimate
$$ \| e^{it L} \|_{1 \to \infty} \le C | t |^{-n/2}, \ \ \ \ \  t \in \RR, \ t \not= 0.$$
Note that by a result of Keel and Tao \cite{KT}, endpoint Strichartz estimate follow from this dispersive estimate.

Strichartz estimates have been studied by several authors. For example,  Burq, Planchon, Stalker and
 A. Tahvildar-Zadeh \cite{BPST} proved such estimates for
Schr\"odinger operators with inverse square potentials, i.e. $L = -\Delta + \frac{c}{|x|^2}$ on $\RR^n$.  Therefore we obtain sharp multiplier results as
well as endpoint Bochner-Riesz estimates for these   operators. It is worth to mention that
if $ -{(n-2)^2/4} < c < 0 $,   the semigroup $\exp(-tL)$ acts on $L^p(\RR^n)$ only for $p \in (p_c', p_c)$ with   $p_c  < \infty$. In particular, the corresponding
heat kernel does not enjoy any good upper bounds such as Gaussian upper bounds. Nevertheless  we obtain sharp spectral multipliers for $L$ on $L^p$
for $p\in (p_c', 2n/(n+2)]$.
We discuss in Part 3 several other examples to which  Theorems A, B, C and D apply. This includes radial
 Schr\"odinger operators with inverse square potentials, the harmonic oscillator, elliptic operators on compact manifolds, Laplacian on asymptotically conic manifolds.

 While this paper was finished we learned that M. Uhl introduced recently in
his PhD Thesis \cite{U} a condition similar to our restriction type condition  and proved a spectral
multiplier result similar to our Theorems 4.1 and  4.2. However the order of
differentiability required in Uhl's result is $\beta > n/2$ and hence it
is less sharp than our Theorems 4.1 and  4.2.

\bigskip

\part{Restriction estimates imply sharp spectral multipliers }

\medskip

\section{Restriction type condition}
\setcounter{equation}{0}

We start by fixing some notation and assumptions.  Throughout this paper,
unless we mention the contrary, $(X,d,\mu)$ is a metric measure  space, that is, $\mu$
is a Borel measure with respect to the topology defined by the metric $d$.
We denote by
$B(x,r)=\{y\in X,\, {d}(x,y)< r\}$  the open ball
with centre $x\in X$ and radius $r>0$. We   often just use $B$ instead of $B(x, r)$.
Given $\lambda>0$, we write $\lambda B$ for the $\lambda$-dilated ball
which is the ball with the same centre as $B$ and radius $\lambda r$.
We set $V(x,r)=\mu(B(x,r))$ the volume of $B(x,r)$ and we say that $(X, d, \mu)$ satisfies
 the doubling property (see Chapter 3, \cite{CW})
if there  exists a constant $C>0$ such that
\begin{eqnarray}
V(x,2r)\leq C V(x, r)\quad \forall\,r>0,\,x\in X. \label{eq2.1}
\end{eqnarray}
If this is the case, there exist  $C, n$ such that for all $\lambda\geq 1$ and $x\in X$
\begin{equation}
V(x, \lambda r)\leq C\lambda^n V(x,r). \label{eq2.2}
\end{equation}
In the sequel we want to consider $n$ as small as possible.
Note that in general one cannot take infimum over such exponents $n$ in \eqref{eq2.2}.
In the Euclidean space with Lebesgue measure, $n$ corresponds to
the dimension of the space.
Observe that if $X$ satisfies (\ref{eq2.1}) and has finite measure then it has finite diameter (see, e.g., \cite{AMP}).
 Therefore if $\mu(X)$ is finite, then we may assume that $X=B(x_0, 1)$ for some $x_0\in X$.

  For $1\le p\le+\infty$, we denote the
norm of a function $f\in L^p(X,{\rm d}\mu)$ by $\|f\|_p$, by $\langle .,. \rangle$
the scalar product of $L^2(X, {\rm d}\mu)$, and if $T$ is a bounded linear operator from $
L^p(X, {\rm d}\mu)$ to $L^q(X, {\rm d}\mu)$, $1\le p, \, q\le+\infty$, we write $\|T\|_{p\to q} $ for
the  operator norm of $T$.
Given a  subset $E\subseteq X$, we  denote by  $\chi_E$   the characteristic
function of   $E$ and  set
$$
P_Ef(x)=\chi_E(x) f(x).
$$
For a given  function $F: {\mathbb R}\to {\mathbb C}$ and $R>0$, we define the
function
$\delta_RF:  {\mathbb R}\to {\mathbb C}$ by putting
 $\delta_RF(x)= F(Rx).$

\medskip

\noindent {\bf 2.1.\, Finite speed propagation for the wave equation.}
Set
\begin{equation*}
\D_\rho=\{ (x,\, y)\in X\times X: {d}(x,\, y) \le \rho \}.
\end{equation*}
Given an  operator $T$ from $L^p(X)$ to $L^q(X)$, we write
\begin{equation}\label{eq2.3}
\supp  K_{T} \subseteq \D_\rho
\end{equation}
 if $\langle T f_1, f_2 \rangle = 0$ whenever $f_k$ is
in~$C(X)$ and has support $\supp f_k \subseteq B(x_k,\rho_k)$ when $k = 1,2$,
and $\rho_1+\rho_2+\rho < {d}(x_1, x_2)$.   Note that  if $T$ is an
integral operator with a  kernel $K_T$, then (\ref{eq2.3}) coincides
with the  standard meaning of $\supp  K_{T}  \subseteq \D_\rho$,
 that is $K_T(x, \, y)=0$ for all $(x, \, y) \notin \D_\rho$.

Given  a non-negative self-adjoint operator $L$
on $L^2({X})$. We say that $L$   satisfies the finite speed
propagation property if
$$
\supp  K_{\cos(t\SL)} \subseteq \D_t \quad \forall t> 0\,.
\leqno{\rm (FS)}
$$
Property (FS) holds for most of second order self-adjoint operators  and  is equivalent to Davies-Gaffney
estimates. See, for example \cite{CGT}, \cite{S2} and  \cite{CouS}.

We recall the following well-known simple lemma.

\begin{lemma}\label{le2.1}
Assume that $L$ satisfies {\rm (FS)} and that $F$ is an even bounded Borel function with Fourier
transform  $\hat{F} \in L^1(\RR)$ and that
$\mbox{\rm supp} \hat{F} \subseteq [-\rho,\, \rho]$.
Then
$$
\mbox{\rm supp} K_{F(\SL)} \subseteq \D_\rho.
$$
\end{lemma}

\noindent {\bf Proof.}
If $F$ is an even function, then by the Fourier inversion formula,
$$
F(\SL) =\frac{1}{2\pi}\int_{-\infty}^{+\infty}
  \hat{F}(t) \cos(t\SL) \;dt.
$$
But $\supp\hat{F} \subseteq [-\rho,\rho]$,
and the lemma follows then from {\rm (FS)}. \hfill{} $\Box$

\medskip

\noindent
{\bf 2.2.\, The Stein-Tomas restriction type condition.}
Assume that $(X,d, \mu)$ satisfies the doubling condition, that is
(\ref{eq2.2}). Consider a non-negative self-adjoint operator $L$  and numbers
$p, s$ and $q$ such that $1\leq p< s\leq \infty$ and $1\leq
q\leq\infty$. We say that  $L$ satisfies the
    {\it  Stein-Tomas restriction type condition}   if:
  for any $R>0$ and all Borel functions $F$ such that $\supp F \subset [0,R]$,
$$
\big\|F(\SL)P_{B(x, r)} \big\|_{p\to s} \leq CV(x,
r)^{{1\over s}-{1\over p}} \big( Rr \big)^{n({1\over p}-{1\over
s})}\big\|\delta_RF\big\|_{q}
\leqno{\rm (ST^q_{p, s})}
$$
  for all
  $x\in X$ and all $r\geq 1/R$.

\begin{remark}\label{re2.2}
 Note that   if condition ${\rm (ST^{q}_{p, s})} $ holds for  some $q\in [1, \infty)$,
 then ${\rm (ST^{\tilde{q}}_{p, s})} $  holds for all $\tilde{q}\geq q$ including the case
 $\tilde{q}=\infty$.
\end{remark}

\begin{proposition}\label {prop2.3} Suppose that $(X, d, \mu)$ satisfies
 property  \eqref{eq2.2}. Let  $1\leq p<2 $ and
$N>n(1/p-1/2)$. Then   ${\rm (ST^{\infty}_{p, 2})}$   is equivalent
to each of the following conditions:
\begin{itemize}
\item[(a)]  For all $x\in X$ and  $r\geq t>0$    we have
$$
\big\|e^{-t^2L}P_{B(x, r)}\big\|_{p\to 2} \leq
CV(x, r)^{{1\over 2}-{1\over p}} \Big({r\over  {t}}\Big)^{n({1\over p}-{1\over 2})}.
\leqno{\rm (G_{p, 2})}
$$

\item[(b)]
For all $x\in X$ and    $r\geq  t >0$   we have
$$
\big\|(I+t\SL)^{-N}P_{B(x, r)}\big\|_{p\to 2} \leq
CV(x, r)^{{1\over 2}-{1\over p}} \left({r\over  {t}}\right)^{n({1\over p}-{1\over 2})}.
\leqno{\rm (E_{p, 2})}
$$
\end{itemize}
\end{proposition}

  \noindent
{\bf Proof.}  We shall show  that
${\rm (E_{p, 2})}\Rightarrow {\rm (ST^{\infty}_{p, 2})}\Rightarrow
 {\rm (G_{p, 2})}\Rightarrow {\rm (E_{p, 2})}$.

 Suppose that   $F$ is   a  Borel function with  $\supp F \subset [0,R]$.
Let $1\leq p<2 $ and $N>n(1/p-1/2)$.  It follows from   ${\rm (E_{p, 2})}$ that
 for every $x\in X$ and $r\geq  1/R$,
\begin{eqnarray*}
\big\|F(\SL)P_{B(x, r)}\|_{p\to 2}
  &=&\big\|F(\SL)\big(I + {\sqrt{L}\over R }\big)^N  \big(I + {\sqrt{L}\over R }\big)^{-N}
 P_{B(x, r)}\big\|_{p\to 2}\\
  &\leq&
 \sup_{\lambda }\big| F(\lambda) \big(1 + {\lambda\over R} \big)^N \big| \, \cdot
\big\|\big(I + {\sqrt{L}\over R }\big)^{-N} P_{B(x, r)}\big\|_{p\to 2}
\\
  &\leq& C2^N V(x, r)^{{1\over 2}-{1\over p}} \big(Rr)^{n({1\over p}-{1\over 2})}
 \|\delta_RF \|_{\infty}.
\end{eqnarray*}
This gives condition   ${\rm (ST^{\infty}_{p, 2})}$.

\smallskip

Next assume  ${\rm (ST^{\infty}_{p, 2})}$. Then
\begin{eqnarray*}
  \left \|e^{-t^2L}P_{B(x, r)} \right\|_{p\to 2}&=& \left\|\int_0^{\infty}     e^{-t^2\lambda}
dE_{L}(\lambda)P_{B(x, r)}  \right\|_{p\to 2}
\\&=&   \left\| \int_0^{\infty}  t^2    e^{-t^2\lambda}
 E_{L}[0, \lambda]P_{B(x, r)}  d\lambda  \right\|_{p\to 2}
 \\ &\le& \int_0^{\infty}    t^{2}e^{-t^2\lambda}
 \left\|\chi_{[0,\sqrt \lambda]} (\sqrt{L}) P_{B(x, r)} \right\|_{p\to 2}\, {d\lambda}.
\end{eqnarray*}
Now if $0<\lambda \leq 1/r^2$, then $B(x,r)\subseteq B(x, \lambda^{-1/2})$.  Therefore,
 by ${\rm (ST^{\infty}_{p, 2})}$
\begin{eqnarray*}
 \|\chi_{[0,\sqrt \lambda]} (\sqrt{L}) P_{B(x, r)} \|_{p\to 2} &\le&
 \|\chi_{[0,\sqrt \lambda]} (\sqrt{L}) P_{B(x, \lambda^{-1/2})} \|_{p\to 2} \\
 &\le& C V(x,\lambda^{-1/2})^{{1\over 2}-{1\over p}} \le C V(x,r)^{{1\over 2}-{1\over p}}.
 \end{eqnarray*}
If $\lambda > 1/r^2$, then by ${\rm (ST^{\infty}_{p, 2})}$
$$ \|\chi_{[0,\sqrt \lambda]} (\sqrt{L}) P_{B(x, r)} \|_{p\to 2}  \le C V(x,r)^{{1\over 2}-{1\over p}}
(r \lambda^{1/2})^{n({\frac{1}{ p}}-{\frac{1}{ 2}})}.$$
This proves ${\rm (G_{p, 2})}$.

\smallskip

To finish the proof  assume that ${\rm (G_{p, 2})}$ holds. Then
\begin{eqnarray*}
\big\|(I+t\SL)^{-N}P_{B(x, r)}\|_{p\to 2}
&\leq&\big\|{(I+t^2L)^{{N/2}} }  (I+t\SL)^{-N}
\big(I+t^2L\big)^{-{N/2}}    P_{B(x, r)}\big\|_{p\to 2}\nonumber\\
 &\leq&
 \sup_{\lambda }\big| {(1+t^2\lambda^2)^{{N/2}}} (1+t\lambda)^{-N}   \big| \, \cdot
 \big\|\big(I+t^2L\big)^{-{N\over 2}}    P_{B(x, r)}\big\|_{p\to 2} \nonumber\\
 &\leq&  C\big\|\big(I+t^2L\big)^{-{N/2}}    P_{B(x, r)}\big\|_{p\to 2}.
\end{eqnarray*}
Next note that  for  $t > 0$
\begin{eqnarray}\label{eq2.5}\nonumber
\left\|(I+t^2L)^{-{N/2}}P_{B(x, r)}\right\|_{p\to 2}
 &=& \left\|C_N\int_0^{\infty}    e^{-s} s^{{N/2}-1}
 e^{-st^2L}P_{B(x, r)}ds\right\|_{p\to 2}
 \\&\leq&C_N\int_0^{\infty}    e^{-s} s^{{N/2}-1}
 \left\|e^{-st^2L}P_{B(x, r)}\right\|_{p\to 2}{ds}.
\end{eqnarray}
Hence if  $s \le r^2/t^2$,  then  by ${\rm (G_{p, 2})}$
$$ \|e^{-st^2L}P_{B(x, r)} \|_{p\to 2} \le V(x,r)^{{1\over 2}-{1\over p}}
\big({r\over ts^{1/2} }\big)^{n{({1\over p}-{1\over 2}})}.$$
If $s >  r^2/t^2$, then  $B(x,r)\subseteq B(x, ts^{1/2})$ and ${\rm (G_{p, 2})}$   implies
that
$$ \|e^{-st^2L}P_{B(x, r)} \|_{p\to 2}  \le \|e^{-st^2L}P_{B(x, ts^{1/2})}
\|_{p\to 2} \le C V(x, ts^{1/2})^{{1\over 2}-{1\over p}} \le
C V(x,r)^{{1\over 2}-{1\over p}}.$$ Using  these estimates in
(\ref{eq2.5})  yields ${\rm (E_{p, 2})}$ for $N>n(1/p-1/2)$.  This
ends the proof.
 \hfill{} $\Box$

\medskip

It is natural to generalise
condition \eqref{I1}
to abstract self-adjoint operators in the following way (see \cite{GHS}).
One says that $L$ satisfies {\it $L^p$ to $L^{p'} $ restriction estimates} if  the spectral
measure $dE_{\sqrt{L}}(\lambda)$ maps  $L^p(X)$ to $L^{p'}(X)$ for some $p<2$, with an operator norm estimate
$$
\big\|dE_{\sqrt{L}}(\lambda)\big\|_{p\to p'}\leq C \lambda^{n({1\over p}-{1\over p'})-{1}}
\leqno{\rm (R_p)}
$$
for all $\lambda\ge 0$, where $n$ is
   as in  condition (\ref{eq2.2}) and
 $p'$ is conjugate of $p$, i.e., ${1/p} +{1/p'}=1$.

 \medskip

\begin{proposition}\label{prop2.4} Fix $1\le p< {2n/(n+1)}$.
Suppose that  there exists a  constant $C>0 $  such that
$C^{-1}r^n \leq V(x, r)\leq C r^n$
for all $x\in X $ and $r>0$.  Then  conditions ${\rm (R_p)}$, ${\rm (ST^{2}_{p, 2})}$ and
 ${\rm (ST^{1}_{p, p'})}$
are equivalent.
\end{proposition}

 \noindent
 {\bf Proof.} The proof  is inspired by  estimates (2.12) of \cite{GHS}. We first  show the implication
${\rm (R_p)}\Rightarrow {\rm (ST^{1}_{p, p'})}$.
Suppose that $F$  is a  Borel function    such that $\supp F \subset [0,R]$
for some $R>0$.  Then by ${\rm (R_p)}$
\begin{eqnarray*}
\big\|F(\SL) P_{B(x,r)} \big\|_{p\to p'}
 &\leq&  \int_0^\infty
|F(\l)| \|dE_{\sqrt{L}}(\l)\|_{p\to
p'}d\l\\
 &\leq&  C\int_0^R |F(\l)| \lambda^{n({{\frac{1}{p}}-{\frac{1}{p'}}})-1}d\l \\
 &\leq&  CR^{n({{\frac{1}{p}}-{\frac{1}{p'}}})-1}\int_0^R |F(\l)|  d\l \\
 &\leq&  C V(x,r)^{{\frac{1}{p'}}-{\frac{1}{p}}}  (rR)^{n({{\frac{1}{p}}-{\frac{1}{p'}}})} \|\delta_RF\|_{1},
\end{eqnarray*}
where in the last inequality we used the assumption that  $V(x, r) \le C  r^n$.

Next we prove that ${\rm (ST^{1}_{p,p'})}\Rightarrow {\rm (ST^{2}_{p, 2})}$. Note that
$  V(x, r)\sim r^n$ for every $x\in X$ and $r>0.$
Letting $r \to \infty$ we obtain from ${\rm (ST^{1}_{p,p'})}$
$$ \big\|F(\SL)\big\|_{p\to p'} \le CR^{n({{1\over p}-{1\over p'}})}\|\delta_RF\|_{1}.$$
By $T^{\ast} T$ argument
$$
 \big\|F(\SL)\big\|^2_{p\to 2}
 =   \big\| |F|^2(\SL)\big\|_{p\to p'}
 \leq  CR^{2n({{1\over p}-{1\over 2}})}\|\delta_RF\|^2_{2}.
$$
Hence
\begin{eqnarray*}
 \big\|F(\SL)P_{B(x, r)}\big\|_{p\to 2} \leq  \big\|F(\SL)\big\|_{p \to 2}
 &\leq& CV(x,r)^{{\frac{1}{2}}-{\frac{1}{p}}}(Rr)^{n({{\frac{1}{p}}-{\frac{1}{2}}})}\|\delta_RF\|_{2}.
\end{eqnarray*}
 This gives ${\rm (ST^{2}_{p, 2})}$.

 Now, we   prove the remaining implication
${\rm (ST^{2}_{p, 2})}\Rightarrow  {\rm (R_p)}$. By
volume estimate   $V(x, r)\geq C^{-1} r^n$  and condition  ${\rm (ST^{2}_{p, 2})}$
\begin{eqnarray}\label{eq2.6}
\big\|F(\SL)P_{B(x, r)} \big\|_{p\to 2} \leq
C R^{n({1\over p}-{1\over 2})}\big\|\delta_RF\big\|_{2}
\end{eqnarray}
  for any $R>0$, all Borel functions $F$ such that $\supp F \subset [0,R]$,  all
  $x\in X$ and $r\geq 1/R$. Taking the limit $r\to \infty$ gives  \smallskip
\begin{eqnarray}\label{eq2.7}
\big\|F(\SL)  \big\|_{p\to 2} \leq
C R^{n({1\over p}-{1\over 2})}\big\|\delta_RF\big\|_{2}.
 \end{eqnarray}
 \noindent
 For $F=\chi_{(\lambda-\varepsilon, \lambda+\eps]}$ and $R=\lambda+\eps$ in  (\ref{eq2.7}) yields
 \begin{eqnarray*}
 \Big\|\varepsilon^{-1}{E_{\sqrt{L}}(\lambda-\eps, \lambda+\varepsilon]  } \Big\|_{p\to p'}
 &=&\varepsilon^{-1}\Big\| {E_{\sqrt{L}}(\lambda-\eps, \lambda+\varepsilon]  } \Big\|^2_{p\to 2} \\
 &\leq&
C \varepsilon^{-1}(\lambda+\eps)^{2n({1\over p}-{1\over 2})}
\big\| \chi_{({\lambda-\eps\over \lambda+\eps}, \, 1]}\big\|^2_{2} \\
 &\leq& C(\lambda+\eps)^{n({1\over p}-{1\over p'})-1}.
 \end{eqnarray*}
 Taking $\varepsilon\to 0$ yields condition ${\rm (R_p)}$ (see Proposition 1, Chapter XI, \cite{Yo}).
\hfill{} $\Box$

\bigskip

\section{Sharp spectral multipliers - compactly supported functions}\label{sec3}
\setcounter{equation}{0}

 In this section  we show that the restriction type condition which we introduce
  in the previous section can be used to obtain  sharp spectral multiplier results
  in the abstract setting of self-adjoint operators acting on homogeneous spaces.
   We first consider the case of compactly supported functions.  We assume here
   that  $(X, d, \mu)$ is  a  metric measure  space satisfying the doubling property
   and recall that $n$ is the  doubling dimension from  condition (\ref{eq2.2}).
   We use the standard notation $W^{\beta,q}(\RR) $ for the Sobolev space
$\|F\|_{W^{\beta,q}}  = \|(I-d^2/dx^2)^{\beta/2}F\|_{q}$. The first result and its proof are
inspired  by  Theorem 1.1 of  \cite{GHS}.

 \begin{theorem}\label{th3.1}
Suppose  that operator $L$ satisfies   property  {\rm (FS)}
and condition  ${\rm (ST^{q}_{p, s})}$ for some $p,s,q$ such that
$1\leq p<s\leq \infty$ and $1\leq q\leq \infty$.
 Next assume that  $F$ is an even function such that  $\supp F\subseteq [-1,1]$
 and $F \in W^{\beta,q}(\RR)$
for some $\beta>n(1/p-1/s)$. Then $F(t\sqrt{L})$ is bounded on $L^p(X)$ for all $t>0$. In addition,
\begin{equation}
\label{eq3.1} \sup_{t>0}\|F(t\sqrt{L})\|_{p\to p} \leq
C\|F\|_{W^{\beta,q}}.
\end{equation}
\end{theorem}

We described the proof of the Theorem \ref{th3.1} at the end of this section.

A standard application of spectral multiplier theorems is Bochner-Riesz means.
Such application is also a good test to check if the considered multiplier result is sharp.
Let us recall that Bochner-Riesz means of order $\delta$ for a non-negative self-adjoint operator $L$
are defined by the formula
\begin{equation}\label{e3.8}
S_R^{\delta} (L)  = \Big(I-{L\over R^2}\Big)_+^{\delta},\ \ \ \ R>0.
\end{equation}
\noindent
The case  ${\delta}=0$ corresponds to the spectral projector $E_{\SL}[0, R]$. For 
${\delta}>0$ we think of (\ref{e3.8}) as a smoothed version of this
spectral projector; the larger ${\delta}$, the more smoothing.
Bochner-Riesz summability  on $L^p$ describes the range of ${\delta}$ for which $S_R^{\delta} (L)$ are bounded on $L^p$,  uniformly in $R$.

In Theorem~\ref{th3.1}, if one chooses $F(\lambda) = (1-\l^2)^{\delta}_+$ then $F\in W^{\beta,q}$ if
and only if ${\delta}>\beta-1/q$. Therefore, we obtain

\begin{coro}\label{coro3.3}
Suppose that  the
 operator $L$ satisfies  the finite speed propagation property {\rm (FS)}
and condition  ${\rm (ST^{q}_{p, s})}$ with some $1\leq p<s\leq
\infty$ and $1\leq q\leq \infty$.
 Then for all  ${\delta}> n(1/p-1/s)-1/q$, we have
\begin{eqnarray}\label {e3.9}
\Big\|\Big(I-{L\over R^2}\Big)_+^{\delta}\Big\|_{p\to p}\leq C
\end{eqnarray}
uniformly in $R>0.$
\end{coro}

As a consequence, we obtain  the following necessary condition for the
restriction condition ${\rm (ST^{q}_{p, s})}$ (see also Kenig, Stanton and Tomas \cite{KST}).

\begin{coro}  \label {coro3.4} Suppose that ${1/q}> n({1/p}-{1/s})$ for some $q\geq 1$ and $1\leq p<s\leq \infty$.
Then condition  ${\rm (ST^{q}_{p, s})}$
implies that $L=0.$
\end{coro}

\noindent {\bf Proof.} Note that if ${1/q}> n({1/p}-{1/s})$ for some
$q\geq 1$ and $1\leq p<s\leq \infty$, then there exist ${\delta}<0$
and $\varepsilon>0$ such that
$S^{\delta}_1(\lambda^2)=(1-\l^2)^{\delta}_+\in W^{n({1\over
p}-{1\over s})+\varepsilon,q}$. By Theorem~\ref{th3.1},
 the operator $S^{\delta}_R(L)$  is bounded
on $L^p(X)$ uniformly
in $R$, i.e.,
$ \|S^{\delta}_R(L) \|_{p\to p}\leq C<\infty$ for some constant $C>0$ independent of $R$.
 However, $S^{\delta}_R(L)$ is a self-adjoint operator,
so $ \|S^{\delta}_R(L) \|_{p'\to p'}<\infty$, and by interpolation,
$ \|S^{\delta}_R(L) \|_{2\to 2}<\infty$. Let $M=1+ \|S^{\delta}_R(L)
\|_{2\to 2}$. It follows  that $S^{\delta}_R(\lambda)>M$ for all
$\lambda\in [R^2(1-M^{1/{\delta}}), \,  R^2].$ By spectral theorem,
this implies that  $E_{{L}}[R^2(1-M^{1/{\delta}}), \,  R^2]=0.$ As
$R$ is arbitrary positive number, this implies $L=0.$ \hfill{}
$\Box$

\medskip

\begin{remark}\label{re3.5}\,
Note that condition ${\rm (ST^{q}_{p, s})}$ allows us to define the operator
$S^{\delta}_R(L)$ even when $\delta<0$ when the function
$\lambda \to S_1^{\delta}(\lambda^2)$ is unbounded.
\end{remark}

We return to the discussion of Bochner-Riesz analysis in Section~\ref{BRS} and we now discuss a
discreet version of Theorem~\ref{th3.1}.

\bigskip

It is not difficult to see that   condition   ${\rm (ST^{q}_{p, s})}$  with  some   $ q< \infty$ implies that
  the set of point spectrum of $L$ is empty.   Indeed, one has for all  $0\leq a< R$ and $x\in X,$
$$
\big\|1\!\!1_{\{a \}  }(\sqrt{L})P_{B(x,r)}\big\|_{p\to s}
\leq C V(x,r)^{{1\over s}-{1\over p}} (rR)^{n({1\over p}-{1\over s})}\big\|1\!\!1_{\{a \} }(R\cdot)\big\|_{q} =0, \ \ \ Rr\ge 1
$$
and therefore $1\!\!1_{\{a\} }(\sqrt{L})=0$. Due to  $\sigma(L)\subseteq [0, \infty)$, it follows that the point spectrum of $L$
is empty.
In particular, ${\rm (ST^{q}_{p, s})}$ cannot hold for any $q<\infty$ for elliptic
operators on compact manifolds or for the harmonic oscillator.
To be able to  study these operators as well,  we introduce a
variation of  condition ${\rm (ST^{q}_{p, s})}$. Following  \cite{CowS, DOS}, for an even Borel function $F$
with  $\supp F\subseteq [-1, 1]$ we define
the norm $\|F\|_{N,q}$ by
$$
\|F\|_{N,q}=\left({1\over 2N}\sum_{\ell=1-N}^{N} \sup_{\lambda\in
 [{\ell-1\over N}, {\ell\over N})} |F(\lambda)|^q\right)^{1/q},
 $$
 where $q\in [1, \infty)$ and $N\in \NN$. For $q=\infty$, we put
 $\|F\|_{N, \infty}=\|F\|_{{\infty}}$.
 It is obvious that $\|F\|_{N,q}$ increases monotonically in $q$.

     Consider a non-negative self-adjoint operator $L$  and numbers
$p, s$ and $q$ such that $1\leq p< s\leq \infty$ and $1\leq
q\leq\infty$.  We shall  say  that   $L$
  satisfies the  {\it Sogge spectral cluster condition}  if: for  a fixed natural number $\kappa$ and
 for all $N\in \NN$ and  all even
  Borel functions $F$  such that\, $\supp F\subseteq [-N, N]$,
$$
\big\|F(\SL)P_{B(x, r)} \big\|_{p\to s} \leq
 CV(x,r)^{{1\over s}-{1\over p}}(Nr)^{n({1\over p}-{1\over s})}\|\delta_NF\|_{N^\kappa,\, q}
\leqno{\rm (SC^{q,\kappa}_{p, s})}
$$
for  all $x\in X$ and $r\geq 1/N$.  For $q=\infty$, ${\rm (SC^{\infty,\kappa}_{p, s})}$ is independent of $\kappa$  so we write
it as ${\rm (SC^{\infty}_{p, s})}$.

\medskip

\begin{remark}\label{re3.55}\, It is easy to check that for $\kappa\geq 1$, ${\rm (SC^{q,\kappa}_{p, s})}$ implies
${\rm (SC^{q,1}_{p, s})}$.
\end{remark}

\begin{theorem}\label{th3.2} Suppose the operator $L$ satisfies   property {\rm (FS)}
and condition  ${\rm (SC^{q,\kappa}_{p, s})}$ for a fixed
$\kappa \in {\mathbb N}$ and some $p,s,q$ such that $1\leq p<s\leq
\infty$ and $1\leq q\leq \infty$. In addition, we assume that for
any $\varepsilon>0$ there exists a constant $C_\varepsilon$ such
that for all $N\in {\mathbb N}$ and all even Borel functions $F$
such that supp $F\subset [-N,N]$,
$$
\|F(\SL)\|_{p\to p}\leq C_\varepsilon N^{\kappa n({1\over p}-{1\over
s})+\varepsilon}\|\delta_N F\|_{N^{\kappa},q}. \leqno{\rm (AB_p)}
$$
Then for any even function $F$ such that $\supp F \subseteq [-1,1]$ and $\|F\|_{W^{\beta,q}}<\infty$
 for some $\beta>\max\{n(1/p-1/s),1/q\}$, the operator $F(t\sqrt{L})$ is bounded on $L^p(X)$ for all $t>0$. In addition,
 \begin{equation}
\label{eq3.2} \sup_{t>0}\|F(t\sqrt{L})\|_{p\to p} \leq
C\|F\|_{W^{\beta,q}}.
\end{equation}
\end{theorem}

  Note that condition ${\rm (SC^{q,\kappa}_{p, s})}$ is weaker than ${\rm (ST^{q}_{p, s})}$ and we need  a priori estimate
  ${\rm (AB_p)}$ in Theorem~\ref{th3.2}.  See also \cite[Theorem 3.6]{CowS} and \cite[Theorem 3.2]{DOS} for related results. Once
  ${\rm (SC^{q,\kappa}_{p, s})}$ is proved, a priori estimate
  ${\rm (AB_p)}$  is not difficult to check in general, see for example the section on the harmonic oscillator.

 Following  \cite{DOS}, we  prove  the following result.

\begin{proposition}\label{le3.3333} Suppose that $\mu(X)<\infty$ and ${\rm (SC^{q,1}_{p, s})}$  for some $p,s,q$ such that $1\leq p<s\leq
\infty$ and $1\leq q\leq \infty$.
Then
\begin{eqnarray*}
\|F(\SL)\|_{p\to p} \leq  C  N^{n({1\over p}-{1\over s}) }\|\delta_N
F\|_{N,q}
\end{eqnarray*}
for all $N\in \NN$ and all Borel  functions $F$ such that supp $F\subseteq [-N, N]$. Therefore,
 for any even function $F$ such that $\supp F \subseteq [-1,1]$ and $\|F\|_{W^{\beta,q}}<\infty$
 for some $\beta>\max\{n(1/p-1/s),1/q\}$, the operator $F(t\sqrt{L})$ is bounded on $L^p(X)$ 
 for all $t>0$ and
$$\sup_{t>0}\|F(t\sqrt{L})\|_{p\to p} \leq
C\|F\|_{W^{\beta,q}}.$$
\end{proposition}

\noindent {\bf Proof.}  Since  $\mu(X)<\infty$, we may assume that $X=B(x_0, 1)$ for some $x_0\in X$ (see \cite{AMP}).
It  follows from H\"older's inequality and condition ${\rm (SC^{q,1}_{p, s})}$
 that
\begin{eqnarray*}
\|F(\SL)\|_{p\to p}&\leq& V(X)^{{1\over p}-{1\over s}}
\|F(\SL)P_{B(x_0, 1)}\|_{p\to s}\\
&\leq & C V(X)^{{1\over p}-{1\over s}}V(X)^{{1\over s}-{1\over
p}}N^{n({1\over p}-{1\over s})}\|\delta_NF\|_{N,\, q}\\
&\leq & C  N^{n({1\over p}-{1\over s})}\|\delta_N
F\|_{N,q}.
\end{eqnarray*}
This means that ${\rm (AB_p)}$ is satisfied and thus the  last assertion follows from Theorem \ref{th3.2}. This proves Proposition~\ref{le3.3333}.
\hfill{}
$\Box$

\medskip

The proof of   Theorems~\ref{th3.1}  and \ref{th3.2}  uses the following
lemma. In the case where  the volume is polynomial this lemma is proved in   \cite{GHS} using a similar argument.

\begin{lemma}\label{le3.6}  Suppose that $T$ is a  linear map such that
for all $x \in X$ and $r>0$ the operator
$ TP_{B(x,r)} $ is bounded  from $L^p (X)$
to $L^s(X)$ for some   $1\le p <s \le \infty $. Assume also that
$$
\supp K_{T} \subseteq \D_{\rho}
$$
for some $\rho>0$. Then
there exists a constant $C=C_{p,s}$ such that
$$
\|T\|_{{p}\to {p}} \le C\sup_{x\in X}\big\{V(x,\rho)^{{1\over p}-{1\over
s}} \|TP_{B(x,\rho)}\|_{p\to s}\big\}.
$$
\end{lemma}

\noindent {\bf Proof.}
We fix $\rho>0$. Then  we choose a sequence $(x_n)  \in X$ such that
$d(x_i,x_j)> \rho/10$ for $i\neq j$ and $\sup_{x\in X}\inf_i d(x,x_i)
\le \rho/10$. Such sequence exists because $X$ is separable.
Second we let $B_i=B(x_i, \rho)$ and define $\widetilde{B_i}$ by the formula
$$\widetilde{B_i}=\bar{B}\left(x_i,\frac{\rho}{10}\right)\setminus
\bigcup_{j<i}\bar{B}\left(x_j,\frac{\rho}{10}\right),$$
where $\bar{B}\left(x, \rho\right)=\{y\in X \colon d(x,y)
\le \rho\}$. Third we put $\chi_i=\chi_{\widetilde B_i}$, where
$\chi_{\widetilde B_i}$ is the characteristic function of the set
${\widetilde B_i}$. Note that for $i\neq j$
 $B(x_i, \frac{\rho}{20}) \cap B(x_j, \frac{\rho}{20})=\emptyset$. Hence
$$
K=\sup_i\#\{j:\;d(x_i,x_j)\le  2\rho\} \le
  \sup_x  {V(x, (2+\frac{1}{20})\rho)\over
  V(x, \frac{\rho}{20})}< C  41^n< \infty.
$$
It is not difficult to see that
$$
\D_{\rho}  \subset \bigcup_{\{i,j:\, d(x_i,x_j)<
 2 \rho\}} \widetilde{B}_i\times \widetilde{B}_j \subset \D_{4 \rho}
$$
so
$$
Tf =\sum_{i,j:\, {d}(x_i,x_j)< 2\rho} P_{\widetilde B_i}T P_{\widetilde
B_j}f.
$$
Hence by H\"older's inequality
\begin{eqnarray*}
\|T f\|_{p}^p=\|\sum_{i,j:\, {d}(x_i,x_j)< 2\rho} P_{\widetilde B_i}T
P_{\widetilde B_j}f\|_{p}^p =\sum_i \|\sum_{j:\,{d}(x_i,x_j)<
2\rho} P_{\widetilde B_i}TP_{\widetilde B_j}f\|_{p}^p   \\
\le CK^{p-1}  \sum_i \sum_{j:\,{d}(x_i,x_j)< 2\rho} \| P_{\widetilde
B_i}TP_{\widetilde B_j}f\|_{p}^p
\\
\le CK^{p-1}  \sum_i   \sum_{j:\,{d}(x_i,x_j)< 2\rho}
\mu(\widetilde{B}_i)^{p({1\over p}-{1\over s})} \|P_{\widetilde
B_i}TP_{\widetilde B_j}f\|_{s}^p
\\ \le CK^{p}   \sum_j     \mu( {B}_j)^{p({1\over p}-{1\over
s})} \|TP_{ \widetilde B_j}f\|_{s}^p
\\\le CK^{p}   \sum_j \mu( {B}_j)^{p({1\over p}-{1\over
s})}\|TP_{\widetilde B_j}\|_{p\to s}^p\|P_{\widetilde
B_j}f\|_{p}^p
\\ \le CK^{p}  \sup_{x\in X}\big\{V(x,\rho)^{p({1\over p}-{1\over
s})} \|TP_{B(x,\rho)}\|^p_{p\to s}\big\}  \sum_j\|P_{\widetilde
B_j}f\|_{p}^p
\\=  CK^{p}  \sup_{x\in X}\big\{V(x,\rho)^{p({1\over p}-{1\over
s})} \|TP_{B(x,\rho)}\|^p_{p\to s} \big\} \|f\|_{p}^p.
\end{eqnarray*}
This finishes the proof of Lemma~\ref{le3.6}. \hfill{}
$\Box$

\medskip

\noindent
{\bf Proof of Theorem~\ref{th3.1}.}
Let   $\eta \in C_c^{\infty}(\mathbb R) $   be even and  such that $\supp \eta\subseteq \{ \xi: 1/4\leq |\xi|\leq 1\}$ and
$$
\sum_{\ell\in \ZZ} \eta(2^{-\ell} \lambda)=1 \ \ \quad \forall
{\lambda>0}.
$$
Then we  set $\eta_0(\lambda)= 1-\sum_{\ell> 0} \eta(2^{-\ell} \lambda)$,
\begin{eqnarray}\label{eq3.7}
F^{(0)}(\lambda)=\frac{1}{2\pi}\int_{-\infty}^{+\infty}
 \eta_0(t) \hat{F}(t) \cos(t\lambda) \;dt
\end{eqnarray}
and
\begin{eqnarray}\label{eq3.8}
F^{(\ell)}(\lambda) =\frac{1}{2\pi}\int_{-\infty}^{+\infty}
 \eta(2^{-\ell}t) \hat{F}(t) \cos(t\lambda) \;dt.
\end{eqnarray}
Note that in virtue of the Fourier inversion formula
$$
F(\lambda)=\sum_{\ell \ge 0}F^{(\ell)}(\lambda)
$$
and by Lemma \ref{le2.1}
$$
\supp K_{F^{(\ell)}(t\SL)} \subset \D_{2^{\ell}t}.
$$
Now by Lemma \ref{le3.6}
\begin{eqnarray}\label{eq3.9}
\big\|F(t\SL)\big\|_{p\to p } &\le& \sum_{\ell \ge 0}\big\|F^{(\ell)}(t\sqrt {L}) \big\|_{p\to p }\nonumber\\
&\le& \sum_{\ell \ge 0}    \sup_{x\in X}\big\{V(x,2^{\ell}t)^{{1\over
p}-{1\over s}}\big\|F^{(\ell)}(t\sqrt {L})P_{B(x,2^{\ell}t)}\big\|_{p\to s }\big\}.
\end{eqnarray}
Since $F^{(\ell)}$ is not compactly
supported we choose a function $\psi \in C_c^\infty(-4, 4)$ such that $\psi(\lambda)=1$ for $\lambda \in (-2,2)$
and note that
\begin{eqnarray} \label{eq3.10}
&&\hspace{-1.5cm}\big\|F^{(\ell)}(t\sqrt {L})P_{B(x,2^{\ell}t)}\big\|_{p\to s} \nonumber\\
& \le&  \big\|\big(\psi
F^{(\ell)}\big)(t\sqrt {L})P_{B(x,2^{\ell}t)}\big\|_{p\to s }
 +\big\|\big((1-\psi)F^{(\ell)}\big)(t\sqrt {L})P_{B(x,2^{\ell}t)}\big\|_{p\to
s }.
\end{eqnarray}
To estimate  the norm $\|\big(\psi
F^{(\ell)}\big)(t\sqrt {L})P_{B(x,2^{\ell}t)}\|_{p\to s }$, we use
condition ${\rm (ST^{q}_{p, s})}$ and the fact that $\psi \in
C_c(-4,4)$  to obtain
$$
\big\|\big(\psi F^{(\ell)}\big)(t\sqrt {L})P_{B(x,2^{\ell}t)}\big\|_{p\to s }
\le CV(x,2^{\ell}t)^{{1\over s}-{1\over p}}2^{\ell n({1\over
p}-{1\over s})} \big\| \delta_{4t^{-1}}\big(\psi
F^{(\ell)}\big) (t\cdot)\big\|_{q}
$$
for all $t>0$. Hence
 \begin{eqnarray}\label{eq3.11}
\sum_{\ell \ge 0}    \sup_{x\in X}\big\{V(x,2^{\ell}t)^{{1\over p}-{1\over
s}}\big\|\big(\psi F^{(\ell)}\big)(t\sqrt {L})P_{B(x,2^{\ell}t)}\big\|_{p\to
s}\big\} &\le& C\sum_{\ell \ge 0}   2^{\ell n({1\over p}-{1\over s})}\big\|  \delta_{4t^{-1}}
\big(\psi F^{(\ell)}\big) (t\cdot)\big\|_{q}\nonumber\\
&\leq& C \sum_{\ell \ge 0}   2^{\ell n({1\over p}-{1\over s})}\| F^{(\ell)} \|_{q}
\\
 &=&C\|F\|_{B_{q,\, 1}^{n({1\over p}-{1\over s})}},\nonumber
\end{eqnarray}
where   the last equality follows from definition of Besov space.
See, e.g., \cite[Chap.~VI ]{BL}.
Recall also that if $\beta>n(1/p-1/s)$ then $W^{\beta,q}\subseteq
B_{q, \, 1}^{n({1/p}-{1/s})}$ and $\|F\|_{B_{q,\, 1}^{n(1/p-1/s)}}\le C_\beta
\|F\|_{W^{\beta,q}}$, see again
\cite{BL}. Hence the forgoing estimates give
\begin{eqnarray}\label{eq3.12}
\sum_{\ell \ge 0}    \sup_{x\in X}\big\{ V(x,2^{\ell}t)^{{1\over p}-{1\over
s}}\big\|\big(\psi F^{(\ell)}\big)(t\sqrt {L})P_{B(x,2^{\ell}t)}\big\|_{p\to
s} \big\}  &\le& C\|F\|_{W^{\beta,q}}.
\end{eqnarray}

Next we show bounds for
   $\big\|\big((1-\psi)F^{(\ell)}\big)(t\sqrt
{L})P_{B(x,2^{\ell}t)}\big\|_{p\to s }$.
Since the function $1-\psi$ is supported outside the
interval $(-2,2)$, we can choose  a function $\phi\in C_c^{\infty}(2,8)$  such that
$$
1=\psi(\l)+\sum_{k\geq 0}\phi(2^{-k}\l)=\psi(\l)+\sum_{k\geq
0}\phi_{k}(\l)\quad \quad \forall \l>0.
$$
Hence
$$
\big((1-\psi)F^{(\ell)}\big)(\lambda)=  \sum_{k\geq 0}
 \big(\phi_{k}F^{(\ell)}\big)(\lambda) \quad \quad \forall \l>0.
$$
It follows from the implication  ${\rm (ST^{q}_{p, s})}\Rightarrow {\rm (ST^{\infty}_{p, s})}$ that
\begin{eqnarray*}
\big\|\big((1-\psi)F^{(\ell)}\big)(t\sqrt {L})P_{B(x,2^{\ell}t)}\big\|_{p\to s }
 &\le&   \sum_{k\geq 0}
\big\|\big(\phi_{k}F^{(\ell)}\big)(t\sqrt {L})P_{B(x,2^{\ell}t)}\big\|_{p\to s}
\\
 &\le&   C \sum_{k\geq 0} V(x,2^\ell t)^{{1\over s}-{1\over p}}2^{n(\ell+k)({1\over p}-{1\over s})}
 \big\|\delta_{ {2^{k+3}t^{-1}} }\big(\phi_kF^{(\ell)})(t \cdot)\big\|_{\infty}.
\end{eqnarray*}
Note that
$\supp F\in [-1, 1]$,  $\supp \phi\subset [2, 8]$ and $\check{\eta}$ is in the Schwartz class so
\begin{eqnarray*}
 \big\|\phi_kF^{(\ell)}\big\|_{\infty}= 2^{l}\big\|\phi_k(F * \delta_{2^{l}}\check{\eta} )   \big\|_{\infty}
\leq C 2^{-M(\ell+k)}\|F\|_{q}
\end{eqnarray*}
and similarly, $ \big\|\phi_kF^{(0)}\big\|_{\infty}\leq C 2^{-Mk}\|F\|_{q}$.
 Therefore
\begin{eqnarray}\label{eq3.13}
 \big\|\big((1-\psi)F^{(\ell)}\big)(t\sqrt {L})P_{B(x,2^{\ell}t)}\big\|_{p\to s }
 &\leq& C \sum_{k\geq 0} V(x,2^\ell t)^{{1\over s}-{1\over p}}2^{n(\ell+k)({1\over p}-{1\over s})}
2^{-M(\ell+k)}\|F\|_{q}\nonumber\\
 &\leq& C  V(x,2^\ell t)^{{1\over s}-{1\over p}}2^{\ell (n({1\over p}-{1\over s})-M)}
 \|F\|_{q}.
\end{eqnarray}
Hence
\begin{eqnarray}\label{eq3.14}\hspace{1cm}
 \sum_{\ell \ge 0}    \sup_{x\in X} \big\{ V(x,2^{\ell}t)^{{1\over
p}-{1\over s}}\big\|\big((1-\psi)F^{(\ell)}\big)(t\sqrt
{L})P_{B(x,2^{\ell}t)}\big\|_{p\to s } \big\}&\le&
 C \sum_{\ell \ge 0}   2^{\ell  (n({1\over p}-{1\over
s})-M )}\|F\|_{q } \nonumber\\
 &\le&  C \|F\|_{q }.
\end{eqnarray}
Now estimate ({\ref{eq3.1}}) follows from \eqref{eq3.14},  (\ref{eq3.9}), (\ref{eq3.10}) and
(\ref{eq3.12}). This completes the proof of Theorem~\ref{th3.1}. \hfill{} $\Box$

\medskip

\noindent {\bf  Proof of Theorem~\ref{th3.2}.}

\smallskip
\noindent
{\it Case} (1).  $t\geq 1/4$.

\smallskip
If  $t\geq 1/4$ then  $\supp \delta_t F
 \subset [-4,4]$. By ${\rm (AB_p)}$,
 \begin{eqnarray*}
\|F(t\SL)\|_{p \to p} &\leq&C 4^{\kappa n({1\over p}-{1\over
s})+\varepsilon}
 \|\delta_4(F(t\cdot))\|_{4^\kappa,q}
\leq C  \|F\|_{\infty}.
\end{eqnarray*}
 Recall that if  $\beta>1/q$,
  then $W^{\beta,q}(\RR) \subseteq L^{\infty}(\RR)\cap C(\RR)$ and $\|F\|_{\infty}\leq C\|F\|_{W^{\beta,q}}$.
Hence
$$
\sup_{t\geq 1/4}\|F(t\SL)\|_{p\to p}\leq C\|F\|_{\infty}\leq
C\|F\|_{W^{\beta,q}}.
$$

\smallskip
\noindent
{\it Case} (2).  $t\leq 1/4$.

\smallskip
Let $\xi\in C_c^\infty$ be an even function   such that $\supp
\xi\subset [-1,1],\, \hat{\xi}(0)=1$ and $\hat{\xi}^{(k)}(0)=0$ for
all $1\leq k\leq [\beta]+2$. Write
$\xi_{N^{\kappa-1}}=N^{\kappa-1}\xi(N^{\kappa-1}\cdot)$ where
$N=8[t^{-1}] +1$ and $[t^{-1}]$ denotes the integer part of
$t^{-1}.$ Following \cite{CowS} we write
$$
F(t\SL)=\big(\delta_tF - \xi_{N^{\kappa-1}}\ast
\delta_tF\big)(\SL)+(\xi_{N^{\kappa-1}}\ast \delta_tF)(\SL).
$$
We first prove that
\begin{eqnarray}
\|\big(\delta_tF - \xi_{N^{\kappa-1}}\ast
\delta_tF\big)(\SL)\|_{p\to p}&\leq& C \|F\|_{W^{\beta,q}}.\label
{eq3.15}
\end{eqnarray}

Observe that $\supp  (\delta_tF - \xi_{N^{\kappa-1}}\ast \delta_tF
)\subseteq [-N, N]$. We apply  ${\rm (AB_p)}$   to obtain
\begin{eqnarray}
 \|\big(\delta_tF - \xi_{N^{\kappa-1}}\ast \delta_tF\big)(\SL)\|_{p \to p}
  \leq C   N^{\kappa n({1\over
p}-{1\over s})+\varepsilon}\big\|\delta_{N}\big(\delta_tF -
\xi_{N^{\kappa-1}}\ast \delta_tF\big)\big\|_{N^\kappa,q}.
\label{ineq3.1}
\end{eqnarray}
Everything then boils down to estimate $\|\cdot\|_{N^\kappa,q}$
norm of $\delta_{N}\big(\delta_tF - \xi_{N^{\kappa-1}}\ast
\delta_tF\big).$ We make the following claim.
  For its proof, see \cite[(3.29)]{CowS} or \cite[Propostion 4.6]{DOS}.

\begin{lemma}\label{le3.7}
Suppose that\ \ $\xi\in C_c^\infty$ is an even function such that
$\supp \xi\subset [-1,1],  \ \ \hat{\xi}(0)=1$ and
$\hat{\xi}^{(k)}(0)=0$ for all $1\leq k\leq [\beta]+2$. Next assume that
 $\supp H\subset [-1,1]$. Then
\begin{eqnarray}\label {eq3.17}
\|H-\xi_N\ast H\|_{N,q}\leq CN^{-\beta}\|H\|_{W^{\beta,q}}
\end{eqnarray}
for all $\beta>1/q$ and any $N\in \NN$.
\end{lemma}

Note that  $\delta_N\big(\delta_tF - \xi_{N^{\kappa-1}}\ast \delta_tF\big)=\delta_{Nt}F - \xi_{N^{\kappa}}\ast \delta_{Nt}F$.
Now, if $\beta>\max \{{n({1/p}-{1/s})}, 1/q\}$
then \eqref{eq3.15} follows from Lemma~\ref{le3.7} and estimate \eqref{ineq3.1}.

\medskip

It remains to show that
\begin{eqnarray}
\|(\xi_{N^{\kappa-1}}\ast \delta_tF)(\SL)\|_{p\to p}&\leq& C
\|F\|_{W^{\beta,q}}. \label {eq3.16}
\end{eqnarray}

Let $ F^{(\ell)} $ be functions defined in
(\ref{eq3.7}) and (\ref{eq3.8}). Following the proof of
Theorem~\ref{th3.1}, we write
$$
(\xi_{N^{\kappa-1}}\ast \delta_tF)(\lambda)=\sum_{\ell \ge 0}
\big(\xi_{N^{\kappa-1}}\ast \delta_tF^{(\ell)}\big)(\lambda),
$$ and by Lemma \ref{le2.1},
 $
\supp K_{  (\xi_{N^{\kappa-1}}\ast \delta_tF^{(\ell)} )(\SL)}
\subset \D_{2^{\ell}t}. $
 Now by Lemma \ref{le3.6}
\begin{eqnarray}\label{eq3.18}
\big\|\big(\xi_{N^{\kappa-1}}\ast \delta_tF \big)(\SL)\big\|_{p\to p
} &\le&
\sum_{\ell \ge 0}\big\|\big(\xi_{N^{\kappa-1}}\ast \delta_tF^{(\ell)}\big) (\SL) \big\|_{p\to p }\nonumber\\
&\le& \sum_{\ell \ge 0}    \sup_{x\in X} \big\{
V(x,2^{\ell}t)^{{1\over p}-{1\over
s}}\big\|\big(\xi_{N^{\kappa-1}}\ast \delta_tF^{(\ell)}\big)
(\SL)P_{B(x,2^{\ell}t)}\big\|_{p\to s } \big\}.
\end{eqnarray}
Take a function $\psi \in C_c^\infty(-4, 4)$ such that $\psi(\lambda)=1$ for $\lambda \in (-2,2)$. Then
\begin{eqnarray} \label{eq3.19}
 \big\|\big(\xi_{N^{\kappa-1}}\ast \delta_tF^{(\ell)}\big) (\SL)P_{B(x,2^{\ell}t)}\big\|_{p\to s}
& \le&  \big\|\big(\delta_t\psi  (\xi_{N^{\kappa-1}}\ast
\delta_tF^{(\ell)})\big) (\SL)
 P_{B(x,2^{\ell}t)}\big\|_{p\to s }\nonumber\\
 &+&\big\|\big((1-\delta_t\psi)  (\xi_{N^{\kappa-1}}\ast \delta_tF^{(\ell)})\big) (\SL)
 P_{B(x,2^{\ell}t)}\big\|_{p\to s } ={ I_{\ell}} +{ I\!I_{\ell}}.
\end{eqnarray}
 Note that ${\rm (SC^{q,\kappa}_{p,
s})}\Rightarrow {\rm (SC^{\infty}_{p, s})}$ and $t\leq 1/4$. Using ${\rm (SC^{\infty}_{p, s})}$
instead of ${\rm (ST^{\infty}_{p, s})}$, we show  as in the proof
of (\ref{eq3.13})
 that
 $I\!I_{\ell}\leq C  V(x,2^\ell t)^{{1\over s}-{1\over p}} 2^{\ell  (n({1\over p}-{1\over
s})-M )}\|F\|_{q}
$
for some large $M> n({1/p}-{1/s})+1$.

Next we estimate  the term $I_{\ell}$. We assume that  $\psi \in
C_c(-4,4)$ so by  ${\rm (SC^{q,\kappa}_{p, s})}$
$$
I_{\ell}\le CV(x,2^{\ell}t)^{{1\over s}-{1\over p}}2^{\ell n({1\over
p}-{1\over s})} \big\| \delta_{N}\big(\delta_t\psi
(\xi_{N^{\kappa-1}}\ast \delta_tF^{(\ell)})\big)\big\|_{N^\kappa,q}.
$$
Observe that (see also \cite[(3.19)]{CowS})
$$
 \big|(\xi\ast \delta_tF^{(\ell)})(\lambda)\big|\leq
\|\xi\|_{L^{q'}}\Big(\int_{\l-1}^{\l+1}|F^{(\ell)}(tu)|^qdu\Big)^{1/q}
$$
so
\begin{eqnarray*}
\big\|\delta_{N}\big(\delta_t\psi  (\xi_{N^{\kappa-1}}\ast
\delta_tF^{(\ell)})\big)\big\|_{N^\kappa,q}&=&\Big({1\over
2N^\kappa}\sum_{i=1-N^\kappa}^{N^\kappa}\sup_{\l\in[{i-1\over
N^\kappa},{i\over N^\kappa} )}|\psi(tN\l)
(\xi_{N^{\kappa-1}}\ast \delta_tF^{(\ell)})(N\l)|^q\Big)^{1/q}\\
&\leq&C   \Big({1\over
N^\kappa}\sum_{i=1-N^\kappa}^{N^\kappa}\sup_{\l\in[{i-1},{i} )}
 \big|(\xi\ast \delta_{tN^{1-\kappa}}F^{(\ell)})(\lambda)|^q\Big)^{1/q}\\
&\leq&C\|\xi\|_{q'} \Big({1\over
N^\kappa}\sum_{i=1-N^\kappa}^{N^\kappa}\sup_{\l\in[{i-1},{i} )}\int_{\l-1}^{\l+1}|F^{(\ell)}(tN^{1-\kappa}u)|^qdu\Big)^{1/q}\\
&\leq&C\Big({1\over
N^\kappa}\sum_{i=1-N^\kappa}^{N^\kappa}\int_{i-2}^{i+1}|F^{(\ell)}(tN^{1-\kappa}u)|^qdu\Big)^{1/q}\\
&\leq&C\Big({1\over
Nt}\int_{-\infty}^{\infty}|F^{(\ell)}(u)|^qdu\Big)^{1/q}
 \leq  C \|F^{(\ell)}\|_{q}.
\end{eqnarray*}
This shows that $
I_{\ell}\le CV(x,2^{\ell}t)^{{1\over
s}-{1\over p}}2^{\ell n({1\over p}-{1\over s})}\|F^{(\ell)}\|_{q}$.
Using the above estimates of $I_{\ell}$ and $I\!I_{\ell}$, together  with (\ref{eq3.18}) and (\ref{eq3.19}),
we can argue as in  (\ref{eq3.12}) and
(\ref{eq3.14}) to obtain estimate
({\ref{eq3.16}}).
This proves Theorem~\ref{th3.2}. \hfill{} $\Box$

\medskip

Now we discuss another condition introduced by C.D. Sogge (see  \cite{Sog1, Sog3, Sog4}).
 We say that $L$ satisfies  {\it   $(p,p')$ spectral cluster estimate} ${\rm ( S_p )}$
 for some $1\leq p<2$ and its conjugate $p'$ if
 the spectral projection
  $E_{\sqrt{L}}[k,k+1)$ maps  $L^p(X)$ to $L^{p'}(X)$
  and
  $$
\big\|E_{\sqrt{L}}[k, k+1)\big\|_{p\to p'} \leq
C (1+k)^{n({1\over p}-{1\over p'})-1} \leqno{\rm (S_p)}
$$
  for all $k\geq 0$.

\begin{proposition}\label{prop3.8}
Suppose that $1\leq p\leq 2n/(n+1)$, $\mu(X)<\infty$ and $ V(x,r) \leq C \min\big(r^n, 1\big)$ for every $x\in X$ and $r>0$.
 Then  conditions ${\rm (S_p)}$  and ${\rm (SC^{2,1}_{p, 2})}$
are
 equivalent.
\end{proposition}

 \noindent
 {\bf Proof.}  We first prove the implication ${\rm (S_p)}\Rightarrow
 {\rm (SC^{2,1}_{p, 2})}$.  Note that for every even Borel function $F$ such that
    $\supp F\subset [-N, N]$,
 \begin{eqnarray*}
 \big\|F(\SL) f\big\|^2_{ 2}
&= &\sum_{k=0}^{N} \big\langle E_{\sqrt{L}}[k, k+1)F(\SL)f,
\, E_{\sqrt{L}}[k, k+1)F(\SL)f\big\rangle \\
&\leq& \sum_{k=0}^{N} \big\|E_{\sqrt{L}}[k,
k+1)F(\SL)\big\|^2_{p\to 2} \|f\|_p^2.
\end{eqnarray*}
Using a $T^{\ast} T$ argument
and condition  ${\rm (S_p)}$  we obtain
\begin{eqnarray*}
 \big\|F(\SL)P_{B(x, r)}\big\|^2_{p\to 2}
 &\leq& C\big\|F(\SL)\big\|^2_{p\to 2} \\
 &\leq& C\sum_{k=0}^{N} \big\|E_{\sqrt{L}}[k, k+1)F(\sqrt{L})\big\|^2_{p\to 2}\\
 &\leq& C\sum_{k=0}^{N} \sup_{\lambda\in [k, k+1)}
\big|F(\lambda)\big|^2 (1+k)^{n({1\over p}-{1\over p'})-{1}}
 \\
 &\leq& CN^{n({1\over p}-{1\over p'})}{1\over N}\sum_{k=0}^{N} \sup_{\lambda\in [k, k+1)}
\big|F(\lambda)\big|^2
 \\
&\leq& CV(x,r)^{2({1\over 2}-{1\over p})}(Nr)^{2n({{1\over p}-{1\over 2}})}
\big\|\delta_NF\big\|_{N,2}^2,
\end{eqnarray*}
and hence condition ${\rm (SC^{2,1}_{p, 2})}$ is satisfied.

\medskip

Next we  prove the implication ${\rm ({SC}_{p, 2}^{2,1})}\Rightarrow
{\rm ({S}_p)}$. By ${\rm ({SC}_{p, 2}^{2,1})}$
 \begin{eqnarray}\label{eq3.21}
  \big\|\chi_{[k,k+1)}(\SL) \big\|_{p\to 2} \leq
 C(1+k)^{n({1\over p}-{1\over 2})}\|\delta_{(1+k)}\, \chi_{[k,k+1)}\|_{1+k,\, 2}.
  \end{eqnarray}
Hence
  \begin{eqnarray*}
  \big\|E_{\SL}[k,k+1) \big\|_{p\to p'}&=&  \big\|E_{\SL}[k,k+1) \big\|^2_{p\to 2}\\
  &\leq& C(1+k)^{2n({1\over p}-{1\over 2})} \big\|\delta_{(1+k)}\, \chi_{[k,k+1)} \big\|_{1+k,\, 2}^2\\
  &\leq& C(k+1)^{n({1\over p}-{1\over p'})-1},
  \end{eqnarray*}
  which shows  ${\rm (S_p)}.$
\hfill{} $\Box$

 \medskip

\begin{proposition}\label{prop3.9}
Assume that $\mu(X)<\infty$.
 Then  conditions ${\rm (SC^{\infty}_{p, 2})}$  and ${\rm (ST^{\infty}_{p, 2})}$
are
 equivalent.
\end{proposition}

 \noindent
 {\bf Proof.}  Since  $\mu(X)<\infty$, we may assume that $X=B(x_0, 1)$ for some $x_0\in X$.
   Observing  that   $\|F\|_{N, \infty}=\|F\|_{\infty}$,
 we have the implication
 ${\rm (ST_{p, 2}^\infty})\Rightarrow {\rm (SC_{p, 2}^\infty)}$. Let us
 show the implication   $ {\rm (SC_{p, 2}^\infty)} \Rightarrow{\rm (ST_{p, 2}^\infty)}$.
  Assume that  supp $F\subseteq [0,R]$. If
$R\geq 1$, then we   let $N=[R]+1$   and
   $ {\rm (ST_{p, 2}^\infty)} $  follows readily.
Now for $0<R<1$,  from  condition ${\rm (SC_{p, 2}^\infty})$
we can take  $N=1$ and the  ball $B(x_0, 1)$
to  obtain
$$
\|F(\SL)\|_{p\to 2}\leq C\|F\|_{\infty}.
$$
 Now for any $x\in X$ and $r>0$, we note that conditions  $1/2-1/p< 0$ and $\mu(X)<\infty$
give that $V(x,r)^{1/2-1/p}\geq C$.
Hence for any $r\geq 1/R,$
$$
\big\|F(\SL)P_{B(x,r)}\big\|_{p\to 2}\leq C\|F\|_{\infty}\leq
CV(x,r)^{{1\over 2}-{1\over p}}(rR)^{n({1\over p}-{1\over
2})}\|\delta_R F\|_{\infty},
$$
this is  ${\rm (ST_{p, 2}^\infty)}.$ The proof of
Proposition~\ref{prop3.9} is finished. \hfill{} $\Box$

\section{Sharp spectral multipliers - singular integral case}\label{sec4}
\setcounter{equation}{0}

\noindent{\bf 4.1.\, Statements.}\,
As in Section~\ref{sec3}
we discuss two type of results corresponding to estimates ${\rm (ST^{q}_{p, 2})}$ or
 ${\rm (SC^{q, \kappa}_{p, 2})}$.
The aim  of this section is to prove  singular integral
 versions of  Theorems \ref{th3.1} and \ref{th3.2}.  We use the same
 assumptions and notation as in Section \ref{sec3}.  Recall that $n$ is a homogeneous
 dimension from  (\ref{eq2.2}). Fix a non-trivial auxiliary function $\eta \in C_c^\infty(0,\infty)$.

 \begin{theorem}\label{th31.1}
  Assume that operator $L$ satisfies   property {\rm (FS)}
and condition  ${\rm (ST^{q}_{p, 2})}$ for some $p,q$ satisfying
$1\leq p<2$ and $1\leq q\leq \infty$. Then for any even bounded Borel
function $F$ such that
$\sup_{t>0}\|\eta\, \delta_tF\|_{W^{\beta, q}}<\infty $ for some
$\beta>\max\{n(1/p-1/2),1/q\}$ the operator
$F(\SL)$ is bounded on $L^r(X)$ for all $p<r<p'$.
In addition,
\begin{eqnarray*}
   \|F(\SL)  \|_{r\to r}\leq    C_\beta\Big(\sup_{t>0}\|\eta\, \delta_tF\|_{W^{\beta, q}}
   + |F(0)|\Big).
\end{eqnarray*}
\end{theorem}

Note that if $q < \infty$ then condition ${\rm (ST^{q}_{p, 2})}$ implies that
$E_L(\{0\})=0$ (and in fact $E_L(\{\lambda\})=0 $) so the term $F(0)$ can be omitted in the above statement.

The next theorem
 is a variation of Theorem~\ref{th31.1} suitable for the operators satisfying
 condition ${\rm (SC^{q,\kappa}_{p, 2})}$. It is a singular integral version of  Theorem \ref{th3.2} above.

\begin{theorem}\label{th31.2}  Suppose the operator $L$ satisfies   property {\rm (FS)},
 conditions   ${\rm (E_{p,2})}$ and   ${\rm (SC^{q,\kappa}_{p, 2})}$ for some $p,q$ such that $1\leq p<2$
and $1\leq q\leq \infty$, and  a fixed natural number $\kappa$. In
addition, we assume that for any $\varepsilon>0$ there exists a
constant $C_\varepsilon$ such that for all $N\in \NN$ and
all even Borel functions $F$ such that supp $F\subset [-N,N]$,
$$
\|F(\SL)\|_{p\to p}\leq C_\varepsilon N^{\kappa n({1\over p}-{1\over
2})+\varepsilon}\|\delta_N F\|_{N^{\kappa},q}. \leqno{\rm (AB_p)}
$$
  Then for any even bounded Borel function $F$
 such that
$\sup_{t>0}\|\eta\, \delta_tF\|_{W^{\beta, q}}<\infty $ for some
$\beta>\max\{n(1/p-1/2),1/q\}$  the operator $F(\SL)$ is bounded on
$L^r(X)$ for all $p<r<p'$. In addition,
\begin{eqnarray*}
   \|F(\SL)  \|_{r\to r}\leq    C_\beta\Big(\sup_{t>0}\|\eta\, \delta_tF\|_{W^{\beta, q}}
   + |F(0)|\Big).
\end{eqnarray*}
\end{theorem}

\medskip

\begin{remark}\label{4.4444}
Suppose that $\mu(X)<\infty$   and ${\rm (SC^{q, \kappa}_{p, 2})}$ holds for some $\kappa\geq 1$. Then
${\rm (SC^{\infty}_{p, 2})}$ and ${\rm (E_{p,2})}$ are satisfied by  Remark~\ref{re3.55} and  Proposition~\ref{prop3.9}.
In addition,   ${\rm (AB_p)}$ holds
by Proposition \ref{le3.3333}.  Therefore,   Theorem~\ref{th31.2} holds in this case without assumptions
${\rm (E_{p,2})}$ and ${\rm (AB_{p})}$.
\end{remark}

\medskip

By a classical dyadic decomposition of $F$, we can write
$F(\sqrt{L})$ as the sum $\sum F_j(\sqrt{L})$. Then we  apply Theorems \ref{th3.1} and \ref{th3.2} to estimate  $\| F_j(\sqrt{L})\|_{r \to r}$.
However, as mentioned in the introduction, this does not automatically imply that the
operator $F(\sqrt{L})$ acts boundedly on $L^r$. See \cite{Carb1, S3,SS} where this problem is discussed in the Euclidean case.
Our proof is almost identical to one in  \cite{Carb1}. Nevertheless we give full details because the changes which are required
to adapt the arguments to the general setting are not trivial.

Note that condition   ${\rm (ST^{q}_{p, 2})}$
implies ${\rm (E_{p,2})}$ (see Proposition~\ref{prop2.3}). Therefore
Theorems~\ref{th31.1}~and~\ref{th31.2} follow from Theorems~\ref{th3.1}~and~
\ref{th3.2} (with $s=2$) and the next result.

 \begin{theorem}\label {th33.1}
Assume that $L$ satisfies  the
 finite speed propagation property {\rm (FS)}
and condition  ${\rm (E_{p_0,2})}$ for some $1\leq p_0<2$. Next
assume that for all even Borel   functions  $F$ such that  $\supp F$ is
compact and $\|F\|_{W^{\beta,q}}< \infty$ for some $\beta $ and $ q$
satisfying $\beta>\max\{n(1/p_0-1/2),1/q\}$ and $1\leq q\leq \infty$,
\begin{eqnarray}\label{e33.1}
\sup_{t>0}\|F(t\SL)\|_{p\to p}\leq C\|F\|_{W^{\beta,q}}, \ \ \ p_0\leq p\leq p'_0.
\end{eqnarray}
  Then for any bounded Borel function $F$  such that
\begin{eqnarray}\label{e333.1}
\sup_{t>0}\|\eta\delta_tF\|_{W^{\beta,q}} < \infty
\end{eqnarray}
for some $\beta>\max\{n(1/p_0-1/2),1/q\}$,  the operator
$F(\sqrt{L})$ is bounded on $L^r(X)$ for all $p_0<r<p_0'$.
\end{theorem}

\medskip

\noindent{\bf 4.2.\, Singular integrals.}\,
This  subsection is devoted to the  proof of Theorem~\ref{th33.1}.
We start with the following lemma.

\begin{lemma}\label{le31.3}  Suppose  that
 operator $L$ satisfies property {\rm (FS)}
and condition  ${\rm (E_{p_0,2})}$ for some $1\leq p_0< 2$.
\begin{itemize}
\item[(a)]
Assume in addition that $F$ is  an even bounded Borel function  such that
$$
\sup_{t>0}\|\eta\delta_tF\|_{C^k}<\infty
$$
for some integer $k> n/2 + 1$ and some non-trivial function $\eta\in C_c^{\infty}(0, \infty)$.
Then the operator  $F(\SL)$ is bounded on $L^{p}(X)$ for all $p_0<p<p_0'$. \\

\item[(b)] Assume in addition that $\psi $  be  an even function in $
{\mathscr S} ({\mathbb{R}})$ such that $\psi(0)=0$. Define the
quadratic functional for $f\in L^2(X)$
\begin{eqnarray*}
{\mathcal G}_L(f)(x)=\Big( \sum_{j\in{\mathbb Z}}  |\psi(2^j\SL)f  |^2\Big)^{1/2}.
\end{eqnarray*}
Then ${\mathcal G}_L$ is bounded on $L^{p}(X)$ for all $p_0<p<p_0'.$
\end{itemize}
\end{lemma}

\noindent
{\bf Proof.}  The finite speed propagation property implies  $L^2-L^2$ off-diagonal estimate
$$ \| P_{X\setminus B(x,2^jr)} e^{-r^2L} P_{B(x,r)} \|_{2\to 2} \le C e^{-c2^{2j}},$$
see \cite{CouS, S2}. It follows from  ${\rm (E_{p_0,2})}$ that
$$  \| P_{X\setminus B(x,2^jr)} e^{-r^2L} P_{B(x,r)} \|_{p_0\to 2} \le C V(x,r)^{\frac{1}{2} - \frac{1}{p_0}}.
$$
Now the Riesz-Thorin interpolation theorem gives for $p \in (p_0,2)$ the following $L^p-L^2$ off-diagonal estimate
$$  \| P_{X\setminus B(x,2^jr)} e^{-r^2L} P_{B(x,r)}
\|_{p\to 2} \le C V(x,r)^{\frac{1}{2} - \frac{1}{p}} e^{-c'
2^{2j}}.$$
Now assertion (a) follows from \cite{Blu}.  The latter  off-diagonal  estimate implies that $L$ has a bounded holomorphic
functional calculus on $L^p$ for $p_0 < p < p_0'$ (see \cite{BlKu}).  It is known that the holomorphic functional calculus  implies
 the quadratic estimate of  assertion  (b) (see \cite{CDMY, Mc}).   \hfill{}$\Box$

\medskip

Throughout the rest of this section, $\Phi $  denotes   an even function, $  \Phi \in {\mathscr S} ({\mathbb{R}})$
such that $\Phi(0)=1$ and whose Fourier transform
$\hat\Phi$ is supported in $[-1,1]$. We take  $\eta \in C_c^{\infty}(\mathbb R) $   even and
such that  $\supp \eta\subseteq \{ \xi: 1/2\leq |\xi|\leq 2\}$ and
$
\sum_{\ell\in \ZZ} \eta(2^{-\ell} \lambda)=1$ for all $
{\lambda>0}.
$
Set
$\eta_\ell(\lambda)=\eta(2^\ell\l)$ and
$$
\widehat{F^{(\ell)}} =  \eta_{-\ell} \widehat{F}, \ \ \ \ \ell\in \ZZ.
$$
Put
$Q_\ell(\l)=\sum_{k\geq0}\eta_{k+\ell}(\l)$.
\begin{proposition}\label{th32.1}
Suppose that operator $L$ satisfies property {\rm (FS)}
and condition  ${\rm (E_{p_0,2})}$ for  $1\leq p_0<2$ and let $p_0<p<2$. Assume in addition that for
an even bounded Borel function $F$ the following estimates hold
\begin{eqnarray}\label {e32.1}
\sum_{k<0}\sup_j\big\|\sum_{\ell\geq
0}F^{(j+\ell)}(\SL)\eta_{j+k}(\SL)\big\|_{p\to p}<\infty
\end{eqnarray}
and
\begin{eqnarray}\label {e32.2}
\sup_j\big\|\sum_{\ell\geq
0}F^{(j+\ell)}(\SL)(I- {\Phi}(2^{j}\SL))Q_j(\SL)\big\|_{p\to
p}<\infty.
\end{eqnarray}
Then $F(\SL)$ is of weak-type $(p,p)$.
\end{proposition}

\noindent
{\bf Proof.}
  Let $f\in L^p(X)$ and $\alpha> \mu(X)^{-1/p}\|f\|_p$. A simple variation of
the Calder\'on-Zygmund decomposition of  $|f|^p$ at height $\alpha$ shows that there
exist  constants $C$ and  $K$ such  that $
f=g+b=g+\sum_{i}b_i $ so that  $\|g\|_p\leq C\|f\|_p,$
 $\|g\|_{\infty}\leq C\alpha$,
  each $b_i$  is supported on  ball $B_i$ of radius $2^{j(i)}$, and  $\# \{i: x\in 8B_i\} \leq K$ for all $x\in X$,
$\int_{X}\ |b_i|^pd\mu\leq C\alpha^p \mu(B_i), $  and
$\sum_{i}\mu(B_i) \leq C{\alpha}^{-p} \|f\|_p^p$.  As  a
consequence,  $\alpha^{p-2}\|g\|_2^2\leq C\|f\|_p^p$. As in \cite{Carb1}, we choose $2^{j(i)}$ rather than
$r(i)$ to be able to sum in $j$.

We define the ``nearly good" and ``very bad" functions ${\tilde g}$ and ${\tilde b}$ by
$$
{\tilde g}=g+ \sum_i2 {\Phi}(2^{j(i)}\SL)b_i-\sum_i {\Phi}^2(2^{j(i)}\SL) b_i \ \ \
 \ {\rm and}\ \ \  \ {\tilde b}=\sum_i(I-{\Phi}(2^{j(i)}\SL))^2b_i.
$$
By Lemma 2.1,     $\supp {\Phi}(2^{j(i)}\SL)b_i\subset 4B_i$, and
by the Calder\'on-Zygmund decomposition,  every  $ x \in X$ belongs to
no more than $K$ balls $4B_i$. Let $N>n(1/p-1/2)$. Now by condition
${\rm (E_{p_0,2})}$
\begin{eqnarray*}
\big\|\sum_i{\Phi}(2^{j(i)}\SL)b_i\big\|_2^2
&&\leq K \sum_i\big\|{\Phi}(2^{j(i)}\SL)b_i\big\|_2^2\nonumber\\
&& \leq
C \sum_i  \sup_\l\big|{\Phi}(2^{j(i)}\l)(1+2^{j(i)}\l)^{N}\big|^2
\big\|(1+2^{j(i)} \sqrt{L})^{-N}b_i\big\|_2^2\nonumber\\
&&\leq C \sum_i \mu(B_i)^{1-{2\over
p_0}}\|b_i\|_{p_0}^2\nonumber\\
&&\leq C \sum_i \mu(B_i)^{1-{2\over
p_0}}\|b_i\|_{p}^2\mu(B_i)^{{2\over p_0}-{2\over p}}\nonumber\\
&&\leq C\alpha^2\sum_i \mu(B_i) \leq C\alpha^{2-p }\|f\|_p^p. \nonumber
\end{eqnarray*}
Replacing $\Phi$ by $\Phi^2$ yields
 $\big\|\sum_i {\Phi}^2(2^{j(i)}\SL) b_i\big\|_2^2\leq C\alpha^{2-p }\|f\|_p^p.$
By the standard $L^2$ argument
$$
\mu\big(\{x:F(\SL)({\tilde g})(x)>\alpha\}\big)
  \leq  C\alpha^{-p}\|f\|_p^p.
$$
 It remains to treat
\begin{eqnarray}\label {e32.4}
 F(\SL)\Big(\sum_i(I-{\Phi}(2^{j(i)}\SL))^2
b_i\Big)
&=&\sum_i\sum_{\ell\geq0}F^{(j(i)+\ell)}(\SL)(I-{\Phi}(2^{j(i)}\SL))^2b_i\nonumber\\
&+&
\sum_i\sum_{\ell<0}F^{(j(i)+\ell)}(\SL)(I-{\Phi}(2^{j(i)}\SL))^2b_i.
\end{eqnarray}
By Lemma 2.1
$$
\supp
F^{(j(i)+\ell)}(\SL)(I-{\Phi}(2^{j(i)}\SL))^2b_i\subseteq B(x_i,
2^{j(i)+\ell+1}+3\cdot 2^{j(i)})\subseteq 8B_i, \ \ \ \forall \ell<0.
$$
Thus the second term of (\ref{e32.4}) is
supported in $\cup8B_i$ and
 $
\sum \mu(8B_i)\leq C\sum\mu(B_i)\leq C\alpha^{-p}\|f\|_p^p.
 $

To treat the first term we  show that
\begin{eqnarray}\label {e32.5}
\big\|\sum_j\sum_{\ell\geq0}F^{(j+\ell)}(\SL)(I-{\Phi}(2^{j}\SL))^2f_j\big\|_p\leq
C\big\|\sum_j|f_j|\, \big\|_{p}.
\end{eqnarray}
   If we apply (\ref{e32.5}) with
$$
f_j=\sum_{i:j(i)=j}b_i,
$$
we see that
\begin{eqnarray*}
 \big\|\sum_i\sum_{\ell\geq0}F^{(j(i)+\ell)}(\SL)(I-{\Phi}(2^{j(i)}\SL))^2b_i\big\|_p
&=&
\big\|\sum_j\sum_{\ell\geq0}F^{(j+\ell)}(\SL)(I-{\Phi}(2^{j}\SL))^2f_j\big\|_p\\
&\leq&C \big\|\sum_j|\sum_{i:j(i)=j}b_i|\big\|_p\\
&\leq& C\big\|\sum_i|b_i|\big\|_p\\
&\leq& C\big(\sum_i\big\|b_i\big\|_p^p\big)^{1/p} \leq C\|f\|_p,
\end{eqnarray*}
which completes the proof.

As in \cite{Carb1} we argue that by duality, (\ref{e32.5}) is equivalent to
\begin{eqnarray}\label {e32.6}
\big \|\sup_j|\sum_{\ell\geq0}\bar{F}^{(j+\ell)}(\SL)(I-{\Phi}(2^{j}\SL))^2h|\big \|_{p'}\leq
C\|h\|_{p'}.
\end{eqnarray}
Write
\begin{eqnarray}\label{e3332.6}
 \sum_{\ell\geq0}\bar{F}^{(j+\ell)}(\SL)(I-{\Phi}(2^{j}\SL))^2
&=& \sum_{k<0}\sum_{\ell\geq0}\bar{F}^{(j+\ell)}(\SL)(I-{\Phi}(2^{j}\SL))^2\eta_{j+k}(\SL)\nonumber\\
 &+&\sum_{\ell\geq0}\bar{F}^{(j+\ell)}(\SL)(I-{\Phi}(2^{j}\SL))^2Q_j(\SL).
\end{eqnarray}
Let  $\tilde{\eta} \in C_c^{\infty}(0, \infty)$ be a non-negative function
satisfying  $\supp \tilde{\eta}\subseteq   [1/4,4]$ and $\tilde{\eta}=1$ on $[1/2,2]$, and let
$\tilde{\eta}_j$ denote the function $\tilde{\eta}(2^j\cdot)$.
By  Lemma~\ref{le31.3} point (b) with $\psi=
\tilde{\eta} $ there exists a constant $C>0$ independent of
$k<0$ such that for $2<p'<p_0'$,
\begin{eqnarray*}
\Big(\sum_j\big\|    \tilde{\eta}_{j+k}(\SL)h \big\|_{p'}^{p'} \Big)^{1/p'} &\leq&
\Big\|\Big( \sum_j\big|  \tilde{\eta}_{j+k}(\SL)h \big|^2\Big)^{1/2}\Big\|_{p'}
 \leq  C\|h\|_{p'},
\end{eqnarray*}
and by Lemma~\ref{le31.3} point (a), it follows that for $2<p'<p_0'$, $\|(I-{\Phi}(2^{j}\SL))^2\|_{p'\to p'} \leq C$ for some constant $C>0$
independent of $j\in {\mathbb Z}$.
Hence
\begin{eqnarray}\label {e32.7}
&&\hspace{-0.8cm}\Big\|\sup_j\big|\sum_{k<0}\sum_{\ell\geq0}\bar{F}^{(j+\ell)}(\SL)
(I-{\Phi}(2^{j}\SL))^2\eta_{j+k}(\SL)h\big|\Big\|_{p'}\nonumber\\
&&\hspace{0.3cm} \leq \sum_{k<0}\Big(\sum_j\big\|  \sum_{\ell\geq0}\bar{F}^{(j+\ell)}(\SL)
(I-{\Phi}(2^{j}\SL))^2\eta_{j+k}(\SL)\tilde{\eta}_{j+k}(\SL)h\big\|_{p'}^{p'} \Big)^{1/p'}
\nonumber\\
&&\hspace{0.3cm}\leq
C\sum_{k<0}\sup_j\Big\{\|(I-{\Phi}(2^{j}\SL))^2\|_{p'\to p'} \Big\| \sum_{\ell\geq0}\bar{F}^{(j+\ell)}(\SL)\eta_{j+k}(\SL) \Big\|_{p'\to
p'}\Big\}\nonumber\\
&&\hspace{7cm} \times
\Big(\sum_j\big\|    \tilde{\eta}_{j+k}(\SL)h\big\|_{p'}^{p'} \Big)^{1/p'} \nonumber\\
&&\hspace{0.3cm}\leq
 C\sum_{k<0}\sup_j\big\|\sum_{\ell\geq0}\bar{F}^{(j+\ell)}(\SL)\eta_{j+k}(\SL)\big\|_{p'\to
p'}
 \|h\|_{p'}
 \leq
 C
 \|h\|_{p'}.
\end{eqnarray}
The last inequality follows from  assumption (\ref{e32.1}). Using assumption
\eqref{e32.2} instead of (\ref{e32.1})   the similar argument  as above gives the following estimate  
\begin{eqnarray*}
 \Big\|\sup_j\big|\sum_{\ell\geq0}\bar{F}^{(j+\ell)}(\SL)(I-{\Phi}(2^{j}\SL))^2Q_j(\SL)h\big|\Big\|_{p'}
 \leq
C\|h\|_{p'}.
\end{eqnarray*}
This shows  (\ref{e32.6}) and ends the
  proof of Proposition~\ref{th32.1}.
  \hfill{}$\Box$

\medskip

\begin{proposition}\label {th32.2}  Suppose that operator $L$ satisfies  property {\rm (FS)}
and condition  ${\rm (E_{p_0,2})}$ for some $1\leq p_0<2$ and fix $p \in (p_0, 2)$. Next
assume that for all even Borel functions  $H$ such that  $\supp H$ is compact
and $\|H\|_{W^{\beta, \infty}}< \infty$ for some $\beta>n/2,$
 \begin{eqnarray}\label{e32.8}
\sup_{t>0}\|H(t\SL)\|_{p\to p}\leq C\|H\|_{W^{\beta, \infty}}.
\end{eqnarray}
 Then for any even bounded   Borel function $F$  such that
  $\|(F\eta_i)^{(j)}(\SL)\|_{p\to p}\leq \alpha(i-j)$ for all $i, j \in \ZZ$ with
$$
\sum_{k\leq 0}(|k|+1)\alpha(k)<\infty,
$$
the operator $F(\SL)$ is of weak-type $(p,p)$.
\end{proposition}

\noindent
{\bf Proof.} By Proposition~\ref{th32.1}, it suffices to verify   (\ref{e32.1}) and (\ref{e32.2}). Note that we may assume $F(0)= 0$ since
$F = F- F(0) + F(0)$.

\medskip

Firstly we show (\ref{e32.1}).   Fix $k\leq 0,\ \ell\geq 0$ and $j\in \ZZ$ and write
 \begin{eqnarray*}
F^{(j+\ell)}(\SL)\eta_{j+k}(\SL)&=&\sum_{i\in
\ZZ}(F\eta_{i+j})^{(j+\ell)}(\SL)\eta_{j+k}(\SL)\\
&=&\sum_{i<k-2}+\sum_{i=k-2}^{k+2}+\sum_{i>k+2}\cdots
=I_{jk\ell}+I\!I_{jk\ell}+I\!I\!I_{jk\ell}.
\end{eqnarray*}
The main term is $I\!I_{jk\ell}$; $I_{jk\ell}$ and $I\!I\!I_{jk\ell}$ are
error terms. By (\ref{e32.8}),
$\|\eta_{j+k}(\SL)\|_{p\to p}\leq C$  so
\begin{eqnarray*}
\|I\!I_{jk\ell}\|_{p\to p}&\leq& C
\sum_{i=k-1}^{k+1}\|(F\eta_{i+j})^{(j+\ell)}(\SL)\|_{p\to p}
 \leq  C\sum_{i=-2}^{2}
\alpha(k+i-\ell).
\end{eqnarray*}
Thus,
\begin{eqnarray}\label {e32.9}
\sum_{k\leq 0}\sum_{\ell\geq 0}\sup_j\|I\!I_{jk\ell}\|_{p\to p}
\leq C \sum_{m\leq 2}(|m|+1)\alpha(m)<\infty.
\end{eqnarray}
 To estimate $\|I_{jk\ell}\|_{p\to p}$ set
$$
G (\l)=\sum_{i<k-2}(F\eta_{i+j})^{(j+\ell)}(2^{-j-k}\l)\eta(\l).
$$
Observe  that
 \begin{eqnarray*}
 {d^\gamma\over
d\l^\gamma}\Big(\sum_{i<k-2}(F\eta_{i+j})^{(j+\ell)}(2^{-j-k}\l)\Big)
&=&\int_{|s|\geq
2^{-k-j+2}}\sum_{i<k-2}F(s)\eta_{i+j}(s)2^{j+\ell}2^{(\ell-k)\gamma}(\check{\eta})^{(\gamma)}
(2^{\ell-k}\l-2^{j+\ell}s)ds.
\end{eqnarray*}
  Now $2^{j+\ell}|s|\geq 2^{\ell-k+2}\geq
 2^{\ell-k+1}|\l|$ for  $\l\in [1/2, 2]$. We may estimate the
integral for each $N\in \NN$ by
$$
\int_{|s|\geq 2^{-k-j+2}}
\|F\|_\infty2^{j+\ell}2^{(\ell-k)\gamma}{C_N\over
(2^{j+\ell}|s|)^N}ds.
$$
If $N$ is chosen sufficiently large, this is dominated by
$\|F\|_\infty2^{(\ell-k)(1-N+\gamma)}$. This yields
\begin{eqnarray}\label {e32.10}
\|G\|_{W^{\gamma,\infty}}\leq C\|F\|_\infty
2^{\varepsilon_0(k-\ell)}
\end{eqnarray}
for some $\varepsilon_0>0$ and all $\gamma\in \NN$.  Then by
 (\ref{e32.8})  $\|I_{jk\ell}\|_{p\to p}\leq C\|F\|_\infty
2^{\varepsilon_0(k-\ell)}$.  Hence
$$
\sum_{k\leq 0}\sum_{\ell\geq 0}\sup_j\|I_{jk\ell}\|_{p\to p}<\infty.
$$
We estimat $I\!I\!I_{jk\ell}$ in the similar way as $\|I_{jk\ell}\|_{p\to p}$.
This proves    (\ref{e32.1}).

\medskip

Secondly we show (\ref{e32.2}).   Fix $\ell\geq 0$ and $j\in \ZZ$. Write
\begin{eqnarray*}
F^{(j+\ell)}(\SL)(I-{\Phi}(2^{j}\SL))Q_j(\SL)&=&\sum_{i\in
\ZZ}(F\eta_{i+j})^{(j+\ell)}(\SL)(I-{\Phi}(2^{j}\SL))Q_j(\SL)\\
&=&\sum_{i\leq0} +\sum_{i> 0}\cdots  = A_{j\ell} +B_{j\ell}.
\end{eqnarray*}

\medskip

We start with the  term $A_{j\ell}$.
By condition (\ref{e32.8}) there exists a positive constant $C$ independent of $j$ such that
$\|(I-{\Phi}(2^{j}\SL))Q_j(\SL)\|_{p\to p}\leq C$.  By assumptions of the proposition, we have
$\|(F\eta_{i+j})^{(j+\ell)}(\SL)\|_{p\to p}\leq \alpha(i-l)$.  Hence
$$
 \sup_j \big\|\sum\limits_{\ell\geq0} A_{j\ell}\big\|_{p\to
p}\leq C\sum\limits_{\ell\geq 0}\sum\limits_{i\leq 0}\alpha(i-\ell)
 \leq C\sum\limits_{m\leq 0}(|m|+1)\alpha(m)<\infty.
$$

\medskip

It remains to treat the term $B_{j\ell}$. More precisely, we want to prove that
\begin{eqnarray*}
&&\sup_j \sum_{\ell\geq0}\|\sum_{i>0}(F\eta_{i+j})^{(j+\ell)}(\SL)(I-{\Phi}(2^{j}\SL))Q_j(\SL)\|_{p\to
p}\\
&& = \sup_j
\sum_{\ell\geq0}\|\sum_{i>0}\sum_{k\in
\ZZ}(F\eta_{i+j})^{(j+\ell)}(\SL)\eta_{j+k}(\SL)(I-\Phi(2^{j}\SL))Q_j(\SL)\|_{p\to
p}<\infty.
\end{eqnarray*}
We write
$$
\sum_{i>0}\sum_{k\in \ZZ}
$$
as
\begin{eqnarray*}
&&\sum_{k\leq
0}\sum_{i=1}^\infty+\sum_{k=1}^{\ell}\sum_{i=1}^{k}+\sum_{k=1}^{\ell}\sum_{i=k+1}^{\ell}
+\sum_{k=1}^{\ell}\sum_{i=\ell +1}^{\infty}+\sum_{k=\ell +1}^{\infty}\sum_{i=1}^{\ell}
+\sum_{k=\ell+1}^{\infty}\sum_{i=\ell+1}^{\infty}\cdots\\
&&\quad\quad\quad\quad\quad\quad\quad\quad\quad\quad\quad\quad\quad\quad\quad
\quad\quad\quad=I_{j\ell}+I\!I_{j\ell}+I\!I\!I_{j\ell}+I\!V_{j\ell}+V_{j\ell}+V\!I_{j\ell}.
\end{eqnarray*}

\noindent {\it Claim 1.} For $k\leq 0$,
$$
\|\sum_{i>0}(F\eta_{i+j})^{(j+\ell)}(\SL)\eta_{j+k}(\SL)\|_{p\to
p}\leq C2^{-\varepsilon(\ell-k)}
$$
for some $\varepsilon>0$.

\medskip

\noindent {\it Claim 2.}
$$
\sum_{i=1}^k\|(F\eta_{i+j})^{(j+\ell)}(\SL)\|_{p\to p}\leq
\sum_{i=1}^k\alpha(i-\ell).
$$

\medskip

\noindent {\it Claim 3.} For $k\leq i$,
$\|(F\eta_{i+j})^{(j+\ell)}(\SL)\eta_{j+k}(\SL)\|_{p\to p}\leq
C2^{-\varepsilon(\ell-k)}$, so that
$$
\sum_{i=k+1}^\ell\|(F\eta_{i+j})^{(j+\ell)}(\SL)\eta_{j+k}(\SL)\|_{p\to
p}\leq C(\ell-k)2^{-\varepsilon(\ell-k)}\leq C
2^{-\varepsilon'(\ell-k)}.
$$

\medskip

\noindent {\it Claim 4.} For $\ell\geq k$,
$$
\|\sum_{i=\ell +1}^\infty(F\eta_{i+j})^{(j+\ell)}(\SL)\eta_{j+k}(\SL)\|_{p\to
p}\leq C2^{-\varepsilon'(\ell-k)}.
$$

\medskip

\noindent {\it Claim 5.}
$$
\sum_{i=1}^\ell\|(F\eta_{i+j})^{(j+\ell)}(\SL)\|_{p\to p}\leq
\sum_{i=1}^\ell\alpha(i-\ell).
$$

\medskip

\noindent {\it Claim 6.} For $k\geq \ell$,
$$
\|\sum_{i=\ell +1}^\infty(F\eta_{i+j})^{(j+\ell)}(\SL)\eta_{j+k}(\SL)\|_{p\to
p}\leq C.
$$

\medskip

\noindent If we agree with the  claims 1-6,  we can finish the proof of the
proposition as follows. Let  $\tilde{\eta} \in C_c^{\infty}(0, \infty)$ be a non-negative function
as in (\ref{e32.7}). We note that by condition (\ref{e32.8}),
\begin{eqnarray}\label {e32.11}
\|\tilde{\eta}_{j+k}(\SL)(I-\Phi(2^{j}\SL))Q_j(\SL)\|_{p\to
p}&\leq&
C\|\tilde{\eta}(1-\delta_{2^{-k}}\Phi)\|_{W^{n,\infty}}\nonumber\\
&\leq& C\min\{1,2^{-k}\}.
\end{eqnarray}
By {\it Claim 1},
$$
\sum_{\ell\geq 0}\|I_{j\ell}\|_{p\to p}\leq C\sum_{\ell\geq
0}\sum_{k\leq0}2^{-\varepsilon(\ell-k)}<\infty;
$$
by {\it Claim 2} and (\ref{e32.11}),
$$
\|I\!I_{j\ell}\|_{p\to p}\leq
C\sum_{k=1}^\ell\sum_{i=1}^k\alpha(i-\ell)2^{-\varepsilon''k};
$$
by {\it Claim 3} and (\ref{e32.11}),
$$
\|I\!I\!I_{j\ell}\|_{p\to p}\leq
C\sum_{k=1}^\ell2^{-\varepsilon'(\ell-k)}2^{-\varepsilon''k}\leq
C2^{-\varepsilon\ell};
$$
by {\it Claim 4} and (\ref{e32.11}),
$$
\|I\!V_{j\ell}\|_{p\to p}\leq
C\sum_{k=1}^\ell2^{-\varepsilon'(\ell-k)}2^{-\varepsilon''k}\leq
C2^{-\varepsilon\ell};
$$
by {\it Claim 5} and (\ref{e32.11}),
$$
\|V_{j\ell}\|_{p\to p}\leq
C\sum_{k=\ell +1}^\infty\sum_{i=1}^\ell\alpha(i-\ell)2^{-\varepsilon''k};
$$
by {\it Claim 6} and (\ref{e32.11}),
$$
\|V\!I_{j\ell}\|_{p\to p}\leq C\sum_{k=\ell +1}^\infty2^{-\varepsilon''k}\leq
2^{-\varepsilon\ell}.
$$
Hence
$$
\sum_{\ell>0}\|I\!I\!I_{j\ell}+I\!V_{j\ell}+V\!I_{j\ell}\|_{p\to p}\leq
C\sum_{\ell>0}2^{-\varepsilon\ell}<\infty,
$$
and
\begin{eqnarray*}
\|I\!I_{j\ell}+V_{j\ell}\|_{p\to p}&\leq&
C\sum_{k=1}^\ell\sum_{i=1}^k\alpha(i-\ell)2^{-\varepsilon''k}+
C\sum_{k=\ell +1}^\infty\sum_{i=1}^\ell\alpha(i-\ell)2^{-\varepsilon''k}\\
&=&C\sum_{i=1}^\ell\alpha(i-\ell)\sum_{k=i}^\infty2^{-\varepsilon
k}\leq C\sum_{i=1}^\ell\alpha(i-\ell)2^{-\varepsilon i}.
\end{eqnarray*}
Thus
\begin{eqnarray*}
\sum_{\ell>0}\|I\!I_{j\ell}+V_{j\ell}\|_{p\to p}&\leq&C \sum_{\ell>0}\sum_{i=1}^\ell
2^{-\varepsilon i}\alpha(i-\ell)\\
&\le&C\sum_{m\leq 0}\alpha(m)2^{-\varepsilon m}\sum_{\ell\geq
-m}2^{-\varepsilon \ell}\leq C\sum_{m\leq 0}\alpha(m)<\infty.
\end{eqnarray*}

\medskip

Now {\it Claims 2} and {\it 5} follow immediately from the definition of
$\alpha(j)$. Similarly to the proof  of (\ref{e32.10}), we  use condition  (\ref{e32.8}) to prove {\it Claims 1, 3, 4} and {\it 6}. So to
establish {\it Claims 1, 3, 4} and {\it 6 } we examine
$$
{d^\gamma\over
d\l^\gamma}\Big(\sum_{i}(F\eta_{i+j})^{(j+\ell)}(2^{-j-k}\l)\Big)
$$
for $|\l|\sim 1$, (where the sum is over a range of $i$
depending on which claim we are proving). In any case,  we obtain
$$
\int\sum_{i}F(s)\eta_{i+j}(s)2^{j+\ell}2^{(\ell-k)\gamma}{d^\gamma \check{\eta} \over
d\l^\gamma}  (2^{\ell-k}\l-2^{j+\ell}s)ds.
$$
For {\it Claim 6}, we just use
$|{d^\gamma \check{\eta} \over
d\l^\gamma} |\leq C$ and the fact that the
integrand is supported in a set of measure $\leq C2^{-(j+\ell)}$ to
estimate the integral by $\|F\|_\infty 2^{(\ell-k)\gamma}$; under
the hypotheses of {\it Claims 1, 3} and {\it 4} we have that (essentially)
$2^{\ell-k}|\l|\geq 2\cdot2^{j+\ell}|s|$ if $|\l|\sim 1$, and thus
we may estimate the integral by
$\|F\|_\infty2^{j+\ell}2^{(\ell-k)(\gamma-N)}$ multiplied by the
measure of the support of the integrand, for each $N\in \NN$. These
measures are $2^{-j}$, $2^{-i-j}$ and $2^{-j-\ell}$ respectively,
and so in each case we can dominate the integral by
$C2^{(\ell-k)(\gamma-N+1)}$ for all $N\in\NN$.
 \hfill{}$\Box$

\medskip

Recall that for $0<\alpha<1$,  $\Lambda_{\alpha}$ is the usual Lipschitz space  as defined for example
in \cite{St12}. As a consequence of Proposition~\ref{th32.2}, we have the following result.

\begin{coro}\label {th32.3}  Assume that operator $L$ satisfies  property {\rm (FS)}
and condition  ${\rm (E_{p_0,2})}$ for some $1\leq p_0<2$. Next
assume that for all even bounded Borel function $F$ such that  $\supp F$ is
compact and   $\|F\|_{W^{\beta, \infty}}< \infty$ for some $\beta>n/2,$
 \begin{eqnarray}\label{e32.12}
\sup_{t>0}\|F(t\SL)\|_{p\to p}\leq C\|F\|_{W^{\beta, \infty}}, \ \ \ p_0<p<p'_0.
\end{eqnarray}
Then for any bounded Borel function $F$  such that for all $i, j \in \ZZ$,
 \begin{eqnarray}\label{e32.13}
\|(F\eta_i)^{(j)}(\SL)\|_{r\to
r}+\|\eta \delta_{2^{-i}}F \|_{\Lambda_\alpha}\leq C
\end{eqnarray}

\noindent for some $0<\alpha<1$ and all $p_0<r\leq 2$,
  $F(\SL)$ is bounded on $L^r$ for $p_0<r<p_0'$.
\end{coro}

\noindent
{\bf Proof.}  First, we note
 that $\|\eta \delta_{2^{-i}}F\|_{\Lambda_\alpha}\leq C$ is
equivalent to $\|(F\eta_i)^{(j)} (\SL)\|_{2\to 2}\leq
2^{\alpha(i-j)}$. By interpolation, there exists some
$\alpha'>0$ such that
$$
\|(F\eta_i)^{(j)} (\SL)\|_{r\to r}\leq C2^{\alpha'(i-j)}
$$
for all $p_0<r\leq2$. From Proposition~\ref{th32.2},
$F(\SL)$ is of weak-type $(r,r)$ for all $p_0<r\leq 2$. By duality and interpolation,
$F(\SL)$ is bounded on $L^r$ for $p_0< r<p_0'$.
 \hfill{}$\Box$

\medskip

\noindent
{\bf Proof of Theorem~\ref{th33.1}.} It is enough   to verify
conditions (\ref{e32.13}) of Corollary~\ref{th32.3}.
First we note that by definition
$$
(F\eta_i)^{(j)}(\l)
 = \int F(2^{-i}t)\eta(t)\widehat{\eta_{i-j}}(2^i\l-t)dt
 = (\delta_{2^{-i}}F\eta)^{(j-i)}(2^i\l),
$$
so
 $\|(F\eta_i)^{(j)}(\SL)\|_{r\to r}=\|(\delta_{2^{-i}}F\eta)^{(j-i)}(2^i\SL)\|_{r\to r}.
 $
 Let $\psi$ be a
$C_c^\infty$ even function which is supported on $[-8,8]$ and $\psi(\l)=1$ on $[-4,4]$. Write
$$
\|(\delta_{2^{-i}}F\eta)^{(j-i)}(2^i\SL)\|_{r\to r}\leq
\|(\psi(\delta_{2^{-i}}F\eta)^{(j-i)})(2^i\SL)\|_{r\to r}+
\|((1-\psi)(\delta_{2^{-i}}F\eta)^{(j-i)})(2^i\SL)\|_{r\to r}.
$$
Observe that  the function $\psi(\delta_{2^{-i}}F\eta)^{(j-i)}$ is  supported on $[-8,8]$ and
\begin{eqnarray*}
\|\psi(\delta_{2^{-i}}F\eta)^{(j-i)}\|_{W^{\beta,q}}&\leq& C\|(\delta_{2^{-i}}F\eta)^{(j-i)}\|_{W^{\beta,q}}\\
&=&C \|{\mathscr F}^{-1}\big((1+\xi^2)^{\beta/2}\widehat{\delta_{2^{-i}}F\eta}(\xi)\eta_{i-j}(\xi)\big)\|_{q}\\
&\leq&C \|{\mathscr F}^{-1}\big((1+\xi^2)^{\beta/2}\widehat{\delta_{2^{-i}}F\eta}(\xi)\big)\|_{q}\\
&=&C \|F(2^{-i}\l)\eta(\l)\|_{W^{\beta,q}}.
\end{eqnarray*}
From condition (\ref{e333.1}), we have that $\|\eta(\l)F(2^{-i}\l)\|_{W^{\beta,q}}<\infty.$
Therefore, we use our assumption (\ref{e33.1}) to obtain
\begin{eqnarray}\label{elll}
\|(\psi(\delta_{2^{-i}}F\eta)^{(j-i)})(2^i\SL)\|_{r\to r}\leq C
\end{eqnarray}
with $C>0$ independent of $i$ and $j$.

We  estimate the term  $\|((1-\psi)(\delta_{2^{-i}}F\eta)^{(j-i)})(2^i\SL)\|_{r\to r}$.
For $k\in {\mathbb N}$ and all $\l\neq 0$, we have by  elementary
calculation,
\begin{eqnarray*}
&&\hspace{-0.8cm}\Big|(1-\psi(\l))\frac{d^k}{d\l^k}\big((\delta_{2^{-i}}F\eta)^{(j-i)})(\l)\big)\Big|\\
&&\leq C2^{(j-i)(k+1)}\big|1-\psi(\l)\big|\int
 \big|F(2^{-i}u)\eta(u)\big|\big(1+2^{j-i}|\l-u|\big)^{-k-1}du
 \leq   C|\l|^{-k},
\end{eqnarray*}
where we use the fact that    $|\l|\geq 4$ and   $|u|\leq 2$.
 We then  apply  (a) of
Lemma~\ref{le31.3} to obtain
 $$
\|((1-\psi)(\delta_{2^{-i}}F\eta)^{(j-i)})(2^i\SL)\|_{r\to r}\leq C.
 $$
This estimate  in combination with (\ref{elll}) shows that
\begin{eqnarray}\label{e33.55}
\|(F\eta_i)^{(j)}(\SL)\|_{r\to r}=\|(\delta_{2^{-i}}F\eta)^{(j-i)}(2^i\SL)\|_{r\to r}\leq C
\end{eqnarray}
 for some constant $C>0$ independent of $i,j.$

Now  we recall that if $1\leq q\leq \infty$ and $\beta-1/q>0$, then
$$
W^{\beta, q} \subseteq  B^{\beta}_{q, \, \infty} \subseteq  B^{\beta-{1\over q}}_{\infty, \, \infty}
\subseteq
\Lambda_{\min\{\beta-{1\over q}, \, {1\over 2}\}}
$$
 and  $\|F\|_{\Lambda_{\min\{\beta-{1\over q}, \, {1\over 2}\}}}\leq C\|F\|_{W^{\beta, q}}$.
See, e.g., \cite[Chap.~VI ]{BL} for more details.  We obtain
\begin{eqnarray}\label{e33.5}
\|\eta \delta_{2^{-i}}F \|_{\Lambda_{\min\{\beta-1/q, \, 1/2\}}}\leq
C\|\eta \delta_{2^{-i}}F\|_{W^{\beta,q}}\leq C'.
\end{eqnarray}
This estimate and (\ref{e33.55}) prove condition   (\ref{e32.13}) of Corollary~\ref{th32.3}.
The   proof of Theorem~\ref{th33.1} is finished.
  \hfill{}$\Box$


\bigskip

\section{Endpoint estimates for Bochner-Riesz means}\label{BRS}
\setcounter{equation}{0}

We have seen in Corollary \ref{coro3.3} that Bochner-Riesz means are bounded on $L^p$ provided the
order $\delta $ satisfies $\delta > \max\big\{n(1/p -1/s) -1/q, 0\big\}$.
In this section we    prove that our restriction type condition
 implies   endpoint estimates   for
 Bochner-Riesz  means. Our approach
is inspired by the results of Christ  and Tao  \cite{C1, C2, T2}.
As in the rest of the paper  we assume  that   $(X, d, \mu)$ is a
metric measure  space satisfying condition
(\ref{eq2.2}) with a  homogeneous dimension $n$.

For any given  $p\in [1,2)$ and $q\in [1, \infty]$  we define
$$
\delta_q(p)=\max\Big\{0,n\Big|\frac{1}{p}-\frac{1}{2}\Big|-\frac{1}{q}\Big\}.
$$
For simplicity we will write $\delta(p)$ instead of $\delta_2(p)$. As in Sections~\ref{sec3} and
\ref{sec4}
we discuss two type of results corresponding to estimates ${\rm (ST^{q}_{p, 2})}$ or
 ${\rm (SC^{q,1}_{p, 2})}$.

 \begin{theorem}\label{th4.1} \,
 Assume that operator $L$ satisfies   property {\rm (FS)}
and condition  ${\rm (ST^{q}_{p, 2})}$ for some $p, q$ satisfying
 $1\leq p <2$ and $1\leq q \leq \infty$.
Then the operator $S_R^{\delta_q(p)} (L)$ is of weak-type $(p,p)$
uniformly in $R$.
\end{theorem}

\medskip

The next theorem
 is a variation of Theorem~\ref{th4.1}. As in Theorems~\ref{th3.2} and ~\ref{th31.2},
 this variation can be used in the  case of operators with nonempty
 pointwise spectrum.

 \medskip

 \begin{theorem}\label{th4.2}  Assume that $\mu(X)<\infty$.
 Assume that operator $L$ satisfies  property {\rm (FS)}
and condition  ${\rm (SC^{q,1}_{p, 2})}$ for some $p, q$ satisfying
 $1\leq p <2$ and $1\leq q \leq \infty$.
Then the operator $S_R^{\delta_q(p)} (L)$ is of weak-type $(p,p)$
uniformly in $R$.
 \end{theorem}

The proofs of Theorem~\ref{th4.1}~and~\ref{th4.2} require three technical lemmas which we discuss first.
We commence with the following observation.

\begin{lemma}\label{le4.3}
 For each $k\leq 0$ there exists a decomposition
\begin{eqnarray}\label{eq4.6}
S_R^{\delta_q(p)}(\l^2)=\eta_k(\l)n_k(\l)+S_R^{\delta_q(p)}(\l^2)n_k(\l)
\end{eqnarray}
such that
\begin{itemize}
\item[(a)]
 Functions  $n_k$ are  even and their Fourier transforms are supported in
$[-2^k/R,\, 2^k/R]$, i.e.,  $\supp \widehat{n_k}\subset
[-2^k/R,\, 2^k/R]$;

\item[(b)] Functions $\eta_k$ are continuous, even and
$\sum_{k=-\infty}^0|\eta_k(\l)|^2\leq C$ with $C$ independent of
$\l$ and $R$;

\item[(c)]  For certain arbitrarily large  $N \in \NN$ there exists a constant $C$ such that
$$|n_k(\l)|\leq
C\Big(1+{2^k|\l|\over R}\Big)^{-N}.$$
\end{itemize}
\end{lemma}

\noindent {\bf Proof.} Following  \cite{T2} we  consider the function  $\Phi(\l)=6\lambda^{-3}
 (\l-\sin \l) $. Note that
$\widehat{\Phi}(t)=3\pi(1-|t|)_+^2$. Set
$n_k(\l)= \Phi^{{N/ 2}} (2^k\l(RN)^{-1})$. We use term $(RN)^{-1}$ rather than $R^{-1}$
to ensure that $\supp \widehat{n_k} \subset [-2^k/R,\, 2^k/R]$.
Next  we write
\begin{eqnarray*}
S_R^{\delta_q(p)}(\l^2)&=&\left({S_R^{\delta_q(p)}(\l^2)\big(1-n_k(\l)\big)\over
n_k(\l)}\right) n_k(\l)+S_R^{\delta_q(p)}(\l^2)n_k(\l)\\[2pt]
&=&\eta_k(\l)n_k(\l)+S_R^{\delta_q(p)}(\l^2)n_k(\l).
\end{eqnarray*}
Verifying  conditions (a) and (c) is straightforward.
By definition of
$n_k(\l)$,
  there exists some constant
$C$   depending only on $N$ such that $n_k(\l)\geq C $  and
$|S_R^{\delta_q(p)}(\l^2)(1-n_k(\l))|\leq C 2^k$ for all   $|\lambda|\leq R$  and $k\le 0$.
 This proves condition  (b) because $\supp S_R^{\delta_q(p)}(\l^2)\subset [-R,R]$. \hfill{} $\Box$

\bigskip

 The following lemma comes from  \cite{T2}.

\begin{lemma} \label {le4.4}
For each $k>0$, there exists a decomposition
$S_R^{\delta_q(p)}(\l^2)=m_k(\l)+\eta_k(\l) n_k(\l)$ such that:
\begin{itemize}
\item[(a)]  Functions $\widehat{m_k}$ and $\widehat{n_k}$ are even and supported
on $[-2^k/R,2^k/R]$;

\item[(b)] $\eta_k$ are continuous and for all $\l>0$,
$\sum_{k=1}^\infty|\eta_k(\l)|^2\leq C$ with $C$ independent of $\l$
and $R$;

\item[(c)] For certain arbitrarily large   $N\in \NN$ there exists a constant $C$ such that
$$
|n_k(\l)|\leq C2^{-\delta_q(p) k}\Big(1+2^k\Big|1-{|\l|\over R}\Big|\Big)^{-N}.
$$
\end{itemize}
\end{lemma}

\noindent {\bf Proof.} For the proof, we refer  to Lemma 2.1 in  \cite{T2}. \hfill{} $\Box$

\bigskip Next we discuss the proof of the last  lemma required in the proofs of Theorems~\ref{th4.1}~and~\ref{th4.2}.

 \begin{lemma}\label{le4.5} Suppose that $L$  is a self-adjoint operator on $L^2(X)$. Assume that
 $\{Q_k\}_{k\in{\mathbb N}}$ is a family of continuous real-valued functions such that
  $\sum_k |Q_k(\lambda)|^2\leq C$ for some constant
 $C$ independent of $\lambda$. Then for any sequence of functions $\{f_k\}_{k\in{\mathbb N}}$ on $X,$
\begin{eqnarray}\label{eq4.22}
\big\|\sum_{k}Q_k(\SL)f_k\big\|_2^2\leq C\sum_{k}\big\|f_k\big\|_2^2.
\end{eqnarray}
 \end{lemma}

\noindent
{\bf Proof.}
  Note that
\begin{eqnarray*}
 \big\|\sum_{k}Q_k(\SL)f_k \big\|_2^2
&=& \Big\langle\sum_{k}Q_k(\SL)f_k,
\sum_{k}Q_k(\SL)f_k \Big\rangle\nonumber\\
&=&\int_{-\infty}^{+\infty} \sum_{k}\sum_{j}Q_k(\lambda)Q_j(\lambda)
 d\big\langle E_{\sqrt{L}}(\lambda) f_k,
f_j\big\rangle.
\end{eqnarray*}

Under the assumption of Lemma~\ref{le4.5} the above integral is a limit of its Riemannian approximations,
see  page 310 of \cite{Yo}. Therefore let us consider Riemannian partition of $[\alpha, \beta]$
$$
\alpha=\lambda_1<\lambda_2<...<\lambda_n=\beta,\ \ \ \ \lambda_l^{'}\in (\lambda_{\ell},\,  \lambda_{\ell+1}]
$$
for some $-\infty<\alpha<\beta<\infty$.
Now to finish the proof of  Lemma~\ref{le4.5}
it is enough to observe that
\begin{eqnarray*}
\left|\sum_{\ell}  \sum_{k}\sum_{j}Q_k(\lambda_{\ell}^{'})
Q_j(\lambda_{\ell}^{'})
 \big\langle E_{\sqrt{L}}(\lambda_{\ell}, \,
\lambda_{\ell+1}] f_k, \, f_j \big\rangle\right| \hspace{5.5cm}\\
 \leq  \sum_{\ell}  \sum_{k}\sum_{j}|Q_k(\lambda_{\ell}^{'})| |Q_j(\lambda_{\ell}^{'})|
\sqrt{\big\langle E_{\sqrt{L}}(\lambda_{\ell}, \, \lambda_{\ell+1}]f_k, \, f_k \big\rangle}
 \sqrt{\big\langle E_{\sqrt{L}}(\lambda_{\ell}, \, \lambda_{\ell+1}] f_j, \,  f_j \big\rangle}  \\
  =   \sum_{\ell}  \bigg(\sum_{k}  |Q_k(\lambda_{\ell}^{'})|
\sqrt{\big\langle E_{\sqrt{L}}(\lambda_{\ell}, \, \lambda_{\ell+1}] f_k,\,  f_k \big\rangle}  \, \bigg)^2  \\
 \leq  \sum_{\ell}  \Big( \sum_{k}  |Q_k(\lambda_{\ell}^{'})|^2\Big)
  \sum_{k}  \big\langle E_{\sqrt{L}}(\lambda_{\ell}, \, \lambda_{\ell+1}] f_k, \, f_k \big\rangle  \\
   \leq  C\sum_{\ell}  \sum_{k}  \big\langle E_{\sqrt{L}}(\lambda_{\ell}, \,
    \lambda_{\ell+1}] f_k, \, f_k \big\rangle\\
   \leq  C\sum_{k}
  \|f_k\|_2^2,\nonumber
\end{eqnarray*}
where we used the fact that
$\sum_{k}|Q_k(\l)|^2\leq C$. \hfill{} $\Box$

\bigskip

\noindent
{\bf Proof  of Theorem~\ref{th4.1}.}
 Assume that  condition  ${\rm
(ST^{q}_{p, 2})}$ holds for  some $1\leq
p<2$ and  $1\leq q\leq \infty$.
  Fix a $f\in L^p$ and $\alpha> \mu(X)^{-1/p}\|f\|_p$,
and apply  the Calder\'on-Zygmund decomposition at height $\alpha$
to $|f|^p.$ There  exist  constants $C$ and  $K$ so that
\begin{enumerate}
\item[(i)] $
f=g+b=g+\sum_{j}b_j;
$

\vskip 1pt
\item[(ii)]   $\|g\|_p\leq C\|f\|_p,$
 $\|g\|_{\infty}\leq C\alpha$;

\vskip 1pt
\item[(iii)] $b_j$  is supported in  $B_j$ and  $\# \{j: x\in 4B_j\} \leq K$ for all $x\in X$;

\vskip 1pt
\item[(iv)]
$\int_{X}\ |b_j|^pd\mu\leq C\alpha^p \mu(B_j), $  and
$\sum_{j}\mu(B_j) \leq C{\alpha}^{-p} \|f\|_p^p.$
\end{enumerate}
Note that by (ii) $\alpha^{p-2}\|g\|_2^2\leq C\|f\|_p^p.$

\smallskip

 Let $r_{B_j}$ be the radius of $B_j$ and denote by $J_k=\big\{j: \, 2^k/R \le r_{B_j}<2^{k+1}/R\big\}$.  Write
\begin{eqnarray*}
f =g+\sum_{j}b_j=
 g+\sum_{k\leq 0}\sum_{j\in J_k}b_j+\sum_{k>0}\sum_{j\in J_k}b_j
 = g+h_1+h_2.
\end{eqnarray*}
By a standard argument it is enough to show that there exists a constant $C>0$ independent
of $R$ and $\alpha$ such that for every $\alpha> \mu(X)^{-1/p}\|f\|_p$,
\begin{equation}\label{eq4.1}
\mu\big(\{x:S_R^{\delta_q(p)}(L)(g)(x)>\alpha\}\big)
  \leq
C\alpha^{-p}\|f\|_p^p
\end{equation}
and  that for $i=1,2,$
\begin{equation}\label{eq4.2}
\mu\big(\big\{x:S_R^{\delta_q(p)}(L)\big(h_i\big)(x)>\alpha\big\}\big)
  \leq
C\alpha^{-p}\|f\|_p^p.
\end{equation}

\medskip

 Note that $\sup_{\l,R>0}\Big(1-{\l\over
R}\Big)^{\delta_q(p)}_+ =1$  and that by by (ii) $\alpha^{p-2}\|g\|_2^2\leq C\|f\|_p^p.$
Hence by   spectral
 theory
\begin{eqnarray}\label{eq4.4}
\mu\big(\{x:S_R^{\delta_q(p)}(L)(g)(x)>\alpha\}\big)&\leq&
\alpha^{-2}\|S_R^{\delta_q(p)}(L)(g)\|_2^2\leq \alpha^{-2}\|g\|_2^2\nonumber\\
& \leq&
C\alpha^{-p}\|f\|_p^p.
\end{eqnarray}

Next we  prove (\ref{eq4.2})   for     $i=1$.
 By  the spectral theorem  and equality (\ref{eq4.6})
\begin{eqnarray*}
 \sum_{k\leq 0}\sum_{j\in
J_k}S_R^{\delta_q(p)}(L)b_j &=&
\sum_{k\leq
0}\eta_k({\SL})\Big(\sum_{j\in
J_k}n_k({\SL})b_j \Big) \\
 &+&  S_R^{\delta_q(p)}(L)\Big(\sum_{k\leq 0}\sum_{j\in
J_k}n_k({\SL})b_j\Big).  \nonumber
\end{eqnarray*}
Applying   the spectral theorem and Lemma~\ref{le4.5}  with
$Q_k(\lambda)=\eta_k(\lambda)$ yields
\begin{eqnarray}\label{eq4.7}
 \Big\|\sum_{k\leq 0}\sum_{j\in
J_k}S_R^{\delta_q(p)}(L)b_j\Big\|^2_2 &\leq& C\sum_{k\leq
0}\Big\|\sum_{j\in J_k}n_k({\SL })b_j\Big\|^2_2+C\Big\|\sum_{k\leq
0}\sum_{j\in J_k}n_k({\SL})b_j\Big\|^2_2.
\end{eqnarray}
Next $\supp \widehat{n_k} \subseteq [-2^k/R,2^k/R]$ so by Lemma~\ref{le2.1},
$$
\supp K_{n_k({\sqrt{L} })} \subseteq \big\{(x,y)\in X\times
X: d(x,y)\leq 2^k/R\big\}.
$$
Hence if $j\in J_k$, then
$\supp n_k(\sqrt{L})b_j\subseteq  4B_j$. Thus by (iii) there exists
constant $C>0$ such that
\begin{eqnarray} \hspace{1cm}
\sum_{k\leq
0}\Big\|\sum_{j\in J_k}n_k({\SL })b_j\Big\|^2_2+\Big\|\sum_{k\leq
0}\sum_{j\in J_k}n_k({\SL})b_j\Big\|^2_2  \le  C\sum_{k\leq 0}\sum_{j\in
J_k}\Big\|n_k({\SL})b_j\Big\|^2_2.\nonumber
\end{eqnarray}
Next, by Proposition~\ref{prop2.3} and Remark~\ref{re2.2},
\begin{eqnarray*}
\Big\|n_k({\SL})b_j\Big\|_2
&=&\Big\|n_k({\SL})P_{B_j}b_j\Big\|_2\nonumber\\
&\leq&\Big\|n_k({\SL})P_{B_j}\Big\|_{p\to2}
\big\|b_j\big\|_p\nonumber\\
&\leq&\Big\| n_k({\SL}) \Big(I+2^k{\SL\over R}\Big)^{N}
\Big\|_{2\to2}
\Big\|\Big(I+2^k{\SL\over R}\Big)^{-N}P_{B_j}\Big\|_{p\to2}\|b_j\|_p\nonumber\\
&\leq&C\Big\|\Big(I+2^k{\SL\over R}\Big)^{-N}P_{B_j}
\Big\|_{p\to2}\|b_j\|_p\nonumber\\
&\leq&C\mu(B_j)^{{1\over 2}-{1\over p}}\alpha\mu(B_j)^{{1\over p}}\nonumber\\
&\leq&C\alpha\mu(B_j)^{1/2}.
\end{eqnarray*}
Hence by (iv)
\begin{eqnarray}\label {eq4.9}\hspace{1cm}
\mu\Big(\Big\{x:\Big|S_R^{\delta_q(p)}(L)\Big(\sum_{k\leq
0}\sum_{j\in J_k}b_j\Big)\Big|>\alpha\Big\}\Big)
&\leq&C\alpha^{-2}\Big\|S_R^{\delta_q(p)}(L)\Big(\sum_{k\leq
0}\sum_{j\in
J_k}b_j\Big)\Big\|^2_2\nonumber\\
&\leq&C \alpha^{-p}\|f\|_p^p.
\end{eqnarray}

\medskip

Now, we prove (\ref{eq4.2})    for    $i=2$.
Let $\Omega^*=\bigcup_{j\in \mathbb{N}}4B_j$. \noindent
By (\ref{eq2.1}) that
$$\mu(\Omega^*) \leq
C\sum_{j}\mu(B_j)\leq C \alpha^{-p}\|f\|_p^p.$$
 Hence it is enough to show that
\begin{eqnarray}\label{eq4.10}
\big\|S_R^{\delta_q(p)}(L)\big(\sum_{k> 0}\sum_{j\in
J_k}b_j\big)\big\|_{L^2(X\backslash\Omega^*)}^2 \leq C\alpha^2\sum_j
\mu(B_j)\leq C \alpha^{2-p}\|f\|_p^p.
\end{eqnarray}

\medskip

Using the decomposition from  Lemma~\ref{le4.4} we write
\begin{eqnarray}\label{eq4.11} \hspace{1.5cm}
S_R^{\delta_q(p)}(L)\Big(\sum_{k> 0}\sum_{j\in J_k}b_j\Big)=
\sum_{k> 0}\sum_{j\in J_k}m_k(\SL)b_j +\sum_{k>
0}\eta_k(\SL)n_k(\SL)\Big(\sum_{j\in J_k}b_j\Big).
\end{eqnarray}
Recall that $\widehat{m_k}$ is even
and supported in $[-2^k/R,2^k/R]$. By Lemma~\ref{le2.1}
$$
\supp K_{m_k(\sqrt{L})} \subset \big\{(x,y)\in X\times
X: d(x,y)\leq {2^k\over R}\big\}.
$$
This implies that if  $x\in X\backslash\Omega^*$, then $
m_k(\sqrt{L})b_j(x)=0$ for any $j\in J_k $ and $k>0$ so it makes
no contribution to (\ref{eq4.10}). By  \eqref{eq4.11} and Lemma~\ref{le4.5}
\begin{eqnarray}\label{eq4.12}
 \big\|S_R^{\delta_q(p)}(L)\big(\sum_{k>
0}\sum_{j\in J_k}b_j\big)\big\|_{L^2(X\backslash\Omega^*)}^2
&\leq &   \|\sum_{k>
0}\eta_k(\SL)n_k(\SL)\big( \sum_{j\in J_k}b_j\big)\big\|^2_2\nonumber\\
&\leq& C\sum_{k> 0}\big\|n_k(\SL)\big( \sum_{j\in J_k}b_j\big)\big\|_2^2.
\end{eqnarray}
Next
 $\widehat{n_k}$ is  even and supported on $[-2^k/R,2^k/R]$ so
by  Lemma~\ref{le2.1}
  $$\supp K_{n_k(\sqrt{L})} \subseteq \big\{(x,y)\in X\times
X: d(x,y)\leq 2^k/R\big\}.
$$
Hence
 $\supp n_k(\sqrt{L})b_j\subseteq 4B_j$ for $j\in J_k$. By (iii)
 there exists a constant $C>0$ such that
\begin{eqnarray}\label{eq4.13}
\sum_{k> 0}\big\|n_k(\SL)\big( \sum_{j\in J_k}b_j\big)\big\|_2^2
&\leq& C\sum_{k> 0} \sum_{j\in
J_k}\big\|n_k(\SL)b_j\big\|_2^2.
\end{eqnarray}
To continue,  fix  $\psi\in C_c^\infty$, even and
supported in $[-2,2]$ such that  $\psi(\l)=1$ for $|\l|\leq
3/2$. Write
\begin{eqnarray}\label{eq4.14}
n_k(\SL)=n_k(\SL)\psi\Big({\SL\over R}\Big)+n_k(\SL)\Big(1-\psi\Big({\SL\over R}\Big)\Big).
\end{eqnarray}
Then  $\supp (1-\psi({\l/
R}))\subset(-\infty,-3R/2]\cup[3R/2,\infty)$ so  for every  $k>0$
and $\delta_q(p)\geq 0$,
\begin{eqnarray}
\Big|n_k(\l)\Big(1-\psi\Big({\l\over R}\Big)\Big)\Big|&\leq& C \Big|1-\psi\Big({\l\over
R}\Big)\Big|2^{-\delta_q(p)
k}\Big(1+2^k\Big|1-{|\l|\over R}\Big|\Big)^{-N}\nonumber\\
&\leq&C \Big(1+\Big|{2^k\l\over R}\Big|\Big)^{-N}\nonumber,
\end{eqnarray}
and we use a similar argument as in the proof of
 (\ref{eq4.2})    for    $i=1$ to conclude that
\begin{eqnarray}\label{eq4.15}
\Big\|n_k(\SL)\Big(I-\psi\Big({\SL\over R}\Big)\Big)b_j\Big\|_2\leq
C\alpha\mu(B_j)^{1/2}.
\end{eqnarray}
On the other hand,
$\supp n_k(\l)\psi(\l/R)\subset [-2R, 2R]$ so by ${\rm (ST^{q}_{p, 2})}$
\begin{eqnarray}
\Big\|n_k(\SL)\psi\Big({\SL\over R}\Big)b_j\Big\|_2
 &=& \Big\|n_k(\SL)\psi\Big({\SL\over R}\Big)P_{B_j}b_j\Big\|_2\nonumber\\
 &\leq& \Big\|n_k(\SL)\psi\Big({\SL\over R}\Big)P_{B_j}\Big\|_{p\to2}\|b_j\|_p\nonumber\\
 &\leq& C \alpha\mu(B_j)^{1\over 2}  2^{kn({1\over
p}-{1\over 2})}
\big\|\delta_{2R}\big(n_k(\l)\psi\big({\l/R}\big)\big)\big\|_{q}
\nonumber.
\end{eqnarray}
Now
\begin{eqnarray*}
\Big\|\delta_{2R}\Big(n_k(\l)\psi\Big({\l\over
R}\Big)\Big)\Big\|_{q}
&\leq&C \Big(\int_0^1|n_k(2R\l)|^qd\l\Big)^{{1/q}}\nonumber\\
&\leq&C2^{-\delta_q(p) k}\Big(\int_0^1\big(1+2^k|1-2\l|\big)^{-Nq}d\l\Big)^{{1/q}}\nonumber\\
&\leq&  C2^{-\delta_q(p) k}  2^{-{k\over q}} = C2^{-nk({1\over
p}-{1\over 2})}.
\end{eqnarray*}
This yields
\begin{eqnarray}\label{eq4.16}
\Big\|n_k(\SL)\psi\Big({\SL\over R}\Big)b_j\Big\|_2 &\leq&C\alpha
\mu(B_j)^{1/2}.
\end{eqnarray}
 By (\ref{eq4.15}) and (\ref{eq4.16})
  \begin{eqnarray*}
\big\|n_k(\SL)b_j\big\|_2\leq C\alpha\mu(B_j)^{1/2}.
\end{eqnarray*}
The rest of  the proof  of (\ref{eq4.2})    for    $i=2$ is   similar to the  case $i=1$.
\hfill{}$\Box$

\medskip

\noindent {\bf Proof  of   Theorem~\ref{th4.2}.}\,   Assume that   condition ${\rm
(SC^{q,1}_{p, 2})}$ holds  for some $1\leq
p<2$ and $1\leq q\leq \infty$.
 The proof of   Theorem~\ref{th4.2} is almost identical to that of Theorem~\ref{th4.1} except some minor
 technical complications so we only give a brief sketch of it.

We apply the Calder\'on-Zygmund decomposition at height
$\alpha$ to $|f|^p$ to get the same decomposition
\begin{eqnarray*}
f =
 g+\sum_{k\leq 0}\sum_{j\in J_k}b_j+\sum_{k>0}\sum_{j\in J_k}b_j
 =g+h_1+h_2
\end{eqnarray*}
as in Theorem~\ref{th4.1}.  The proof of week type estimates for $g$ and $h_1$ uses the simple
observation that ${\rm
(SC^{q,1}_{p, 2})}\Rightarrow {\rm (SC^{\infty}_{p, 2})}
\Leftrightarrow {\rm (ST^{\infty}_{p, 2})}\Rightarrow {\rm (E_{p,
2})}$ (see Remark~\ref{prop3.9} and Proposition~\ref{prop2.3}) and is essentially the same
as the corresponding argument in the proof of Theorem~\ref{th4.1}.

It remains to show that
\begin{equation}\label{eq4.18}
\mu\big(\big\{x:S_R^{\delta_q(p)}(L) (h_2) (x)>\alpha\big\}\big)
  \leq
C\alpha^{-p}\|f\|_p^p.
\end{equation}
To show (\ref{eq4.18})
 we note that if $\mu(X)$ is finite, then we may assume $X=B(x_0, 1)$ for some $x_0\in X$.
Thus  the radius of each $B_j$ in the Calder\'on-Zygmund  decomposition  satisfies  $2^k/R\leq 4$.
If $R\leq 4$,  then  $k\leq 4$. Hence one can use the same argument as in the proof of Theorem~\ref{th4.1}
 for $i=1$.

   \bigskip

   Next we consider the remaining case $R>4$. Using the decomposition described in  Lemma~\ref{le4.4}
   it is not difficult to note that
 to finish the proof it is enough to show that
\begin{eqnarray}\label{eq4.19}
\Big\|n_k(\SL)\Big(I-\psi\Big({\SL\over R}\Big)\Big)b_j\Big\|_2&\leq& C\alpha\mu(B_j)^{1/2},\\[3pt]
\Big\|n_k(\SL)\psi\Big({\SL\over R}\Big)b_j\Big\|_2 &\leq& C
\alpha\mu(B_j)^{1/2},\label{eq4.20}
\end{eqnarray}
where $n_k$ is defined in Lemma~\ref{le4.4} and $\psi$ is a function in (\ref{eq4.14}). The proof of
(\ref{eq4.19}) is similar to that of (\ref{eq4.15}).

To prove (\ref{eq4.20}) set
 $N=[2R]+1$.
By  condition ${\rm (SC^{q,1}_{p, 2})}$
\begin{eqnarray}
\Big\|n_k(\SL)\psi\Big({\SL\over R}\Big)b_j\Big\|_2 &\leq&C
\mu(B_j)^{({1\over 2}-{1\over p})} 2^{kn({1\over p}-{1\over 2})}
\big\|\delta_{N}\big(n_k(\l)\psi
\big({\lambda/R}\big)\big)\big\|_{N,q}\|b_j\|_p\nonumber.
\end{eqnarray}
Next (assuming that $\sup |\psi| =1$)
 \begin{multline*}
 2^{1/q}\big\|\delta_{N}\big(n_k(\l)\psi \big({\lambda/R}\big)\big)\big\|_{N,q}
  \leq   \Big(\frac{1}{N}\sum_{\ell=1-N}^{N}\sup_{\l\in
[\ell-1,\ell)}|n_k(\l)|^q\Big)^{{1/q}}\\
 \leq  \Big(\frac{1}{N}\sum_{\ell=1-N}^{-[R]-2}\sup_{\l\in
[\ell-1,\ell)}|n_k(\l)|^q\Big)^{{1/q}}
 +
\Big(\frac{1}{N}\sum_{\ell=-[R]-1}^{-[R]+3}\sup_{\l\in
[\ell-1,\ell)}|n_k(\l)|^q\Big)^{{1/q}} \\
 + \Big(\frac{1}{N}\sum_{\ell=-[R]+4}^{0}\sup_{\l\in
[\ell-1,\ell)}|n_k(\l)|^q\Big)^{{1/q}}
 +
 \Big(\frac{1}{N}\sum_{\ell=1}^{[R]-3}\sup_{\l\in
[\ell-1,\ell)}|n_k(\l)|^q\Big)^{{1/q}}\\
 + \Big(\frac{1}{N}\sum_{\ell=[R]-2}^{[R]+2}\sup_{\l\in
[\ell-1,\ell)}|n_k(\l)|^q\Big)^{{1/q}}
 + \Big(\frac{1}{N}\sum_{\ell=[R]+3}^{N}\sup_{\l\in
[\ell-1,\ell)}|n_k(\l)|^q\Big)^{{1/q}} \\
  =  I+I\!I+I\!I\!I+I\!V+V+V\!I.
\end{multline*}
Let $M$ be a  sufficiently large natural number. By
Lemma~\ref{le4.4}
\begin{eqnarray*}
V&=&\Big(\frac{1}{N}\sum_{\ell=[R]-2}^{[R]+2}\sup_{\l\in
[\ell-1,\ell)}|n_k(\l)|^q\Big)^{{1/q}}\\
&\leq&C 2^{-\delta_q(p)
k}\Big(\frac{1}{N}\sum_{\ell=[R]-2}^{[R]+2}\sup_{\l\in
[\ell-1,\ell)}\Big(1+2^k\Big|{\l\over R}-1\Big|\Big)^{-Mq}\Big)^{{1/q}}\\
&\leq&C 2^{-\delta_q(p) k}R^{-{1/q}} \leq C 2^{-\delta_q(p) k} 2^{-{k\over q}} \big({2^k\over R}\big)^{{1/q}}
\\
&\leq&C 2^{-nk({1\over p}-{1\over 2})}
\end{eqnarray*}
and
\begin{eqnarray*}
V\!I&\leq&C \Big(\frac{1}{N}\sum_{\ell=[R]+3}^{N}\sup_{\l\in
[\ell-1,\ell)}\Big|2^{-\delta_q(p) k}\Big(1+2^k\Big({\l\over R}-1\Big)\Big)^{-M}\Big|^q\Big)^{{1/q}}\\
&\leq&C 2^{-\delta_q(p)
k}\Big(\int_{R}^\infty\Big(1+2^k\Big({\l\over R}-1\Big)\Big)^{-Mq}R^{-1}d{\l }\Big)^{{1/q}}\\
&\leq&C 2^{-nk({1\over p}-{1\over 2})}.
\end{eqnarray*}
A similar argument as in  $V$   shows that $I\!I\leq C
2^{-nk(1/p-1/2)}$; the similar argument as in  $V\!I$   shows that each  of
$I,$ $I\!I\!I$ and $I\!V$ is less than $C 2^{-nk(1/p-1/2)}$. Thus
\begin{eqnarray*}
\big\|\delta_{N}\big(n_k(\l)\psi \big({\lambda/R}\big)\big)\big\|_{N,q}\leq
C2^{-nk({1\over p}-{1\over 2})}.
\end{eqnarray*}
This  finishes the proof of Theorem~\ref{th4.2}.
  \hfill{}$\Box$

 \bigskip

\part{Dispersive and  restriction estimates}

\section{Dispersive and Strichartz estimates}
\setcounter{equation}{0}

Let  $(X, d,\mu)$ be   a   metric measure  space.
Next let $L$ be a non-negative
self-adjoint
operator acting on $L^2(X)$. In virtue of the spectral theory,
we can define the semigroup $\exp(-zL)$ for all $z \in \CC$ with $\Re z  \ge 0$ and such that
 $$\|\exp(-zL)\|_{2 \to 2} \le 1. $$
 We   say that  the operator $L$ satisfies a  dispersive type
estimates if there exist constants $n$ and $C$ such that
\begin{equation}  \label{eq5.1}
\|\exp\left(isL\right)\|_{1\to\infty}\le C|s|^{-n/2}, \quad \forall \, s\in
\mathbb{R}\setminus \{0\}.
\end{equation}
Of course, the standard Laplacian on $\RR^n$ satisfies the dispersive estimate.
Such estimates are of importance in analysis and PDE. In particular, they imply
endpoint Strichartz estimates (see  Keel and Tao \cite{KT}).  We   refer to
Strichartz endpoint estimates for the corresponding Schr\"odinger equation as
\begin{equation}\label{str-end}
\int_{\RR}  \| e^{it L} f \|_{{\frac{2n}{n-2}}}^2 dt \le C \| f \|_{2}^2,  \, \, f \in L^2.
\end{equation}
This endpoint estimate together with the obvious fact
$$\| \exp(itL) f \|_{{L_t^\infty}L_x^2} \le \| f \|_{L^2}$$
give $L_t^pL_x^q$ Strichartz estimates. See \cite{KT} for more details.
Our aim  will be   to  explain how  sharp spectral multipliers follow from dispersive or Strichartz estimates.

It is natural to consider the dispersive estimate  (\ref{eq5.1}) in
conjunction with the smoothing condition
\begin{equation}\label{eq5.2}
\|\exp(-tL)\|_{1 \to \infty}\le Kt^{-n/2}, \ \forall\, t>0.
\end{equation}

Note  that if the self-adjoint contractive semigroup $\exp(-tL)$ on $L^2(X)$ is in addition
uniformly bounded on $L^\infty(X)$,  which includes the case of  a sub-Markovian  semigroup, then
\begin{equation*}
\|\exp(-(t+is)L)\|_{1\to \infty}\le \|\exp(-tL/2)\|_{\infty \to \infty}
\|\exp(-isL)\|_{1 \to \infty} \le C {\ |s|}^{-n/2}.
\end{equation*}
Together with (\ref{eq5.2}) this yields
\begin{equation*}
\|\exp(-(t+is)L)\|_{1\to \infty}\le C\min\{t^{-n/2},|s|^{-n/2}\} \le
C^{\prime}|t+is|^{-n/2}
\end{equation*}
for all $t>0, s\in\mathbb{R}$.
Hence
\begin{equation}\label{eq5.3}
 \|\exp(-zL)\|_{1 \to \infty} \le C |z|^{-n/2}
\end{equation}
for all $\Re z \ge 0$. Of course this estimate  implies
(\ref{eq5.1}) and (\ref{eq5.2}) and the argument above shows that if
semigroup $\exp(-tL)$  is uniformly bounded on $L^\infty(X)$
then it  is equivalent to conjunction  (\ref{eq5.1}) and
(\ref{eq5.2}). It turns out however that  this equivalence holds
without the boundedness assumption on $L^\infty(X)$. This fact will be
used in the next subsection  in which we will not assume  uniform
boundedness of semigroup $\exp(-tL)$ on $L^\infty(X)$.

\medskip
\begin{lemma}\label{le5.1} Suppose that $L$  is a non-negative self-adjoint operator on $L^2(X)$. Then
the dispersive estimate {\rm (\ref{eq5.1})} is equivalent to {\rm
(\ref{eq5.3})}.
\end{lemma}

\noindent{\bf Proof.}  All what we need is to prove that
(\ref{eq5.1}) is enough to get (\ref{eq5.3}) on the positive
half-plane. Fix $f, g \in L^1(X) \cap L^2(X)$ and consider the function
$$H(z) = z^{n/2} \langle \exp(-zL)f, g \rangle.$$
The  analyticity of the semigroup on $L^2$ implies analyticity of
$H$ on the open right half-plane and continuous on the boundary. Now
for $z = is$ with $s \in \RR$,
 the dispersive estimate (\ref{eq5.1}) gives
$$
| H(is) | \le C \| g \|_{1} \| f \|_{1}.
$$
For all $z $ with $\Re z \ge 0$, we have
$$ | H(z) | \le | z |^{n/2}  \| g \|_{2} \| f \|_{2}.$$
Therefore, we can apply the Phragm\'en-Lindel\"of  theorem and conclude that
\begin{equation} \label{eq5.4}
| H(z) | \le C \| g \|_{1} \| f \|_{1}
\end{equation}
for all $z$ with $\Re z \ge 0$. From this and  the   density  of $L^1(X) \cap L^2(X)$
in $L^1(X)$ we obtain the lemma.
\hfill{} $\Box$

\section{From dispersive and Strichartz estimates  to sharp multipliers}\label{sst}
\setcounter{equation}{0}

We continue with the assumption that $(X, d,\mu)$ is   a   metric measure  space. In this section, let us
  start with  the following proposition.
\begin{proposition}\label{th6.1}
Let  $L$ be a non-negative self-adjoint operator on $L^2$. Assume
that $L$ satisfies the dispersive estimate {\rm (\ref{eq5.1})}. Then
for all $1\le p < \frac{2n}{n+2}$ and  all $\lambda \ge 0 $
\begin{equation}\label{eq6.1}
\|d E_{\sqrt{L}}(\lambda)\|_{p \to p' }\le C\lambda^{n(\frac{1}{p} -\frac{1}{p'}) -1},
\end{equation}
where $p'$ is again the conjugate  exponent of $p$.
\end{proposition}

\noindent{\bf Proof.}  We first prove that $L$ satisfies
\begin{equation}\label{eq6.2}
\|F(\sqrt{L} ) \|_{p \to p'} \le C R^{n(\frac{1}{p} -\frac{1}{p'})} \| \delta_R F \|_{1}
\end{equation}
for all bounded $F \in L^1$ with  $\supp F \subseteq [0, R]$ and all $p$ with $1 \le p < \frac{2n}{n+2}$.
 This estimate is  very similar  to ${\rm (ST^1_{p,p'})}$ studied in Part 1.  Note that we do not
  consider here $X$ to be
a doubling space (even a metric $d$ is not needed).

Consider the case where $R = 1$ and fix $F$ with support contained in $[0,1]$.
Set $G (\lambda) = F(\sqrt{\lambda}) e^{\lambda}$. By the inverse Fourier transform, we have (up to a constant)
$$ G(\lambda) = \int_{\RR} \hat{G}(\xi) e^{i \xi \lambda} d \xi.$$
This gives
\begin{equation}\label{eq6.3}
 F(\sqrt{L}) = \int_{\RR} \hat{G}(\xi) \exp(-(1 - i \xi) L) d\xi.
 \end{equation}
This equality follows  immediately  from Fubini's theorem if $\hat{G} \in L^1$.  One may start by proving  (\ref{eq6.3})  for smooth functions
$F_n$ and then use standard approximation arguments to obtain the equality for all $F$ as above.

Now the dispersive estimate together with Lemma \ref{le5.1} imply
that for any $p \in [1,2]$
\begin{eqnarray}\label{eeeee}
\|F(\sqrt{L} ) \|_{p \to p'} &\le&  \int_{\RR} | \hat{G}(\xi) | \|\exp(-(1 - i \xi) L)  \|_{p \to p'}  d\xi\nonumber\\
&\le& C  \int_{\RR} | \hat{G}(\xi) |  (1 + \xi^2)^{-\frac{n}{4}(\frac{1}{p} -\frac{1}{p'})} d\xi \label{eqp7.1} \nonumber\\
&\le& C \| \hat{G} \|_{\infty}  \int_{\RR} (1 + \xi^2)^{-\frac{n}{4}(\frac{1}{p} -\frac{1}{p'})} d\xi.
\end{eqnarray}
Now we note that
$$ \| \hat{G} \|_{\infty}  \le \| G \|_{1} \le C \| F \|_{1} $$
and
$$  \int_{\RR} (1 + \xi^2)^{-\frac{n}{4}(\frac{1}{p} -\frac{1}{p'})} d\xi < \infty$$
for $p < \frac{2n}{n+2}$. This shows (\ref{eq6.2}) when $R = 1$.
Now, for general $R > 0$ and $F$ with support in $[0,R]$ we reproduce the
 previous arguments with the function  $\delta_R F$ and the  operator $L' =
\frac{L}{R^2}$. This leads to (\ref{eq6.2}).
Now we argue as in the proof of  Proposition \ref{prop2.4}.  Fix
$\lambda \ge 0$  and $\eps > 0$ small. We use (\ref{eq6.2}) to obtain
 \begin{eqnarray*}
 \Big\|\varepsilon^{-1} E_{\sqrt{L}}(\lambda-\varepsilon, \lambda+\varepsilon]   \Big\|_{p\to p'}
 &=&\varepsilon^{-1}\Big\|  1\!\!1_{(\lambda-\varepsilon, \lambda+\varepsilon] } (\sqrt{L}) \Big\|_{p\to p'} \\
 &\leq&
C \varepsilon^{-1}(\lambda+\varepsilon)^{n({1\over p}-{1\over {p'}})}
\big\| \chi_{({\lambda-\varepsilon\over \lambda+\varepsilon}, \, 1]}\big\|_{1} \\
 &\leq& C(\lambda+\varepsilon)^{n({1\over p}-{1\over p'})-1}.
 \end{eqnarray*}
Letting $\eps \to 0$ we obtain
$$  \Big\|  dE_{\sqrt{L}}(\lambda)   \Big\|_{p\to p'} \le C\lambda^{n({1\over p}-{1\over p'})-1},$$
which implies the estimate of the proposition. \hfill $\Box$

\begin{coro}\label{der}  Suppose that  $(X, d,\mu)$
  satisfies  the doubling property (\ref{eq2.2}) and  there exists a   positive constant $C > 0$ such that
$V(x, r)\leq C r^n$
for every $x\in X $ and $r>0$.  Assume that  $L$ satisfies  the
 finite speed propagation property {\rm (FS)} and
  the dispersive estimate {\rm (\ref{eq5.1})}.
Fix $p \in [1, \frac{2n}{n+ 2}]  
 $ and suppose next that $G$ is a distribution, $\supp G \subset [0,R]$ and
that  $\int_{\RR} | \hat{G}(\xi) |  (1 + \xi^2/R^2)^{-\frac{n}{4}(\frac{1}{p} -\frac{1}{p'})} d\xi < \infty$.
Then the operator $G(L)$ is well defined as an operator acting from $L^p$ to $L^{p'}$ and
$$
\|G(L)\|_{p\to p'} \le C\int_{\RR} | \hat{G}(\xi) |  (1 + \xi^2/R^2)^{-\frac{n}{4}(\frac{1}{p} -\frac{1}{p'})} d\xi.
$$
\end{coro}

\noindent{\bf Proof.} Corollary~\ref{der} follows from estimates \eqref{eeeee}.
\hfill $\Box$

\medskip

Proposition \ref{th6.1} does not yield the optimal results for the standard Laplace operator.
However in the abstract setting we can  include the endpoint $p = \frac{2n}{n+2}$ when $n > 2$ in the following way.   We
start with the Strichartz estimate (\ref{str-end}) and repeat  the previous proof to get
\begin{eqnarray*}
\|F(\sqrt{L}) \|_{2 \to \frac{2n}{n-2}} &\le& \int_\RR | \hat{G}(\xi) |
 \|\exp(-(1 - i \xi) L)  \|_{2 \to \frac{2n}{n-2}}   d\xi\\
&\le& \| \hat{G} \|_{2} \left( \int_\RR \|\exp(i \xi L) \|_{2 \to \frac{2n}{n-2}}^2 d\xi \right)^{1/2}\\
&\le& C \| F \|_{2}
\end{eqnarray*}
for all $F$ with support in $[0,1]$. For $F$ supported in $[0,R]$ we apply the
previous estimate with  $\delta_RF$ and $L' = \frac{L}{R^2}$ to get
$$ \|F(\sqrt{L}) \|_{2 \to \frac{2n}{n-2}} \le C R \| \delta_R F \|_2.$$
Therefore,
\begin{equation}\label{ST-str}
 \|F(\sqrt{L} )\|_{\frac{2n}{n+2} \to \frac{2n}{n-2}} \le C R^2 \| \delta_R F \|_2^2.
 \end{equation}
 Similar arguments as above give (\ref{eq6.1}) for $p = \frac{2n}{n+2}$.

 We can now extend this easily to $p <  \frac{2n}{n+2}$ if the smoothing property (\ref{eq5.2})
 is satisfied. More precisely,
 fix  $p <   \frac{2n}{n+2}$ and assume that
 \begin{equation}\label{eq5.2p2}
\|\exp(-tL)\|_{p \to \frac{2n}{n+2} }\le Kt^{-\frac{n}{2} (\frac{1}{p} - \frac{n+2}{2n} )}, \ \forall\, t>0.
\end{equation}
 We introduce as before
$G_R (\lambda) = (\delta_RF)( \sqrt{\lambda} ) e^\lambda$ for $F$ supported in $[0,R]$. Then for $q = p'$ we have
 \begin{eqnarray*}
 \| F(\sqrt{L}) f \|_{q} &= & \|\int_\RR  \widehat{G_R}(\xi)
 \exp(-( \frac{1}{R^2} - i \frac{\xi}{R^2})L) f d\xi \|_{q}\\
 &\le&   \| \exp(-\frac{L}{R^2}) \|_{\frac{2n}{n-2} \to q} \int_\RR  | \widehat{G_R}(\xi) |
  \|\exp(-( i \frac{\xi}{R^2})L) f \|_{{\frac{2n}{n-2} }} d\xi\\
 &\le& C R^{n( \frac{1}{2} -\frac{1}{q})  -1} \| G_R \|_{2} \left( \int_{\RR}
  \|\exp(-( i \frac{\xi}{R^2})L) f \|_{{\frac{2n}{n-2} }}^2 d\xi \right)^{1/2}\\
 &\le& C R^{n( \frac{1}{2} -\frac{1}{q})} \| \delta_R F \|_{L^2(\RR)} \| f \|_{L^2(X)},
 \end{eqnarray*}
 where  we used the Strichartz estimate  (\ref{str-end}) to obtain the last inequality.
As in the last proposition, this gives  the Stein-Tomas restriction (\ref{eq6.1}) for all
 $p \in [1, \frac{2n}{n+2}]$.  We  have proved
 \begin{proposition}\label{pro6.1}
Let  $L$ be a non-negative self-adjoint operator on $L^2$. Assume
that $L$ satisfies the Strichartz estimate {\rm (\ref{str-end})}   for some $n > 2$.
Fix $p$ such that $1\le p \le  \frac{2n}{n+2}$
and assume that  the smoothing property
 {\rm (\ref{eq5.2p2}) }  is satisfied.
Then for all   $\lambda \ge 0 $
\begin{equation}\label{eq6.11}
\|d E_{\sqrt{L}}(\lambda)\|_{p \to p' }\le C\lambda^{n(\frac{1}{p} -\frac{1}{p'}) -1},
\end{equation}
where $p'$ is again the conjugate  exponent of $p$.
\end{proposition}

We mentioned above  that  the dispersive estimate implies the  Strichartz estimate (\cite{KT}).
For this reason we formulate the results below for  the case where $L$
satisfies the endpoint Strichartz estimate.

\begin{theorem}\label{th6.2} Suppose that  $(X, d,\mu)$
  satisfies  the doubling property (\ref{eq2.2}) and  there exists a   positive constant $C > 0$ such that
$V(x, r)\leq C r^n$
for every $x\in X $ and $r>0$.
Assume that  $L$ satisfies  the
 finite speed propagation property {\rm (FS)} and the Strichartz  estimate {\rm (\ref{str-end})}
 with the same $n$ as in the doubling property.  Assume also that $n > 2$. Fix $p \in [1, \frac{2n}{n+ 2}]  
 $ and assume that {\rm (\ref{eq5.2p2})} holds.      For every
 compactly supported bounded function $F$ such that
 $$ \|   F \|_{W^{\beta,2}} < \infty$$
 for some $\beta > n (\frac{1}{p} - \frac{1}{2})$, the operator
 $F(t\sqrt{L})$ is bounded on $L^p$ for all $t>0$ and
 $$\sup_{t > 0} \| F(t\sqrt{L}) \|_{p\to p} \le C \| F \|_{W^{\beta,2}}.$$
 \end{theorem}

\noindent{\bf Proof.}  If $F$ is supported in $[0, R]$ then by
Proposition \ref{pro6.1}, we have
\begin{eqnarray*}
\| F(\sqrt{L}) \|_{p\to p'} &=& \| \int_0^R F(\lambda) dE_{\sqrt{L}}(\lambda) \|_{p\to p'} \\
&\le& C \int_0^R |F(\lambda) | \lambda^{n(\frac{1}{p} -\frac{1}{p'}) -1} d\lambda\\
&\le& C R^{n (\frac{1}{p} -\frac{1}{p'})} \| \delta_R F \|_{1}.
\end{eqnarray*}
Hence by the $T^{\ast} T$ argument
$$ \| F(\sqrt{L}) \|_{p \to 2} \le C R^{n(\frac{1}{p} -\frac{1}{2})} \| \delta_R F \|_{2}.$$
Combining this with our assumption on the volume yields ${\rm (ST^{2}_{p, 2})}$.
We then apply Theorem \ref{th3.1} and obtain the result for $p \in [1, \frac{2n}{n+ 2}]$.
   \hfill $\Box$

\medskip

We have seen that the assumptions of the previous theorem imply ${\rm (ST^{2}_{p, 2})}$, we
can then apply Theorem \ref{th4.1} to obtain endpoint
estimate for Bochner-Riesz means. In addition, by applying  Theorem \ref{th31.1} we obtain
under the assumptions of the previous theorem the following result.
\begin{theorem}\label{th66.1}
Fix $p \in [1,   \frac{2n}{n+2}].$
For any even bounded Borel function $F$ such that
$\sup_{t>0}\|\eta\, \delta_tF\|_{W^{\beta, 2}}$ $<\infty $ for some  $\beta>\max\{n(1/p-1/2),1/2\}$ and some
non-trivial function $\eta \in C_c^\infty(0,\infty)$, the operator
$F(\SL)$ is bounded on $L^r(X)$ for all $r \in (p, p')$.
In addition,
\begin{eqnarray*}
   \|F(\SL)  \|_{r\to r}\leq    C_\beta\Big(\sup_{t>0}\|\eta\, \delta_tF\|_{W^{\beta, 2}}
   + |F(0)|\Big).
\end{eqnarray*}
\end{theorem}

\medskip

\noindent
\begin{remark}\label{Remark 7.5} In the general  setting of doubling spaces, we can replace the dispersive
 estimate {\rm (\ref{eq5.1})} by
\begin{equation}\label{disdo}
\| P_{B(x,r)} \exp\left(isL\right) P_{B(x,r)} \|_{p \to p'} \le
 C V(x,r)^{\frac{1}{p'} -\frac{1}{p}} \left( \frac{r}{\sqrt{|s|}} \right)^{n(\frac{1}{p} -\frac{1}{p'} )}.
\end{equation}
Here the constant $n$ is the same as in the doubling condition.
The arguments  in the proof of (\ref{eq6.2})  show that  for $p < \frac{2n}{n+2}$,
$$\| P_{B(x,r)} |F|^2 P_{B(x,r)} \|_{p \to p'} \le C V(x,r)^{\frac{1}{p'} -\frac{1}{p}}
 (Rr)^{n(\frac{1}{p} -\frac{1}{p'} )} \| \delta_R F \|_{2}^2
 $$
for all $F$ supported in $[0,R]$. The $T^{\ast} T$ argument gives ${\rm (ST^2_{p, 2})}$
 for all $p <  \frac{2n}{n+ 2}$.
Therefore we obtain the same conclusion as in  Theorems \ref{th6.2}  and \ref{th66.1} for $L$ 
satisfying   property {\rm (FS)} and (\ref{disdo})  on any doubling space.

 Note also that under these two  assumptions, we obtain from Theorem \ref{th4.1}  endpoint estimates for
  the Bochner-Riesz mean  $S^{\delta(p)}_R(L)$ on $L^p$ ($p\in [1, \frac{2n}{n+2})$)  for
 $\delta(p) =\max\big\{0,n\big|\frac{1}{p}-\frac{1}{2}\big|-\frac{1}{2}\big\}$.
 \end{remark}
 
 In the following proposition we consider local dispersive estimate. 
 
\medskip

\begin{proposition}\label{prop7.7}
Let  $L$ be a non-negative self-adjoint operator on $L^2$. Assume
that $L$ satisfies the local dispersive estimate 
$$
\|\exp(itL) \|_{1 \to \infty }\le Ct^{-n/2}
$$
for all $0<|t|\leq 2$. Assume also that 
$$ \|\exp(-tL) \|_{\infty \to \infty }\le C \ \  \ {\rm and}\ \ \ \|\exp(-tL) \|_{1\to \infty }\le Ct^{-n/2}
$$
for all $t>0$.
Then
for all $1\le p < \frac{2n}{n+2}$ and  all $k>0 $
\begin{equation}\label{eq6}
\|E_L{[k, k+1)}  \|_{p \to p' }\le C(1+k)^{\frac{n}{2}(\frac{1}{p} -\frac{1}{p'}) -1},
\end{equation}
where $p'$ is again the conjugate  exponent of $p$.
\end{proposition}

\noindent{\bf Proof.}
Set $\widehat{G_k}(\xi)=3\pi (1-|\xi|)_+^2e^{-ik \xi}$ so that ${G_k}(\lambda)=6 (\lambda -k)^{-2} -6{\sin (\lambda-k)\over 
(\lambda -k)^{3}}.$
Note that there exists a positive constant $c$ such that ${G_k}(\lambda)e^{-\lambda/(2k)} \ge c\chi_{[k, k+1)}(\lambda)$ for 
all $\lambda \in [k, k+1)$. So 
$$
c^2\|E_L{[k, k+1)}f\|_{2 }^2  \le \|{G_k}{(L)}e^{-L/(2k)}f\|_{2 }^2.
$$
It follows from the above inequality and the T*T argument that 
$$
c^2\|E_L{[k, k+1)}\|_{p \to p' }  \le \| |G_k|^2({L})e^{-L/k}\|_{p\to p' }.
$$
Next  (\ref{eq5.3}) holds for $|{\rm Im} z|\leq 2$. Hence
\begin{eqnarray*}\nonumber
\||G_k|^2({L})e^{-L/k} \|_{p \to p'} &\le&  \int_{-2}^2 | \widehat{G_k}\ast \widehat{G_k} | \|\exp(-(1/k - i \xi) L)  \|_{p \to p'}  d\xi\\
&\le& C  \int_{-2}^2  (k^{-2} + \xi^2)^{-\frac{n}{4}(\frac{1}{p} -\frac{1}{p'})} d\xi \label{eqp7.1} \\
&\le& C  \int_{\RR} (k^{-2}+ \xi^2)^{-\frac{n}{4}(\frac{1}{p} -\frac{1}{p'})} d\xi\\
& \le& C
(1+k)^{\frac{n}{2}(\frac{1}{p} -\frac{1}{p'})-1}.
\nonumber
\end{eqnarray*}
This proves   estimate (\ref{eq6}).
\hfill $\Box$

\bigskip

\part{Applications}

\section{Standard Laplace operator and compact manifolds}\label{sec51}
\setcounter{equation}{0}

As mentioned in the introduction, the restriction estimates  ${\rm (R_p)}$
 for standard Laplace operator on ${\mathbb R}^n$ are valid for
$1\leq p\leq {2(n+1)/(n+3)}.$  As a consequence of Theorem~\ref{th4.1}, we obtain
alternative proof of Theorem 1.1 of \cite{T2} by Tao, Theorem 1 \cite{C2} and main result of
\cite{C1} described by  M. Christ. These results can be stated in the following way:

 \begin{proposition}
 For all $n\geq 2$ and $1\leq p\leq {2(n+1)/(n+3)}$, the operator
 $(I-\Delta)_+^{\delta(p)}$ is of weak-type
 $(p,p).$
 \end{proposition}

\medskip

\noindent
{\bf Proof.} This result is straightforward from   Proposition~\ref{prop2.4}
and Theorem~\ref{th4.1}.
  \hfill{}$\Box$

 \bigskip

Similarly, using  Theorem~\ref{th4.2} one can obtain
alternative proof of Theorem 1.2 of \cite{T2}. Our proof shows that
this  result  holds for all operators on compact manifolds  which satisfy property  (FS) and condition ${\rm(S_p)}$ as in the following
proposition.

 \begin{proposition} Suppose the operator   $L$ satisfies  ${\rm (FS)}$ and condition ${\rm (S_p)}$ for some $1\leq p\leq {2(n+1)/(n+3)}$.
Then the operator
 $(I-L/R^2)_+^{\delta(p)}$ is of weak-type
 $(p,p)$ uniformly in $R.$
 \end{proposition}

\medskip

\noindent {\bf Proof.} This result is straightforward from
Proposition~\ref{prop3.8} and Theorem~\ref{th4.2}.
  \hfill{}$\Box$

  \medskip

In both cases of compact manifolds with or without boundaries,
examples which satisfy condition ${\rm (S_p)}$   are described in
\cite{Sog1, Sog4} by C.D. Sogge.

We mentioned here  endpoint Bochner-Riesz summability results.
From Theorems~\ref{th31.1} and ~\ref{th31.2} we have more
general spectral multiplier results for the operators
considered in the previous propositions.

 \medskip

\section{Asymptotically conic manifolds}\label{sec52}
\setcounter{equation}{0}

Scattering manifolds  or  asymptotically conic manifolds
 are defined as  the interior of a compact manifold with boundary $M$, and the metric $g$
 is smooth on $M^\circ$ and has the form
\[g=\frac{dx^2}{x^4}+\frac{h(x)}{x^2} \]
in a collar neighbourhood near $\pl M$, where $x$ is a smooth boundary defining function
for $M$ and $h(x)$ a smooth one-parameter family of metrics on $\pl M$; the function $r:=1/x$
near $x=0$ can be thought of as a radial coordinate near infinity and the metric there is
asymptotic to the exact metric cone $((0,\infty)_r\x \pl M, dr^2+r^2h(0))$.

 In this subsection we consider  the following classical operators:
\begin{itemize}
\item Schr\"odinger operators, i.e. $-\Delta + V$ on $\RR^n$,
where $V$ smooth and decaying sufficiently at infinity;

\item The Laplacian with respect to  metric perturbations of
the flat metric on $\RR^n$, again decaying sufficiently at infinity;

\item The Laplacian on asymptotically conic manifolds.
\end{itemize}

\begin{proposition}\label{prop8.1} Let $(M,g)$ be an
asymptotically conic manifold of dimension $n \geq 3$, and let $x$ be a smooth boundary defining function of
$\pl M$. Let $L:= - \Delta+V$ be a Schr\"odinger operator
with $V\in x^3C^\infty(M)$ and assume that $L$ has no $L^2$-eigenvalues and that $0$ is not a resonance.
Then
\begin{itemize}
\item[(i)]
  For any $\lambda_0 > 0$ there exists a constant $C>0$ such that the spectral measure $dE(\lambda)$
for $\sqrt{L}$ satisfies

\begin{equation}\label{eq8.1}
\| dE_{\sqrt{L}}(\lambda) \|_{L^p(M) \to L^{p'}(M)} \leq C
\lambda^{n({1\over p}-{1\over p'}) - 1}
\end{equation}

\noindent
for $1 \leq p \leq 2(n+1)/(n+3)$ and $0 < \lambda \leq \lambda_0$.

\item[(ii)]  If $(M,g)$ is nontrapping, then there exists $C>0$ such that
 \eqref{eq8.1}  holds for all $\lambda > 0$.

\end{itemize}
\end{proposition}

\smallskip

 Proposition~\ref{prop8.1} was proved in \cite[Theorem 1.2]{GHS}. This proposition
 has useful consequence to  establish  the convergence of the Riesz means up to the critical exponent
$\delta(p)=\max\big\{0,n\big| {1/p}- {1/2}\big|- {1/2}\big\}$ for all $1\leq p\leq 2(n+1)/(n+3)$.

\medskip

\begin{coro} Let $(M,g)$   be nontrapping and the operator $L$ satisfies  all assumptions of
Proposition~\ref{prop8.1}.   Let $1\leq p\leq 2(n+1)/(n+3)$. Then
\begin{itemize}
\item[(i)]
 $S_R^{\delta(p)} (L)$ is of weak-type $(p,p)$ uniformly in $R$.

\item[(ii)] For any even bounded Borel
function $F: [0, \infty) \to \CC$ such that
$\sup_{t>0}\|\eta\, \delta_tF\|_{W^{\beta, 2}}<\infty $ for some  $\beta>\max\{n(1/p-1/2),1/2\}$
 and some non-trivial function $\eta \in C_c^\infty(0,\infty)$, the operator
$F(\SL)$ is bounded on $L^r(X)$ for all $p<r<p'$ with
\begin{eqnarray*}
   \|F(\SL)  \|_{r\to r}\leq    C_\beta\Big(\sup_{t>0}\|\eta\, \delta_tF\|_{W^{\beta, 2}}
   + |F(0)|\Big).
\end{eqnarray*}
\end{itemize}
\end{coro}

\noindent {\bf Proof.} This result is a straightforward application of
Propositions~\ref{prop8.1}, ~\ref{prop2.4} and
Theorems~\ref{th4.1},  ~\ref{th31.1}.
  \hfill{}$\Box$

\bigskip
\section{Schr\"odinger operators with rough potentials}\label{sec54}
\setcounter{equation}{0}
This section is devoted to  Shr\"odinger operators  $-\Delta + V$ for which we prove new spectral multiplier results.

\medskip

\noindent
{\bf 10.1.\ Schr\"odinger operators with   inverse-square potential.}\
 We start with inverse square potentials, that is $V(x) = \frac{c}{| x |^2}$.  Fix $n > 2$ and assume that
$  -{(n-2)^2/4}< c $.     Define by quadratic form method
$L  = -\Delta + V$ on $L^2(\RR^n, dx)$.
The classical Hardy  inequality
\begin{equation}\label{hardy1}
- \Delta\geq  {(n-2)^2\over 4}|x|^{-2},
\end{equation}
shows that  for all $c > -{(n-2)^2/4}$,  the self-adjoint operator $L$  is   non-negative.
Set  $p_c^{\ast}=n/\sigma$, $\sigma= \max\{ (n-2)/2-\sqrt{(n-2)^2/4+c}, 0\}$. If $c \ge 0$ then the semigroup $\exp(-tL)$
is pointwise bounded by the Gaussian semigroup and hence act on all $L^p$ spaces with $1 \le p \le \infty$.  If $ c < 0$, then $\exp(-tL)$
acts as a uniformly bounded semigroup on  $L^p(\RR^n)$ for
$ p \in ((p_c^{\ast})', p_c^{\ast})$ and the range $((p_c^{\ast})', p_c^{\ast})$ is optimal (see for example  \cite{LSV}).

It is proved in \cite{BPST} that the solution $u(t) = e^{-itL}f$  of the corresponding Schr\"odinger  equation
$$i \partial_t u  + L u  = 0, \, \, u(0) = f$$
satisfies Strichatrz estimate (\ref{str-end}).  The smoothing property (\ref{eq5.2p2}) is proved in \cite{Ass}.
Therefore, we obtain from
Proposition \ref{pro6.1} that $L$ satisfies restriction estimate ${\rm (R_p)}$
 for all $p \in ((p_c^{\ast})',  \frac{2n}{n+2}]$. If $c \ge 0$, then $p = (p_c^{\ast})' = 1$ is included.
From this and  Theorems~\ref{th4.1} and ~\ref{th6.2}
we obtain
\begin{theorem}\label{th10.1}
Suppose that $n > 2$  and $ -{(n-2)^2/4} < c $  and that
 $p\in ((p_c^{\ast}){'}, 2n/(n+2)]$ where $p_c^{\ast}=n/\sigma$ and $\sigma= \max\{ (n-2)/2-\sqrt{(n-2)^2/4+c}, 0\}$
 and $(p_c^{\ast}){'}$ its dual exponent.
Then
\begin{itemize}
\item[(i)] $S_R^{\delta(p)} (L)$
is of weak-type $(p,p)$ uniformly in $R$.

\item[(ii)]
  For any even bounded Borel
function $F: [0, \infty) \to \CC$ such that
$\sup_{t>0}\|\eta\, \delta_tF\|_{W^{\beta, 2}}<\infty $ for some  $\beta>\max\{n(1/p-1/2),1/2\}$
 and some non-trivial function $\eta \in C_c^\infty(0,\infty)$, the operator
$F(\SL)$ is bounded on $L^r(X)$ for all $p<r<p'$ with
\begin{eqnarray*}
   \|F(\SL)  \|_{r\to r}\leq    C_\beta\Big(\sup_{t>0}\|\eta\, \delta_tF\|_{W^{\beta, 2}}
   + |F(0)|\Big).
\end{eqnarray*}
\end{itemize}
\end{theorem}

\medskip

\noindent
{\bf 10.2.\ Scattering operators.}\
Assume now that $n= 3$ and   $V$ is a real-valued measurable function such that
\begin{eqnarray}\label{eq10.1}
\int_{{\mathbb R}^6} \frac{|V(x)|\, |V(y)|}{|x-y|^2}dx dy < (4\pi)^2 \quad \ \ \mbox{and} \quad \ \
\sup_{x\in{\mathbb R}^3}
\int_{{\mathbb R}^3} \frac{ |V(y)| }{|x-y|}dy< 4\pi.
\end{eqnarray}
The following proposition is a consequence of Proposition \ref{th6.1}
and the main result in Rodnianski and Schlag \cite{RS} which gives  the dispersive
estimate for $\exp(it (-\Delta + V))$ on $\RR^3$.

\begin{proposition}\label{prop10.1}
Suppose that $L =-\Delta  + V$ on ${\mathbb R}^3$
 with a real-valued    $V$ which satisfies   \eqref{eq10.1}. Then $L$ satisfies ${\rm (R_p)}$
  for all $1 \le p < 6/5$.
\end{proposition}

In the special case  $p=1$, Proposition~\ref{prop10.1} was obtain in \cite[Theorem 7.15]{DOS}
for compactly supported function $V\geq 0$   which satisfies   \eqref{eq10.1}.
The following result is a consequence of Theorem~\ref{th4.1}, \ref{th31.1} and
Proposition~\ref{prop10.1}.

\begin{coro}\label{coro10.2}  Suppose that $L =-\Delta + V$ on ${\mathbb R}^3$
and that $V$ satisfies assumption of Proposition~\ref{prop10.1}.
Assume also that  $1\leq p < 6/5$. Then
\begin{itemize}
\item[(i)]
 $S_R^{\delta(p)} (L)$ is of
weak-type $(p,p)$ uniformly in $R$.

\item[(ii)]
  For any even bounded Borel
function $F: [0, \infty) \to \CC$ such that
$\sup_{t>0}\|\eta\, \delta_tF\|_{W^{\beta, 2}}<\infty $ for some  $\beta>\max\{3(1/p-1/2),1/2\}$
 and some non-trivial function $\eta \in C_c^\infty(0,\infty)$, the operator
$F(\SL)$ is bounded on $L^r(X)$ for all $p<r<p'$.
\end{itemize}
\end{coro}

If $n \ge 3$ and   potential $V \in W^{s,2}(\RR^n)$  for some $s > \frac{n}{2} -1$ and has fast decay,   Bourgain \cite{B2} proved the dispersive estimate
for $\exp(it (-\Delta + V))$ .
Our results apply for $L = -\Delta + V$ and allow to obtain  sharp spectral multiplier results.  We also refer to
Rodnianski and Schlag  \cite{RS} for more  references on dispersive  estimates
for Schr\"odinger operators.

We also mention that Strichartz estimates are proved for a class of elliptic operators with variable coefficients
 by J. Marzuola, J. Metcalfe and D. Tataru \cite{MMT}  (see Theorem 1.20).
Therefore, the same reasoning as for the Theorem \ref{th10.1} allows us to obtain sharp spectral multipliers and endpoint
Bochner-Riesz summability for these elliptic operators.
\bigskip

\medskip

\noindent
{\bf 10.3.\ The harmonic oscillator.} In this section we focus  on Schr\"odinger operators such as the  harmonic oscillator 
$-\Delta + |x|^2$ on $L^2(\RR^n)$
for $n\ge 2$. As in    \cite{KoT} we can also consider 
Schr\"odinger operators $
L= - \Delta + V$
with a  positive potential $V$ which satisfies the following condition 
\begin{equation}\label{eq111.01}
V \sim |x|^2, \quad |\nabla V| \sim |x|, \quad |\partial_x^2 V| \le 1.
\end{equation}
We apply Theorems~\ref{th3.2}~and~\ref{th31.2}  and the results from \cite{KoT} to prove 
sharp results on  Bochner-Riesz summability and singular spectral multipliers for $L$.   
 Bochner-Riesz  summability results for the harmonic oscillator were  obtained  first  by Kardzhov  \cite{Kar}. Here  we describe an 
alternative proof. The corresponding  singular integral multiplier is a new result. The following theorem is the main goal of this section. 
\begin{theorem}\label{coro12.01} Assume that potential $V$ satisfies condition \eqref{eq111.01} and
set $L= -\Delta +V$.   Let $1\leq p\leq 2n/(n+2)$. Then
\begin{itemize}
\item[(i)]
For any even function $F$ such that $\supp F \subseteq [-1,1]$ and $\|F\|_{W^{\beta,2}}<\infty$
 for some $\beta>\max\{n(1/p-1/2),1/2\}$, the operator $F(t\sqrt{L})$ is bounded on $L^p(X)$ for all $t>0$ and
$$
\label{eq3.2} \sup_{t>0}\|F(t\sqrt{L})\|_{p\to p} \leq
C\|F\|_{W^{\beta,2}}.
$$
\item[(ii)] For any even bounded Borel
function $F: [0, \infty) \to \CC$ such that
$\sup_{t>0}\|\eta\, \delta_tF\|_{W^{\beta, 2}}<\infty $ for some  $\beta>\max\{n(1/p-1/2),1/2\}$
 and some non-trivial function $\eta \in C_c^\infty(0,\infty)$, the operator
$F(\SL)$ is bounded on $L^r(X)$ for all $p<r<p'$.
In addition,
\begin{eqnarray*}
   \|F(\SL)  \|_{r\to r}\leq    C_\beta\Big(\sup_{t>0}\|\eta\, \delta_tF\|_{W^{\beta, 2}}
   + |F(0)|\Big).
\end{eqnarray*}
\end{itemize}
\end{theorem}

\noindent
{\bf Proof.}  It follows from Theorem 4 in  \cite{KoT} that  for all $\lambda \ge 0$  and 
 all $1\leq p\leq 2n/(n+2)$
\begin{eqnarray}\label{chen111}
\|E_L[\l^2,\l^2+1)\|_{p\to 2} \leq C(1+\l)^{n({1\over p}-{1\over 2})-1}.
\end{eqnarray}
Take a function $F$ with support in $[-N, N]$.  
 We have as in the proof of Proposition~\ref{prop3.8}
\begin{eqnarray}\label{e10.555}
\|F(\SL) \|_{p\to2}^2 &\le&   \sum_{\ell = 1}^{N^2} \| E_{{\SL}}[\frac{\ell -1}{N}, \frac{\ell}{N} ) F(\SL) \|_{p\to2}^2 \nonumber\\
&\le&  \sum_{\ell = 1}^{N^2}  \sup_{\lambda \in [\frac{\ell -1}{N}, \frac{\ell}{N} ) } |F(\lambda)|^2 \| E_{{\SL}}[\frac{\ell -1}{N}, \frac{\ell}{N} ) \|_{p\to2}^2\nonumber \\
&=& \sum_{\ell = 1 }^{N^2}  \sup_{\lambda \in [\frac{\ell -1}{N}, \frac{\ell}{N} ) } |F(\lambda)|^2 \| E_{{L}}[(\frac{\ell -1}{N})^2,  (\frac{\ell}{N})^2 ) \|_{p\to2}^2.
\end{eqnarray}
Now  we observe that for all $\ell=1,2,\cdots, N^2$
\begin{eqnarray*}
 \big\|E_{{L}}[(\frac{\ell -1}{N})^2,  (\frac{\ell}{N})^2 ) \big\|_{p\to2}^2 &\leq&  
 \big\|E_{{L}}[(\frac{\ell -1}{N})^2,  (\frac{\ell-1}{N})^2+2 ) \big\|_{p\to2}^2.
\end{eqnarray*}
This, in combination with  (\ref{chen111}) and (\ref{e10.555}),  shows  that for $1\leq p\leq  2n/(n+2)$
\begin{eqnarray*}
\|F(\SL) \|_{p\to2}^2 &\le& \sum_{\ell = 1 }^{N^2} \sup_{\lambda \in [\frac{\ell -1}{N},
\frac{\ell}{N} ) } |F(\lambda)|^2
\Big(2+\frac{\ell-1}{N}\Big)^{2n({1\over p}-{1\over 2})-2}\\
&\le&C N^{2n({1\over p}-{1\over 2})}\frac{1}{N^2}\sum_{\ell = 1 }^{N^2}
\sup_{\lambda \in [\frac{\ell -1}{N}, \frac{\ell}{N} ) }
|F(\lambda)|^2 \\
&\le&C N^{2n({1\over p}-{1\over 2})}\frac{1}{N^2}\sum_{\ell = 1 }^{N^2}
\sup_{\lambda \in [\frac{\ell -1}{N^2}, \frac{\ell}{N^2} ) }
|F(N\lambda)|^2.
\end{eqnarray*}
This proves  ${\rm (SC^{2,\kappa}_{p, 2})}$ for $\kappa=2$ and $p$
such that  $1\leq p\leq 2n/(n+2)$. It remains to show ${\rm (AB_p)}$
and then apply Theorems \ref{th3.2} and \ref{th31.2}. Now condition
${\rm (AB_p)}$ follows from the following lemma. \hfill{}$\Box$

 \begin{lemma}\label{lem12.01}
   Let $L=-\Delta +V$, where $V  \in
   L^1_{\mbox{\small{\rm loc}}}(\RR^n)$ and  $V   \ge  0$.
   Suppose that  for some $\kappa>0$ and any $\epsilon >0$
     \begin{equation}\label{Vi}
     \int_{\RR^n} (1+V(x))^{n(1-\kappa)/2-\epsilon}d x < \infty.
     \end{equation}
   Then condition ${\rm (SC^{2,\kappa}_{p, 2})}$ implies  ${\rm (AB_p)}$.
 \end{lemma}

 The proof of Lemma \ref{lem12.01} is a straightforward modification of the proof
 of Lemma 7.9 of \cite{DOS} so we skip it here.
\hfill{}$\Box$

\bigskip

\section{Operators $\Delta_n+  {c\over r^{2}}$ acting on $L^2((0,\infty), r^{n-1}dr)$}
\setcounter{equation}{0}

In this section  we consider a class of Schr\"odinger operators
on $L^2((0,\infty), r^{n-1} dr)$. These operators
generate semigroups but do not  have the classical Gaussian upper bound for the  heat kernel.

Fix  $n > 2$ and   $ c>  -{(n-2)^2/4} $ and consider the space $L^2((0,\infty), r^{n-1}dr)$.
For  $f,g \in C_c^\infty(0,\infty)$
we define the quadratic form

\begin{equation}
Q_{n,c}^{(0,\infty)}  (f,g)=\int^{\infty}_{0}f'(r)g'(r)r^{n-1} dr
+\int_0^\infty  \frac{c}{r^2} f(r)g(r) r^{n-1}dr.
\label{eq9.1}
\end{equation}

\noindent Using the Friedrichs extension one can define the operator  $L_{n,c}=
\Delta_n+c/{r^2}$ as the unique self-adjoint operator corresponding
to $ Q_{n,c}^{(0,\infty)} $, acting on $L^2((0,\infty), r^{n-1}dr)$. In the sequel we will write
$L$ instead of $L_{n,c}$, which is formally given by the following formula
$$
Lf=(\Delta_n+\frac{c}{r^2} ) f=-\frac{d^2}{dr^2}f-\frac{n-1}{r}\frac{d}{dr}f +\frac{c}{r^2}f.
$$

\noindent
The classical Hardy  inequality (\ref{hardy1})
shows that  for all $c > -{(n-2)^2/4}$,  the self-adjoint operator  $L$
is non-negative. Such operators can be seen as radial  Schr\"odinger operators
with inverse-square potentials.  It follows by Theorem 3.3 of \cite{CouS} that
$L$ satisfies  Davies-Gaffney estimate, which in turns implies property (FS).

Now for $   -{(n-2)^2/4}<c<0,$ we set $p_c^{\ast}=n/\sigma$  where $\sigma=(n-2)/2-\sqrt{(n-2)^2/4+c}$
and    $(p_c^{\ast}){'}$ its dual exponent.   Note that $2<{2n\over n-2}<p_c^{\ast}$.
Liskevich, Sobol
and Vogt \cite{LSV}  proved that for  all $t>0$ and all
$p\in ((p_c^{\ast}){'}, p_c^{\ast})$,

$$
\|e^{-tL}\|_{p\to p}\leq C.
$$
They also proved that range $((p_c^{\ast}){'}, p_c^{\ast})$ is optimal
 and  that for all $p\not\in ((p_c^{\ast}){'}, p_c^{\ast})$, the semigroup does not even act on
$L^p((0,\infty), r^{n-1}dr)$
(see also  \cite{DS, CouS, HS}).

 \medskip

\begin{proposition}\label{prop9.1}
Suppose that $n > 2$  and $  -{(n-2)^2/4}<c $. For $c<0$, set
 $p\in \big((p_c^{\ast}){'}, {2n\over n+1}\big)$ where $p_c^{\ast}=n/\sigma$ and $\sigma=(n-2)/2-\sqrt{(n-2)^2/4+c}$
 and $(p_c^{\ast}){'}$ its conjugate exponent. For $c\ge 0$, set $p\in \big[1, {2n\over n+1}\big).$
Then  for any $R>0$ and all Borel functions $F$ such that $\supp F \subset [0,R]$,
\begin{equation}\label{eq9.2}
\big\|F(\SL)P_{B(x_B, r_B)} \big\|_{p\to 2} \leq CV(x_B,
r_B)^{{1\over 2}-{1\over p}} \big(Rr_B)^{n({1\over p}-{1\over
2})}\big\|\delta_RF\big\|_{2}
\end{equation}
  for   all
  $x_B\in {\mathbb R}_+$ and $r_B\geq 1/R$.
\end{proposition}

\medskip

\noindent
{\bf Proof.}\, In \cite{HS} the explicit formula for the resolvent of the operator $L=\Delta_n+\frac{c}{r^2}$
 is described.
Based on this formula we calculate explicitly the spectral projections $dE_{\sqrt{L}}(\lambda)$
(This calculation was shown to us by Andrew Hassell).
Define the number $n' = n'(n, c)$ to be the maximum positive root of the equation
$(n'/2 - 1)^2 = (n/2 - 1)^2 + c$ so that $p_c^{\ast}=2n/(n-n')$. Hence, $p_c^{\ast}=n/\sigma$
where $\sigma=(n-2)/2-\sqrt{(n-2)^2/4+c}$.
Next let $K_{\alpha}$ and $I_\alpha$ are modified Bessel function,
see \cite[\S 9.6.1 p. 374]{AS} or \cite[\S 1.14 p. 16]{Tra}. Set
$\ell(x)=x^{-n/2+1}I_{n'/2-1}(x)$ and $k(x)=x^{-n/2+1}K_{n'/2-1}(x)$.
Then by (4.2) and Section 6.1 of \cite{HS} the resolvent kernel for $L=\Delta_n+\frac{c}{r^2}$ is given by
\begin{equation}\label{laponker}
R(\lambda)(x,y)=K_{(L+\lambda^2)^{-1}}(x,y) = \left\{ \begin{array}{ll}
   \nu\lambda^{d-2}  k(\lambda y) \ell(\lambda x)   & \mbox{ if $y\ge x$},\\[6pt]
   \nu\lambda^{d-2}  \ell(\lambda y)k(\lambda x)  & \mbox{ if $x >y $}
           \end{array}
    \right.
\end{equation}
for some constant $\nu$.
Next recall that $x^{-\alpha}I_\alpha$ is an even analytic function and that
$$
K_\alpha(x) = {\pi\over 2}\ \frac{ I_{-\alpha} (x) - I_\alpha (x) }{ \sin (\alpha \pi)},
$$
see \cite[\S 9.6.1 p. 374]{AS}. Hence by the limiting absorption principle (that is Stone's Theorem)
for all $x \le y$
\begin{eqnarray*}
K_{dE_{\sqrt{L}}(\lambda)}(x,y)&=&\frac{i\lambda}{\pi}(R(i\lambda)(x,y)-R(-i\lambda)(x,y))\\
&=&
 \nu \frac{i\lambda}{\pi}x^{1-n/2}y^{1-n/2}I_{n'/2-1}(i\lambda x)\frac{I_{{1-n'/2}} (i\lambda y) - I_{n'/2-1} (i\lambda y)}{\sin ((n'/2-1) \pi)}\\
 &&-
 \nu\frac{i\lambda}{\pi} x^{1-n/2}y^{1-n/2}I_{n'/2-1}(-i\lambda x)\frac{I_{{1-n'/2}} (-i\lambda y) - I_{n'/2-1} (-i\lambda y)}{\sin ((n'/2-1) \pi)}.
\end{eqnarray*}
Recall next that  $x^{-\alpha}I_\alpha$ is an even analytic function so
$$
I_{n'/2-1}(i\lambda x)I_{{1-n'/2}} (i\lambda y)=I_{n'/2-1}(-i\lambda x)I_{{1-n'/2}} (-i\lambda y)
$$
and
$$
I_{n'/2-1}(-i\lambda x)I_{{n'/2-1}} (-i\lambda y)=e^{i\pi(n'-2)}I_{n'/2-1}(i\lambda x)I_{{n'/2-1}} (i\lambda y).
$$
Thus
\begin{eqnarray}\label{res}\nonumber
K_{dE_{\sqrt{L}}(\lambda)}(x,y)&=&  \nu \frac{i\lambda}{\pi}
x^{1-n/2}y^{1-n/2}\big( e^{i\pi(n'-2)}- 1\big) \frac{I_{n'/2-1}(i\lambda x)I_{n'/2-1} (i\lambda y)}{\sin ((n'/2-1) \pi)}\\
&=&
{2i \nu\over \pi} e^{i\pi(n'/2-1)}\lambda
x^{1-n/2}y^{1-n/2} I_{n'/2-1}(i\lambda x)I_{n'/2-1} (i\lambda y).
\end{eqnarray}
 We prove equality \eqref{res} under assumption that $x \le y$ but similar argument shows that
\eqref{res} holds for all $x$ and $y$.
Now if we set  $\ell(x)=x^{-n/2+1}I_{n'/2-1}(ix)$ then by \eqref{res}
\begin{equation}\label{eq9.3}
dE_{\sqrt{L}}(\l)f(x)=C \l^{n-1} \ell(\l x)\int_0^\infty \ell(\l y) f(y)y^{n-1}dy
\end{equation}
 for some constant $C$. The function  $\ell = \ell_{n,c}(\lambda)$ for $n > 2$ satisfies the following estimates
(\cite{AS, Tra}):

\begin{equation*}
|\ell(\lambda)|\le  \left\{ \begin{array}{ll}
\lambda^{{n'-n\over 2}}  & \mbox{if} \quad \lambda \le 1\\[6pt]
\lambda^{{1-n\over 2}}&    \mbox{if} \quad 1 \le \lambda.
 \end{array}
    \right.
\end{equation*}

 \noindent  By \eqref{eq9.3}
$$
K_{F(\SL)}(x,y)=C\int_0^\infty F(\l)\ell(\l x) \ell(\l y) \lambda^{n-1}d\l.
$$

\medskip

Let us   prove our estimate (\ref{eq9.2}). We consider only the case $  -{(n-2)^2/4}<c <0$. The proof is similar for the case
$c\ge 0$.   For every $B=B(x_B,
r_B)$, one writes

$$
F(\SL)P_{B}f(x)=C\int_0^\infty\Big(\int_0^\infty   F(\l)\ell(\l x)
\ell(\l y) \lambda^{n-1}d\l \Big)\chi_B(y)f(y)y^{n-1}dy.
$$

\noindent Hence

\begin{eqnarray*}
 \|F(\SL)P_{B}f\|_2^2
&=&C\int_0^\infty\Big|\int_0^\infty\Big(\int_0^\infty
F(\l)\ell(\l x)
\ell(\l y) \lambda^{n-1}d\l\Big)\chi_B(y)f(y)y^{n-1}dy\Big|^2x^{n-1}dx\\
&=&C\int_0^\infty\Big|\int_0^\infty \ell(\l x)\l^{n-1} \int_0^\infty
F(\l)\ell(\l y) \chi_B(y)f(y)y^{n-1}dyd\l\Big|^2x^{n-1}dx.
\end{eqnarray*}

\noindent Note that the following  Plancherel type equality is satisfied
$$
\int_0^\infty\Big|\int_0^\infty
 F(\l)\ell(\l x) \lambda^{n-1} d\l \Big|^2x^{n-1}dx=\int_0^{\infty}|F(\l)|^2 \l^{n-1} d\l,
$$
which  yields
\begin{eqnarray*}
\|F(\SL)P_{B}f\|_2^2&=&C\int_0^\infty|F(\l)|^2\Big|\int_0^\infty
\ell(\l y)\chi_B(y)f(y)y^{n-1}dy\Big|^2\l^{n-1}d\l\\
&\leq&C\int_0^\infty|F(\l)|^2\|\ell(\l y)\chi_B(y)\|_{p'}^2\|f\|_{p}^2\l^{n-1}d\l.
\end{eqnarray*}

\smallskip

\noindent {\bf Case I:} $B=B(x_B,r_B)$ and $x_B\leq 2r_B$. Hence

\begin{eqnarray*}
\big\|\ell (\lambda y)\chi_B\big\|^{p'}_{p'}
&\leq&C\int_0^{2r_B}\big|\ell(\l {y})\big|^{p'}y^{n-1}dy\\
&\leq&C \int_0^{\infty}\big|\ell(\l {y})\big|^{p'}y^{n-1}dy\\
&\leq&C \l^{-n}\int_0^{\infty}\big|\ell(  {y})\big|^{p'}y^{n-1}dy\\
&\leq& C\l^{-n}
\end{eqnarray*}

\noindent for all  $2n/(n-1) < p' < p_c^{\ast}=2n/(n-n')$  (For
$c\geq 0$ this condition should be replaced  by $2n/(n-1) < p' $).
Thus
\begin{eqnarray*}
\|F(\SL)P_{B}f\|_2^2
&\leq&C\int_0^\infty|F(\l)|^2\l^{n-1-\frac{2n}{p'}}d\l\|f\|_{p}^2\\
&\leq&CR^{n-\frac{2n}{p'}}\|\delta_RF\|_2^2\|f\|_{p}^2.
\end{eqnarray*}
\noindent When $x_B\leq 2r_B$, we have that $V(B)\approx r_B^n$ and
$V(B)^{1/2-1/p}r_B^{n(1/p-1/2)}\approx1$. This gives
\begin{eqnarray*}
\|F(\SL)\chi_Bf\|_2^2 &\leq& C V(B)^{2({1\over 2}-{1\over
p})}r_B^{2n({1\over p}-{1\over 2})}R^{n-{2n\over p'}}\|\delta_R
F\|_2^2\|f\|_{p}^2\\
&\leq&C V(B)^{2({1\over 2}-{1\over p})}(Rr_B)^{2n({1\over p}-{1\over
2})}\|\delta_RF\|_2^2\|f\|_{p}^2.
\end{eqnarray*}
This proves Case I.

\smallskip

\noindent {\bf Case II:} $B=B(x_B,r_B)$ and $x_B> 2r_B$. Then
\begin{eqnarray*}
\|\ell(\l y) \chi_B\|^{p'}_{p'}
&\leq&\int_{x_B-r_B}^{x_B+r_B}\big|\ell(\l {y})\big|^{p'}y^{n-1}dy\\
&\leq&C \int_{x_B-r_B}^{x_B+r_B}(\l y)^{\frac{(1-n)p'}{2}}  y^{n-1}dy\\
&\leq&C\l^{p'\frac{1-n}{2}} r_Bx_B^{n-1+p'\frac{1-n}{2}}.
\end{eqnarray*}
Hence
\begin{eqnarray*}
\|F(\SL)P_{B}f\|_2^2
&\leq&C\int_0^\infty|F(\l)|^2r_B^{{2\over p'}}x_B^{2(n-1+\frac{(1-n)p'}{2})/p'  }d\l\|f\|_{p}^2\\
&\leq&Cr_B^{{2\over p'}}x_B^{{2(n-1)\over p'}-n+1}R\|\delta_R
F\|_2^2\|f\|_{p}^2.
\end{eqnarray*}
\noindent Note that if $x_B>2r_B$ then $V(B)\approx x_B^{n-1}r_B$.
 Thus
\begin{eqnarray*}
\|F(\SL)P_{B}f\|_2^2 &\leq& C r_B^{{2\over
p'}}x_B^{{2(n-1)\over p'}-n+1}R\|\delta_R
F\|_2^2\|f\|_{p}^2\\
&\leq& C V(B)^{2({1\over 2}-{1\over p})}(Rr_B)^{2n({1\over
p}-{1\over 2})}(Rr_B)^{1-2n({1\over p}-{1\over 2})}\|\delta_R
F\|_2^2\|f\|_{p}^2\\
&\leq&C V(B)^{2({1\over 2}-{1\over p})}(Rr_B)^{2n({1\over p}-{1\over
2})}\|\delta_RF\|_2^2\|f\|_{p}^2
\end{eqnarray*}

 \noindent
according to the condition  $Rr_B \ge 1$ and $p'>{2n/(n-1)}$. This
proves Case II, and then the proof of Proposition~\ref{prop9.1}
 is complete.
\hfill{}$\Box$

\medskip

\begin{coro}   Suppose that $n > 2$  and $  -{(n-2)^2/4}<c $. For $c<0$, set
 $p\in \big((p_c^{\ast}){'}, {2n\over n+1}\big)$ where $p_c^{\ast}=n/\sigma$ and $\sigma=(n-2)/2-\sqrt{(n-2)^2/4+c}$
 and $(p_c^{\ast}){'}$ its conjugate exponent. For $c\ge 0$, set $p\in \big[1, {2n\over n+1}\big).$
Then
\begin{itemize}
\item[(i)]
   $S_R^{\delta(p)} (L)$
is of weak-type $(p,p)$ uniformly in $R$.

\item[(ii)] For any bounded Borel
function $F: [0, \infty) \to \CC$ such that
$\sup_{t>0}\|\eta\, \delta_tF\|_{W^{\beta, 2}}<\infty $ for some  $\beta>\max\{n(1/p-1/2),1/2\}$
 and some non-trivial function $\eta \in C_c^\infty(0,\infty)$, the operator
$F(\SL)$ is bounded on $L^r(X)$ for all $p<r<p'$ with
\begin{eqnarray*}
   \|F(\SL)  \|_{r\to r}\leq    C_\beta\Big(\sup_{t>0}\|\eta\, \delta_tF\|_{W^{\beta, 2}}
   + |F(0)|\Big).
\end{eqnarray*}
\end{itemize}
\end{coro}

\medskip

\noindent {\bf Proof.} This result is straightforward from
Proposition~\ref{prop9.1} and Theorems~\ref{th4.1}, ~\ref{th31.1}.
  \hfill{}$\Box$

  \medskip

\begin{remark}\label{re9.3}
 Note that   for the standard Laplacian $\Delta$ on
 ${\mathbb R}^n$    Stein-Tomas estimate  ${\rm (R_p)}$ holds  if and only if
 $1\leq p \leq {2(n+1)/(n+3)}$. Surprisingly, if $n > 2$  and $  -{(n-2)^2/4}<c<0,$ then the restriction estimate
 ${\rm (ST^{2}_{p, 2})}$ for the operator $\Delta_n+\frac{c}{r^2}$ holds only
 for all $p\in ((p_c^{\ast}){'}, 2n/(n+1))$ where
 $p_c^{\ast}=n/\sigma$ and $\sigma=(n-2)/2-\sqrt{(n-2)^2/4+c}$.
 \end{remark}

\section{Examples of $(1,2)$-restriction type conditions}\label{sec55}
\setcounter{equation}{0}

In \cite{DOS} the following {\it Plancharel   condition} is introduced:
  for any  $R>0$ and all even
   Borel functions $F$ such that $\supp F\subseteq [0,R]$,
 \begin{equation}\label{eq11.1}
    \int_X  |K_{F(\sqrt{L})} (x,y)|^2d\mu(x) \le C V(y,R^{-1})^{-1}
     \Vert \delta_RF \Vert_{q}^2
     \end{equation}
   for some $q \in [2,\infty]$.   Here $m\geq 2$ denotes the order of operator $L$ and for operators which
   we consider here $m=2$(see \cite[(3.1)]{DOS}). Note that for every $x\in X$ and $r \ge 1/R,$
\begin{eqnarray}
 \big\| F(\SL) P_{B(x,r)}f\big\|_{2}&\leq& \big\|\int_X K_{F(\SL)}(z,w)
 \chi_{B(x,\, r)}(w) f(w)d\mu(w)\big\|_{2}\nonumber\\
 &\leq&\int_X \big\|K_{F(\SL)}(\cdot,\, w)\big\|_{2}\chi_{B(x,\, r)}(w) |f(w)|d\mu(w)
 \nonumber\\
 &\leq&C   \Vert \delta_RF \Vert_{q} \int_X V(w,R^{-1})^{-1/2}\chi_{B(x, \, r)}(w) |f(w)|d\mu(w)
 \nonumber\\
 &\leq&C  V(x,r)^{-1/2} (Rr)^{n/2}  \Vert \delta_RF \Vert_{q}  \|f\|_1\nonumber,
 \end{eqnarray}
 where in the last inequality we used the doubling condition (\ref{eq2.2}).
  Therefore,  for every $x\in X$ and $r \ge1/R$
 \begin{equation}\label{eq11.2}
   \Vert {F(\sqrt{L})}P_{B(x,r)}\Vert_{1 \to 2}  \le C V(x,r)^{-1/2} (Rr)^{n/2}
     \Vert \delta_RF \Vert_{q},
     \end{equation}
and so  the condition  (\ref{eq11.1}) is just a  slightly stronger
version of condition ${\rm (ST^q_{1,2})}$.

 It was noted in \cite{DOS} that condition  \eqref{eq11.1} with $q=2$ holds for homogeneous
 sub-Laplacian acting on homogeneous
 Lie groups. It is shown in \cite{S1} that  \eqref{eq11.1} with $q=2$ holds also
 for ``quasi-homogeneous" subelliptic and elliptic operators.
As we  note  condition \eqref{eq11.1} is stronger than
condition ${\rm (ST^2_{1,2})}$ so this implies the following result.

\begin{proposition}\label{prop11.1}
Let $L$ be a homogeneous sub-Laplacian or  ``quasi-homogeneous" operator acting on homogeneous
 Lie group with  homogeneous dimension $d$. Then the Riesz mean
 $S_R^{(d-1)/2}(L)$ of order $(d-1)/2$ is of weak-type $(1,1)$ uniformly in $R$.
\end{proposition}

\noindent
{\bf Proof.}
 Proposition~\ref{prop11.1} follows directly from Theorem~\ref{th4.1}.
  \hfill{}$\Box$

\medskip

We believe that in this generality Proposition~\ref{prop11.1} is a
new result. However in the case of Heisenberg group it follows from
the result obtained  by M\"uller, Stein and Hebisch that the the
Riesz means of order $\delta> {(d_e-1)/2}$ is bounded on $L^1$,
where $d_e<d$ is the topological dimension of the Heisenberg(see
\cite{MS, H}). Therefore it is likely that
Proposition~\ref{prop11.1} is not a genuine endpoint result.

\bigskip

\noindent
{\bf Acknowledgements:} A. Sikora was partly  supported by
Australian Research Council  Discovery Grant DP 110102488. L. Yan
was  supported by  NNSF of China (Grant No.  10925106),
 Guangdong Province Key Laboratory of Computational Science
 and the Fundamental Research Funds for the Central Universities (Grant No. 09lgzs610).
  Part of this work was done while
 E.M. Ouhabaz was visiting Macquarie University. His visit was partly supported
 by ARC Discovery Grant DP 110102488 and CNRS. He wishes
 to thank X.T. Duong for the invitation. A. Sikora and E.M. Ouhabaz would like to thank
 M. Cowling, X.T. Duong, A. Hassell  and A. McIntosh for fruitful discussions.

\end{document}